\documentclass[twoside,11pt]{article}

\usepackage{jmlr2e}
\usepackage{amsmath,color,colordvi,caption,subcaption,float,url}
\usepackage{hyperref}
\graphicspath{ {./Figures/} }
\usepackage{footmisc}

\usepackage{multirow}

\newcommand{\be}{\begin{equation}}
\newcommand{\ee}{\end{equation}}
\newcommand{\ba}{\begin{eqnarray}}
\newcommand{\ea}{\end{eqnarray}}
\newcommand{\bas}{\begin{eqnarray*}}
\newcommand{\eas}{\end{eqnarray*}}

\def\ba{\mathbf a}
\def\bb{\mathbf b}
\def\be{\mathbf{e}}
\def\bd{\mathbf d}

\def\bw{\mathbf{w}}
\def\bx{\mathbf x}
\def\by{\mathbf y}
\def\bz{\mathbf z}

\def\G{{\cal G}}

\def\N{{\cal N}}
\def\P{{\cal P}}
\def\T{{\cal T}}
\def\R{\mathbb{R}}
\def\CoSaMP{\texttt{CoSaMP}}
\def\HTP{\texttt{HTP}}
\def\FHTP{\texttt{FHTP}}
\def\GraSP{\texttt{GraSP}}
\def\GraHTP{\texttt{GraHTP}}
\def\IHT{\texttt{IHT}}
\def\AIHT{\texttt{AIHT}}
\def\NHTP{\texttt{NHTP}}
\def\SIP{\texttt{SIP}}
\def\support{\mbox{supp}}
\def\cs{\mbox{cs}}
\def\dist{\mbox{dist}}
\def\supp{\rm{supp}}

\firstpageno{1}

\begin{document}

\title{Global and Quadratic Convergence of Newton Hard-Thresholding Pursuit\thanks{~This research was partially supported by the National Natural Science Foundation of China (11971052, 11926348), ``111'' Project of China (B16002), and The Alan Turing
Institute.}}

\author{\name Shenglong Zhou \email shenglong.zhou@soton.ac.uk \\
       \addr School of Mathematical Sciences\\
       University of Southampton\\
       Southampton SO17 1BJ,  UK
       \AND
       \name Naihua Xiu \email nhxiu@bjtu.edu.cn \\
       \addr Department of Applied Mathematics \\
             Beijing Jiaotong University \\
             Beijing, China
       \AND
       \name Hou-Duo Qi \email hdqi@soton.ac.uk \\
       \addr School of Mathematical Sciences\\
       and CORMSIS\\
       University of Southampton\\
       Southampton SO17 1BJ,  UK}


\maketitle

\begin{abstract}
Algorithms based on the hard thresholding principle have been well studied with sounding theoretical guarantees in the compressed sensing and more general sparsity-constrained optimization. It is widely observed in existing empirical studies that when a restricted Newton step was used (as the debiasing step),  the hard-thresholding algorithms tend to meet halting conditions in a significantly low number of iterations and are very efficient.
{{
Hence, the thus obtained Newton hard-thresholding algorithms call for stronger
theoretical guarantees than for their simple hard-thresholding counterparts.
}} 
This paper provides a theoretical justification for the use of the restricted Newton step.
We build our theory and algorithm, Newton Hard-Thresholding Pursuit (\NHTP), for the sparsity-constrained optimization. Our main result shows that \NHTP\ is quadratically convergent under the standard assumption of restricted strong convexity and smoothness.
We also establish its global convergence to a stationary point under a weaker assumption.
In the special case of the compressive sensing, \NHTP\ effectively reduces to some of the existing hard-thresholding algorithms with a Newton step.  Consequently, our fast convergence result justifies why those algorithms perform better than without the Newton step. The efficiency of \NHTP\ was demonstrated on both synthetic and real data in compressed sensing and sparse logistic regression.
\end{abstract}


{\bf Running Title:} Newton Hard-Thresholding Pursuit
\section{Introduction}

In this paper, we are mainly concerned with numerical methods for
the sparsity constrained optimization
\begin{equation}\label{SCO}
\min_{\bx\in \mathbb{R}^n }\ \  f(\bx) , \qquad {\rm s.t.} \qquad \|\bx\|_0\leq s,
\end{equation}
where $f: \mathbb{R}^n \mapsto \mathbb{R}$ is continuously differentiable,
$\|\bx\|_0$ is the $l_0$ norm of $\bx$, counting the number of nonzero elements in $\bx$,
and $s$ is a given integer regulating the sparsity level in $\bx$ (i.e., $\bx$ is $s$-sparse).
This problem has been well investigated by \cite{bahmani2013greedy} (from statistical learning perspective)
and \cite{Beck13} (from optimization perspective).
Problem (\ref{SCO}) includes the widely studied Compressive Sensing (CS) (see, e.g., \cite{Elad2010, Zhao2018}) as a special
case:
\begin{equation} \label{CS}
  \min_{\bx\in \mathbb{R}^n } \ \ f(\bx) =  f_{\cs}(\bx) = \frac 12 \| A \bx - \bb \|^2 , \qquad {\rm s.t.} \qquad \|\bx\|_0\leq s,
\end{equation}
where $A$ is an $m \times n$ sensing matrix, $\bb\in \mathbb{R}^n$ is the observation and
$\| \cdot\|$ is the Euclidean norm in $\mathbb{R}^n$.
Problem (\ref{SCO}) has also been a major model in high-dimensional statistical recovery
(\cite{agarwal2010fast}, \cite{negahban2012unified}), 
nonlinear compressive sensing (\cite{blumensath2013compressed}), 
and learning model-based sparsity (\cite{BBR16}).
An important class of algorithms makes use of the gradient information together with the hard-thresholding technique.
We refer to \cite{bahmani2013greedy}, \cite{YLZ18} (for (\ref{SCO})) and
\cite{needell2009cosamp}, \cite{foucart2011hard} (for CS)
for excellent examples of such methods 
and their corresponding theoretical results.
In terms of the numerical performance, it has been widely observed that whenever a restricted
Newton step is used in the so-called debiasing step, those algorithms appear
to take a significantly low number of iterations to converge,
{{see \cite{foucart2011hard} and Remark 8 in Section~\ref{Subsect-CS}}}.
Yet, their theoretical guarantee appears no better than their pure gradient-based counterparts.
Hence, there exists an intriguing gap between the exceptional empirical experience and the best
convergence theory. This paper aims to provide a theoretical justification for their
efficiency by establishing the quadratic convergence of such methods under standard
assumptions used in the literature.

In the following, we give a selective review of past work that directly motivated our research,
followed by a brief explanation of our general framework that shares a similar structure with several existing algorithms.

\subsection{A selective review of past work} \label{Subsect-Review}

There exists a large number of computational algorithms that can be applied to (\ref{SCO}).
For instance, many of them can be found in Google Scholar from the many papers citing
\cite{Wright2007}, \cite{needell2009cosamp}, \cite{Elad2010} and also in the latest book by \cite{Zhao2018}. We opt to conduct a bit technical review on a small number of papers that
directly motivated our research.
Those reviewed papers more or less suggest the following algorithmic framework that largely obeys
the principle laid out in \cite{needell2009cosamp} and follow the recipes for hard-thresholding
methods in \cite{KC2011}.
Given the $k$th iterate $\bx^k$, update it to the next iterate $\bx^{k+1}$ by the following
steps:
\begin{equation} \label{General-Framework}
\left\{  \begin{array}{ll}
\mbox{Step 1 (Support Identification Process)}: & \ T_k = \texttt{SIP} (h (\bx^k)), \\ [0.6ex]
\mbox{Step 2 (Debiasing)}: & \  \widetilde{\bx}^{k+1} = \arg\min \ \left\{
  q_k(\bx): \ \ \bx |_{T_k^c} = 0 \right\} ,\\ [0.6ex]
\mbox{Step 3 (Pruning)}: & \ \bx^{k+1} \in \P_s( \widetilde{\bx}^{k+1} ).
\end{array}
\right.
\end{equation}
We put the three steps in the perspectives of some existing algorithms and explain the notation
involved.
For the case of CS,
the well-known \CoSaMP\ (Compressive Sample Matching Pursuit) of \cite{needell2009cosamp} chose
the identification function $h(\bx)$ to be the gradient function $\nabla f(\bx)$
and the support identification process \SIP\ is chosen to be the union of the best $2s$ support of $h(\bx^k)$ (i.e., the $2s$ indices that are from the $2s$ largest elements of $h(\bx^k)$ in magnitude)
 and $\support(\bx^k)$, which are the indices of nonzero elements in $\bx^k$.
In this case, the number of indices in $T_k$ is below $3s$ (i.e., $|T_k | \le 3s$).
In the \HTP\ (Hard Thresholding Pursuit) algorithm of \cite{foucart2011hard},
$h(\bx)$ is set to be $(\bx - \eta \nabla f(\bx))$, where $\eta >0$ is a step size.
$T_k$ is chosen to be the best $s$ support of $h(\bx^k)$. Hence, $|T_k| = s$.
In \AIHT\ of \cite{Blu12} (Accelerated Iterative Hard Thresholding), $T_k$ is chosen as
in \HTP.
For the general nonlinear function $f(\bx)$,
the \GraSP\ of \cite{bahmani2013greedy} (Gradient Support Pursuit)
chose $T_k$ as in \CoSaMP\ so that $T_k| \le 3s$.
The \GraHTP\ of \cite{YLZ18} (Gradient Hard Thresholding Pursuit) chose $T_k$ as in \HTP\ for
CS.

Once $T_k$ is chosen, Step 2 (debiasing step) attempts to provide a better estimate for
the solution of (\ref{SCO}) by solving an optimization problem within a restricted subspace
obtained by setting all elements of $\bx$ indexed by $T_k^c$ to zero. Here $T_k^c$ is the
complementary set of $T$ in $\{1, \ldots, n\}$. Step 3 (pruning step) simply applies the
hard-thresholding operator, denoted as $\P_s$, to $\widetilde{\bx}^{k+1}$.
{{To be more precise, $\P_s(\bx)$ contains all minimal $2$-norm distance solutions from $\bx$ under the
$s$-sparsity constraint: 
\[ 
 \P_s(\bx) ={\rm argmin}_{\bz} \left\{  \|\bx-\bz\|\ | \ \|\bz\|_0\leq s \right\},
\] 
}}
which can be obtained by retaining the $s$ largest elements in magnitude from $\bx$ and  
setting the remaining to zero. 
The great flexibility in choosing $T_k$ and the objective function $q_k(\bx)$ in the debiasing
step makes it possible to derive various algorithms in literature.
For instance, if we choose $T_k = \{1, \ldots, n \}$ (hence $T_k^c = \emptyset$) and
$q_k(\bx)$ to be the first-order approximation of $f$ with a proximal term at $\bx^k$:
\[
  q_k(\bx) := f(\bx^k) + \langle \nabla f(\bx^k), \ \bx - \bx^k \rangle + \frac 1{2\eta}
  \| \bx - \bx^k \|^2,
\]
then we will recover the popular (gradient) hard-thresholding algorithms,
see, e.g., \cite{blumensath2008gradient, blumensath2009iterative} and \cite{Beck13}
for the iterated hard-thresholding algorithms,
and \cite{bahmani2013greedy} for the restricted gradient descent
 and \cite{YLZ18} for \GraHTP.
The \CoSaMP\ is recovered if $T_k$ is chosen as in \CoSaMP\ and $q_k(\bx) = f_{\cs}(\bx)$.
More existing methods can be interpreted this way and we omit the details here.

Instead, we focus on the algorithms that make use of the second order approximation in $q_k(\bx)$. \cite{bahmani2013greedy} proposed the restricted Newton step, which is equivalent
to choosing $q_k(\bx)$ to be a restricted second-order approximation to $f(\bx)$ at $\bx^k$:
\begin{equation} \label{Restriced-QP}
q_k(\bx) := f(\bx^k) + \langle \nabla f(\bx^k_{T_k}), \ \bx_{T_k} - \bx^k_{T_k} \rangle
 + \frac 1{2} \langle \bx_{T_k} - \bx^k_{T_k}, \ \nabla^2_{T_k} f(\bx^k) (\bx_{T_k} - \bx^k_{T_k})
 \rangle
\end{equation}
where the notation $\bx_{T_k}$ denotes the restriction of $\bx$ to the indices in $T_k$,
$\nabla f(\bx^k_{T_k})$ is the (partial) gradient of $f(\bx)$ with respect to the variables indexed by $T_k$ and evaluated at $\bx^k_{T_k}$,  and
$\nabla^2_{T_k} f(\bx^k) $ is the principle submatrix of the Hessian matrix $\nabla^2 f(\bx^k)$
indexed by $T_k$. In the case of CS (\ref{CS}), the restricted Newton step is equivalent to
minimizing $q_k(\bx) = f_{\cs}(\bx)$ restricted on the subspace defined by $\bx|_{T_k^c} = 0$. Hence,
the restricted Newton step recovers \CoSaMP. We note that in both cases, $|T_k | \le 3s$ (i.e., $T_k$ is relatively large). In the \HTP\ algorithm, \cite{foucart2011hard} managed to choose
$T_k$ of size $s$ by making use of the hard thresholding technique, which is further investigated
by \cite{Blu12} by the name of accelerated iterative hard-thresholding.

The benefit of using the Newton step has been particularly witnessed for the case of CS. \cite{foucart2011hard} compiled convincing numerical evidence that \HTP\ took a significantly low number of iterations to converge when proper step-size $\eta$ is used.
However, the existing theoretical guarantee for \HTP\
is no better than their greedy counterparts (e.g., simple iterative hard-thresholding algorithms (\IHT)). That is, the theory ensures that the distance between each iterate to any given reference (sparse) point is bounded by the sum of two terms. The first term converges linearly
and the second term is a fixed approximation error that depends on the choice of the reference point. We refer to the latest paper of \cite{SL18} for many of such a result, which is often called statistical error-bound guarantee.
The discrepancy between being able to offer better empirical performance than many simple \IHT\ algorithms and only sharing similar theoretical guarantee with them invites an intriguing question: why is it so? A positive answer will inevitably provide a deep understanding of the
Newton-type HTP algorithms and lead to new powerful algorithms. This is exactly what we are going to achieve in this paper.

A different line of research for (\ref{SCO}) was initiated by \cite{Beck13} from an optimization
perspective. The convergence results established were drastically contrasting to the statistical
error bound result mentioned above.
It is proved that any accumulation point of the generated sequence by the \IHT\ method is
one kind of stationary point (i.e., $\eta$-stationarity, to be defined later).
In the particular case of CS, the whole sequence converges to an {{ $\eta$-stationary point}}
under the $s$-regularity assumption of the sensing matrix $A$
({{ i.e.,  any $s$ columns of $A$ are linearly independent}}).
It is known that $2s$-regularity is a minimal condition that any two $s$-sparse vectors can be
distinguished and 
{{
it is often assumed by many quantities
}}
related to the restricted isometry property (RIP)
of \cite{candes2005decoding}.
The fact that the $s$-regularity is weaker than the $2s$-regularity means that many hard-thresholding algorithms actually converge to an $\eta$ stationary point of (\ref{SCO}).
Hence, the quality of those algorithms can be measured not only by their statistical error bounds,
but also by the quality of the $\eta$ stationary point (e.g., whether a stationary point is
optimal).  We refer to \cite{Beck13, beck2015minimization} for more discussion on the $\eta$
stationarity in relation to the global optimality.

Similar convergence results to \cite{Beck13} have also been established in the literature of CS.
\cite{blumensath2010normalized} showed that the normalized \IHT\ with an adaptive step-size rule converges to a local minimum of (\ref{CS}) provided that the $s$-regularity holds.
This leads us to ask the following question: 
{{when the Newton step is used in a framework of \IHT\ (such as \HTP\ algorithm of \cite{foucart2011hard}), we would like to know whether the resulting algorithm enjoys
the following fast quadratic convergence:
}}
\begin{equation} \label{QConvergence}
\bx^k \ \ \rightarrow \ \ \bx^* \qquad \mbox{and} \qquad
\| \bx^{k+1} - \bx^* \| \le c \| \bx^k - \bx^* \|^2 \quad \mbox{for sufficiently large}\ k,
\end{equation}
where $c$ is a constant solely dependent on the objective function $f$ (independent of the iterates $\bx^k$ and its limit $\bx^*$).
This fast convergence result would justify the stronger numerical performance of various Newton-type methods reviewed in the first part of the subsection.
Although, it is expected in optimization that Newton's method (\cite{NW1999}) will usually lead to quadratic convergence, the problem (\ref{SCO}) is not a standard optimization problem and
it has a combinatorial nature. Hence, quadratic convergence does not follow from any existing theory from optimization.
We also note that both \cite{bahmani2013greedy} and \cite{YLZ18} listed the restricted Newton
step as a possible variant for the debiasing step, but it was not theoretically investigated.

We finish this brief review by noticing that there are researches that exclusively studied the role of Newton's method for (\ref{SCO}) (see, e.g., \cite{dai2009subspace}, \cite{yuan2017newton}, \cite{chen2017fast}).
However, as before, they did not offer any better theoretical guarantees than their simple greedy counterparts. Furthermore, their algorithms do not follow the general framework of (\ref{General-Framework}) and hence
their results cannot be used to explain the efficiency of Newton's method that follows (\ref{General-Framework}).
In this paper, we will design an algorithm, that also makes uses of a restricted Newton step in the debiasing step (Step 2) and analyse its role in convergence.
We will show that our algorithm enjoys the quadratic convergence (\ref{QConvergence}) as well as others.
We will particularly relate it to
\HTP\ of \cite{foucart2011hard} so as to justify the strong empirical performance of similar
algorithms.

\subsection{Our approach and main contributions}

The first departure of our proposed Newton step from the one of \cite{bahmani2013greedy} is that
we employ a different quadratic function, denoted as $q_k^N(\bx)$:
\begin{eqnarray} \label{qkN}
q_k^N(\bx) &:=& \mbox{the second-order Taylor expansion of $f(\bx)$ at $\bx^k$, then set} \ \bx|_{T_k^c} = 0 \\
&=& \langle \nabla_{T_k} f(\bx^k), \ \bx_{T_k} - \bx^k_{T_k} \rangle
+ \frac 1{2} \langle \bx_{T_k} - \bx^k_{T_k}, \ \nabla^2_{T_k} f(\bx^k) (\bx_{T_k} - \bx^k_{T_k}) \rangle  \nonumber \\
&& -\ \langle \bx_{T_k}, \ \nabla^2_{T_k, T_k^c} f(\bx^k) (\bx^k_{T_k^c}) \rangle
   + \ (\mbox{constant term independent of $\bx$}), \nonumber
\end{eqnarray}
where $\nabla_{T_k} f(\bx^k) := ( \nabla f(\bx^k))_{T_k}$ and $\nabla^2_{T_k, T_k^c} f(\bx^k)$ is
the submatrix whose rows and columns are from the Hessian matrix $\nabla^2 f(\bx^k)$
indexed by $T_k$ and $T_k^c$ respectively.
For the case of CS problem (\ref{CS}),
it is straightforward to verify that $q_k^N(\bx) = q_k(\bx)$ in (\ref{Restriced-QP}).
Therefore, the Newton step will become the one used in \CoSaMP\ or \HTP\ depending on how $T_k$ is
selected. In this paper, we choose $T_k$ to be the best $s$ support of $\bx^k - \eta \nabla f(\bx^k)$.
That is, $T_k$ contains a set of indices that define the $s$ largest absolute values in $\bx^k - \eta \nabla f(\bx^k)$ with $\eta$ being steplength.
For the case of CS, it is the same as that in the algorithm $\HTP^{\mu}$ of
\cite{foucart2011hard}.
For the general nonlinear function $f$, however, $q_k^N(\bx)$ and $q_k(\bx)$ are different.
The function $q_k(\bx)$ in (\ref{Restriced-QP})  is obtained in such a way that we first
restrict $f(\bx)$ to the subspace $\bx|_{T_k^c} = 0$ and then approximate it by the second-order Taylor expansion (i.e., restriction and approximation).
In contrast, the function $q_k^N(\bx)$ is obtained in the opposite way.
We first approximate $f(\bx)$ by its second-order Taylor expansion and then restrict the approximation to the subspace $\bx_{T_k^c} = 0$ (i.e., approximation and restriction).
We will see that our way of construction will allow us quantitatively bound the error $\| \bx^{k+1} - \bx^* \|$
in terms of $\| \bx^k - \bx^* \|$, eventually leading to the quadratic convergence in (\ref{QConvergence}).

Our second innovation is to cast the Newton step as a Newton iteration for a nonlinear equation:
\begin{equation} \label{Feta}
F_{\eta} (\bx, \ T_k) = 0,
\end{equation}
where $F_{\eta}(\cdot, T_k): \mathbb{R}^n \mapsto \mathbb{R}^n$ is a function reformulated from
the $\eta$ stationarity condition. We defer its technical definition to the next section.
A crucial point we would like to make is that this new interpretation of the Newton step
offers a fresh angle to examine it and will allow
us to develop new analytical tools mainly from optimization perspective and eventually
establish the promised quadratic convergence.

It is known that Newton's method is a local method. A commonly used technique for globalization is
the line search strategy, which is adopted in this paper.
Therefore, we will have a Newton iterate with varying step-size.
This agrees with the empirical observation that adaptive step-size in \HTP\ often works more efficiently
than other variants. Putting together the three techniques (quadratic approximation $q_k^N(\bx)$,
nonlinear equation (\ref{Feta}), and the line search strategy) will result in our proposed algorithm
termed as Newton Hard-Thresholding Pursuit (\NHTP) due to the Newton-step and the way how $T_k$
is selected  
being the two major components of the algorithm.
We finish this section by summarizing our major contributions.

\begin{itemize} 
	\item[(i)] We develop the new algorithm \NHTP, which largely follows the general framework of
	(\ref{General-Framework}) with Step 3 (pruning) to be replaced by a globalization step.
	The new step is achieved through the Armijo line search.
	We will establish its global convergence to an $\eta$ stationary point
	under the restricted strong smoothness of $f$.

	\item[(ii)] If $f$ is further assumed to be restricted strongly convex at the one of the accumulation points	of \NHTP, the Armijo line search steplength will eventually becomes $1$.
	Consequently, \NHTP\ will become the restricted Newton method and leads to its convergence at a quadratic rate. This result successfully extends the
	classical quadratic convergence result of Newton's method to the sparse case.
	For the case of CS, \NHTP\ reduces to some known algorithms including the \HTP\ family of \cite{foucart2011hard}, with properly chosen step-sizes. The quadratic convergence result
	resolves the discrepancy between the strong numerical performance of \HTP\ (and its alike) and its existing linear convergence guarantee.
	
	\item[(iii)] Rigorously establishing the quadratic convergence of \NHTP\ is a major contribution
	of the paper.
	As far as we know, it is the first paper that establishes both the global and
	the quadratic convergence for an algorithm  that employs both the Newton step and the gradient
	step (through the hard thresholding operator) for (\ref{SCO}).
	The developed framework of analysis is innovative and will open possibility to prove that
	other Newton-type HTP methods may also enjoy the quadratic convergence.
	In our final contribution of this paper, we show experimental results in CS and the 
	logistic regression, with both synthetic and real data, to illustrate the way \NHTP\ works.
	
\end{itemize}

\subsection{Organization}

In the next section, we will describe the basic assumptions on the objective function $f$
and their implications. We will also develop a theoretical foundation for the Newton method to be used
in a way that it also solves a system of nonlinear equations.
Section \ref{Section-NHTP} includes the detailed description of \NHTP\ and its global and quadratic
convergence analysis.
We will particularly discuss its implication to the CS problem and compare with the
methods of \HTP\ family \cite{foucart2011hard}.
Since some of the proofs are quite technical, 
we move all of the proofs to the appendices in order to avoid interrupting presentation of
the main results.
We report our numerical experiments 
 in Section \ref{Section-Numerical} and
conclude the paper in Section \ref{Section-Conclusion}.

\section{Assumptions, Stationarity and Interpretation of Newton's Step} \label{Section-Assumptions}

\subsection{Notation}

For easy reference, we list some commonly used notation below.

\begin{table}[h!]
	\begin{tabular}{ll}
$:=$ & means ``define'' \\
$\bx$ & a column vector and hence $\bx^\top$ is a row vector. \\
$x_i$ & the $i$th element of a vector $\bx$. \\
$x_{(i)}$ & the $i$th largest absolute value among the elements of $\bx$. \\
$\supp(\bx)$ & the support set of $\bx$, namely, 
the set of indices of nonzero elements of $\bx$. \\
$T$ &  index set from $\{1, 2, \ldots, n\}$ \\
$|T|$ &  the number of elements in $T$ (i.e., cardinality of $T$). \\
$ T^c$ & the complementary set of $T$ in $\{1,2,\ldots,\}\setminus T$. \\
$\bx_T$ & the sub vector  of $\bx$  containing elements  indexed on $T$. \\
$\nabla_T f(\bx)$ & $=(\nabla f(\bx))_T$. \\
$\nabla^2 f(\bx)$ & the Hessian matrix of function $f(\cdot)$ at $\bx$. \\
$\nabla^2_{T,J} f(\bx)$ & the submatrix of the Hessian matrix whose rows and columns \\
                         & are respectively indexed by $T$ and $J$. \\
$\nabla^2_{T} f(\bx)$ & $ =\nabla^2_{T,T} f(\bx)$. \\
$\nabla^2_{T:} f(\bx)$ & the submatrix of the Hessian matrix containing rows indexed by $T$.\\
$\langle \bx, \by \rangle$ & the standard inner product for $\bx, \by \in \Re^n$.\\
$\| \bx\|$ & the norm induced by the standard inner product (i.e., Euclidean norm). \\
$\|\bx\|_\infty$ & $= \max \{ |x_i|\}$ (the infinity norm of $\bx \in \Re^n$). \\
$\|A\|_2$ &        the spectral norm of the matrix $A$. \\
$\|A\|$   &        may refer to any norm of $A$ equivalent to $\|A\|_2$. \\
\end{tabular}
\end{table}
$\P_s(\bx)$ has been defined in Section~\ref{Subsect-Review}. 
It is important to  note that $\P_s(\bx)$ may have multiple {{best $s$-sparse approximations}}.
For example, for $\bx^\top =(1,2,-1,0)$ and $s=2$, $\P_s(\bx)$ contains two {{best $s$-sparse approximations}}:
$(1,2,0,0)$ and $(0,2,-1,0)$. 



\subsection{Basic assumptions and stationarity}

In order to study the convergence of various algorithms for the problem (\ref{SCO}), some kind of regularities needs to be assumed. They are more or less analogous to the RIP for CS (see \cite{candes2005decoding}). Those regularities often share the property of strong restricted convexity/smoothness, see \cite{agarwal2010fast}, \cite{SSZ2010}, \cite{jalali2011learning}, \cite{negahban2012unified}, \cite{bahmani2013greedy}, \cite{blumensath2013compressed}, and \cite{YLZ18}.
We state the assumptions below in a way that is conducive to our technical proofs.

\begin{definition} \label{def-RSCS}(Restricted strongly convex and smooth functions)
Suppose that $f: \mathbb{R}^n \mapsto \mathbb{R}$ is a twice continuously differentiable function
whose Hessian is denoted by
$\nabla^2 f(\cdot)$. Define
\[
   M_{2s}(\bx) := \sup_{\by \in \mathbb{R}^n}\left\{
    \langle \by, \;  \nabla^2 f(\bx) \by \rangle \ \bigg| \
     | \supp(\bx) \cup \supp(\by) | \le 2s, \ \| \by \| = 1
   \right\}
\]
and
\[
m_{2s}(\bx) := \inf_{\by \in \mathbb{R}^n}\left\{
\langle \by, \;  \nabla^2 f(\bx) \by \rangle \ \bigg| \
| \supp(\bx) \cup \supp(\by) | \le 2s, \ \| \by \| = 1
\right\}
\]
for all $s$-sparse vectors $\bx$.
\begin{itemize}
	\item[(i)] We say $f$ is restricted strongly smooth (RSS) if there exists a constant $M_{2s} >0$
	such that $M_{2s}(\bx) \le M_{2s}$ for all $s$-sparse vectors $\bx$. In this case, we say $f$ is $M_{2s}$-RSS. $f$ is said to be locally RSS at $\bx$ if $M_{2s}(\bz) \le M_{2s}$ only holds for
	those $s$-sparse vectors $\bz$ in a neighborhood of $\bx$.
	
	\item[(ii)] We say $f$ is restricted strongly convex (RSC) if there exists a constant $m_{2s} >0$
	such that $m_{2s}(\bx) \geq m_{2s}$ for all $s$-sparse vectors $\bx$. In this case, we say $f$ is
	$m_{2s}$-RSC. $f$ is said to be locally RSC at $\bx$ if $m_{2s}(\bz) \ge m_{2s}$ only holds for those
	$s$-sparse vectors $\bz$ in a neighborhood of $\bx$.
	
\item[(iii)] We say that $f$ is locally restricted Hessian Lipschitz continuous at an $s$-sparse vector $\bx$ if
there exists a Lipschitz constant  $L_f$ and a neighborhood {{$\N_s(\bx):=\{\bz\in\R^n:\supp(\bx)\subseteq\supp(\bz),\|\bz\|_0\leq s\}$ such that
\[
  \| \nabla^2_{T:} f(\by) - \nabla^2_{T:} f(\bz)  \| \le L_f \| \by - \bz\|, \ \ \ \ \forall~ \by,\bz\in \N_s(\bx), 
\]
for any index set $T$ with $|T|\leq s$ and $T \supseteq\supp(\bx)$.}}
\end{itemize}

\end{definition}

\noindent{\bf Remark 1.}
We note that the definition of $M_{2s}(\bx)$ and $m_{2s}(\bx)$ is taken from the definition of the restricted stable Hessian (RSH) of \cite[Def.~1]{bahmani2013greedy}. If $m_{2s}(\bx)$ is bounded away
from zero, the RSH is equivalent to the RSC and RSS putting together. Under the assumption of
twice differentiability, RSS and RSC become that of \cite{negahban2009unified}, \cite{SSZ2010} and \cite[Def.~1]{YLZ18}.
{{
The local condition (iii) is a technical condition required for proving the quadratic convergence of our algorithm.
Typical examples of such function satisfying (iii) include the quadratic function (\ref{CS}) and the quartic 
function studied in \cite{Beck13}:
\[
  f(\bx) = \sum_{i=1}^\ell  \left(
   \bx^\top A_i \bx - c_i
  \right)^2,
\]
where $A_i, i=1, \ldots, \ell$ are $n \times n$ symmetric matrices and $c_i$, $i=1, \ldots, \ell$ are given.
}}
By a standard calculus argument, $M_{2s}$-RSS implies
\begin{equation} \label{Eq-Lipschitz}
\left\{
 \begin{array}{l}
   \| \nabla f(\bx) - \nabla f(\by) \| \le M_{2s} \| \bx - \by \|, \\ [0.6ex]
   f(\bx) - f(\by) - \langle \nabla f(\bx),\; \bx - \by \rangle \le \frac{M_{2s}}{2} \| \bx - \by \|^2 ,
 \end{array} \ \
 \begin{array}{l}
 \forall\ \bx, \by,  \ | \support(\bx) | \le s \\
 | \support(\bx) \cup \support(\by) | \le 2s.
 \end{array}
\right .
\end{equation}
The properties in (\ref{Eq-Lipschitz}) ensure that any optimal solution of (\ref{SCO}) must be an
$\eta$-stationary point, which is a major concept introduced to {{the sparse optimization}} (\ref{SCO}) by
\cite{Beck13}. We state the concept below.

{{
\begin{definition} ($\eta$-stationarity) \cite[Def.~2.3]{Beck13} 
An $s$-sparse vector $\bx^*$ is called an $\eta$-stationary point of (\ref{SCO}) if it satisfies the
following relation
\[
 \bx^* \in \P_s( \bx^* - \eta \nabla f(\bx^*)).
\]
\end{definition} 
}}
\cite{Beck13} called it the $L$-stationary point because $\eta$ is very much related to the Lipschitz
constant $M_{2s}$ defined in (\ref{Eq-Lipschitz}). Lemma 2.2 in \citep{Beck13} states that an $s$-sparse vector $\bx^*$ is an $\eta$-stationary point if and only if
\begin{eqnarray}\label{agradient}
\nabla_\Gamma f(\bx^*)=0,\ \ \ \
\|\nabla_{\Gamma^c} f(\bx^*)\|_\infty\leq  x^*_{(s)}/\eta.
\end{eqnarray}
where $\Gamma :=\supp(\bx^*)$. By invoking the proofs of \cite[Lemma~2.4 and Thm.~2.2]{Beck13}
under the condition of (\ref{Eq-Lipschitz}), the existence of $\eta$-stationary point is ensured.

\begin{theorem} \label{Thm-Beck}
(Existence of $\eta$-stationary point) \cite[Thm.~2.2]{Beck13}.
Suppose that there exists a constant $M_{2s} >0$ such that (\ref{Eq-Lipschitz}) holds.
Let $\eta < 1/M_{2s}$ and $\bx^*$ be an optimal solution of (\ref{SCO}). Then
\begin{itemize}
	\item[(i)] $\bx^*$ is an $\eta$-stationary point;
	\item[(ii)] $\P_s(\bx^* - \eta \nabla f(\bx^*))$ contains exactly one element.
\end{itemize}	
Consequently, we have
\begin{equation} \label{Fix-Point-Eq}
  \bx^* = \P_s(\bx^* - \eta \nabla f(\bx^*)) .
\end{equation}
\end{theorem}

We would like to make a few remarks on the significance of Thm.~\ref{Thm-Beck}.\\

\noindent {\bf Remark 2.} The characterization of the optimal solution $\bx^*$ as a solution of
the fixed-point equation (\ref{Fix-Point-Eq}) immediately suggests a simple iterative procedure:
\[
  \bx^{k+1} \in \P_s ( \bx^k - \eta \nabla f(\bx^k)), \ \ \ k=0,1, 2, \ldots .
\]
{{Indeed, for the special case of CS, we have $\nabla^2 f(\bx) = A^\top A$ and the fact
		about the relationship between the spectral norm $\|A\|_2$ and the quantity $M_{2s}$:
$$\| A\|_2^2=\sup_{\by\in\mathbb{R}^n, \|\by\|=1} \langle \by, A^\top A \by\rangle \geq \sup_{\|\by\|_0\leq 2s, \|\by\|=1} \langle \by, A^\top A \by\rangle  = M_{2s}.$$}}
\noindent
When $\| A \|_2 < 1$, the unit length choice of $\eta = 1$, which satisfies 
$ 1 < 1/\| A\|_2^2 \le 1/M_{2s}$, recovers the \IHT\ of \cite{blumensath2008gradient}.
Moreover, any stationary point of
$\{ \bx^k\}$ is an $\eta$-stationary point and satisfies the fixed-point equation (\ref{Fix-Point-Eq}).
For the case $\| A \|_2 \ge 1$, the same conclusion holds as long as $\eta < 1/M_{2s}$,
see \cite[Remark~3.2]{Beck13}.\\

\noindent{\bf Remark 3.}
The fixed-point equation characterization also measures how far an $s$-sparse point $\bx$ is from being
an $\eta$-stationary point (and hence a possible candidate for an optimal solution of (\ref{SCO})) by
computing
\begin{equation} \label{Halting-Condition}
 h(\bx, \eta) := \dist( \bx, \; \P_s( \bx - \eta \nabla f(\bx))),
\end{equation}
which defines the shortest Euclidean distance from $\bx$ to the set $\P_s( \bx - \eta \nabla f(\bx))$.
If $h(\bx, \eta)$ is below a certain tolerance level (e.g., small enough), we may stop at $\bx$.
This halting criterion is different from those commonly used in CS literature such as in \CoSaMP, \GraSP, and \HTP.

Our next remark is about a differentiable nonlinear equation reformulation of the fixed-point equation (\ref{Fix-Point-Eq}) and it will give rise to a nice interpretation of the Newton step obtained from
minimizing $q_k^N(\bx)$ in (\ref{qkN}).  This remark is the main content of the next subsection.

\subsection{Nonlinear equations and new interpretation of Newton's step}

Given a point $\bx\in \mathbb{R}^n$ and $\eta>0$, we define the collection of all index sets of best $s$-support of the vector {{$\bx - \eta \nabla f(\bx)$ by
\begin{equation} \label{Tu}
 \T(\bx; \eta) := \left\{
 T \subset \{1, \ldots, n \} \ \bigg| \ 
  \begin{array}{l}
 |T| = s, ~T\supseteq \supp(\bz),\\
 \exists \ \bz\in\P_s( \bx - \eta \nabla f(\bx))
 \end{array}
 \right\}.
\end{equation}
That is, each $T$ in $\T(\bx; \eta)$ includes $s$ indices that define the locations of the $s$ largest absolute values among the elements of $\bx - \eta \nabla f(\bx)$. 
}}
Then for any given $T \in \T(\bx; \eta)$, we define the corresponding nonlinear equation:
\begin{eqnarray}\label{station}
F_{\eta}(\bx;T):=
\left[\begin{array}{c}
\nabla_{ T}  f (\bx)\\
\bx_{ T^c }\\
\end{array}\right]=
0.
\end{eqnarray}
One advantage of defining the function $F_\eta(\bx; T)$ is that it is continuously differentiable with respect to $\bx$ once $T$ is selected. Moreover, we have the following characterization of the fixed-point equation
(\ref{Fix-Point-Eq}) in terms of $F_{\eta}$.

\begin{lemma} \label{Lemma-Feta}
{{Suppose $\eta>0$ is given.  A  point $\bx \in \mathbb{R}^n$  is an $\eta$-stationary point if and only if 
$$
  F_{\eta}(\bx;T) = 0,  \qquad  \exists \ T\in \T (\bx;\eta).
$$
Furthermore, }}a   point $\bx \in \mathbb{R}^n$ satisfies the fixed point equation
(\ref{Fix-Point-Eq}) if and only if	
$$
  F_{\eta}(\bx;T) = 0,  \qquad  \forall \ T\in \T (\bx;\eta).
$$
\end{lemma}

\noindent{{\bf Remark 4 (Deriving a new stopping criterion)}}
This result is instrumental and crucial to our algorithmic design.
Bearing in mind that it is impossible to solve all the nonlinear equations associated with all
possible $T \in \T(\bx; \eta)$ to {{ get a solution  that satisfies the fixed point equation}}, 
our hope is that solving one such equation would lead to our desired results.
To monitor how accurately the equation (\ref{station}) is solved, we develop a new
stopping criterion that involves the gradient of $f$ in both parts indexed by $T$ and $T^c$.

{{
  \begin{itemize}
  	\item[(i)] We note that $\bx_{T^c} = 0$ in (\ref{station}) is easily satisfied. Hence, 
  	the magnitude $\|F_{\eta}(\bx;T)\|$ of the residual actually measures the gradient of $f$ on the $T$ part.
  	
  	\item[(ii)] Now suppose $\bx$ satisfies (\ref{station}). If follows from the definition of $T$ that
  	\[
  	| x_j | =  | x_j - \eta \nabla_j f(\bx) | \ge | x_i - \eta \nabla_i f(\bx) | = \eta | \nabla_i f(\bx) | , 
  	\qquad \forall \ j \in T \ \ \mbox{and} \ \ \forall \ i \in T^c .
  	\]
  	This, together with $\bx_{T^c} = 0$ and $|T| = s$, leads to
  	\[
  	  x_{(s)} = \min_{j \in T} | x_j | \ge \eta | \nabla_i f(\bx) | 
  	\]
  	or equivalently
  	\begin{equation} \label{Extra-Condition}
  	|\nabla_i f(\bx) | \le \frac 1{\eta} x_{(s)}, \qquad \forall \ i \in T^c
  	\end{equation}
  	This is the gradient condition on the $T^c$ part that an $\eta$-stationary point has to satisfy.
  	Therefore, a measure on the violation of this condition indicates
  	how close it is approximated on the $T^c$ part.
  \end{itemize} 
}}


\noindent
Consequently, a natural tolerance function {{ to measure how far $\bx$ is from being an $\eta$-stationary point}} is
\begin{equation} \label{Tolerance-Function}
\mbox{Tol}_{\eta}(\bx;\; T) := \| F_\eta (\bx; T) \| + \max_{i \in T^c} \left\{ \max \Big( | \nabla_i f(\bx) | - x_{(s)}/\eta, \ 0  \Big) \right\}.
\end{equation}
It is easy to see that the halting function $h(\bx, \eta) = 0$ in (\ref{Halting-Condition}) implies
that there exists $T \in \T(\bx; \eta)$ such that $\mbox{Tol}(\bx; \; T) =0$ {{and vice versa. }}
Our purpose is to quickly find this correct $T$.

We now turn our attention to the solution methods for (\ref{station}).
Suppose $\bx^k$ is the current approximation to a solution of (\ref{station}) and $T_k$ is chosen
from $\T(\bx^k; \eta)$. 
Then Newton's method for the nonlinear equation (\ref{Feta}) takes the following form to get the next iterate
$\widetilde{\bx}^{k+1}$:
\begin{equation} \label{Newton-Method}
F'_{\eta} (\bx^k; T_k) (\widetilde{\bx}^{k+1} - \bx^k) = - F_{\eta} (\bx^k; T),
\end{equation}
where $F'_{\eta} (\bx^k; T_k)$ is the Jacobian of $F_{\eta}(\bx; T_k)$ at $\bx^k$ and it assumes the following form:
\begin{equation} \label{Jacobian-Matrix}
F'_{\eta} (\bx^k; T_k) = \left[
 \begin{array}{cc}
  \nabla^2_{T_k} f(\bx^k) \ & \ \nabla^2_{T_k, T_k^c} f(\bx^k)  \\ [0.6ex]
  0 & I_{n-s}
 \end{array}
\right] .
\end{equation}
Let $\bd_N^k := \widetilde{\bx}^{k+1} - \bx^k$ be the Newton direction. Substituting (\ref{Jacobian-Matrix}) into
(\ref{Newton-Method}) yields
\begin{equation}\label{Newton-Direction}
\left\{
\begin{array}{rcl}
\nabla_{T_k}^{2} f (\bx^k)  (\bd_N^k)_{T_{k}} &=&
   \nabla_{   T_{k}, T_{k}^c }^{2} f (\bx^k) \bx^{k}_{T_{k}^c} - \nabla_{ T_{k}}  f (\bx^k) \\ [1.5ex]
(\bd_N^{k})_{ T^c_{k}}&=&-\bx^{k}_{{ T }_{k}^c}.
\end{array}
\right.
\end{equation}
At this point, it is interesting to observe that the next iterate $(\widetilde{\bx}^{k+1} = \bx^k + \bd_N^k)$  is exactly the one we would get for the restricted Newton step from minimizing the restricted
quadratic function $q_k^N(\bx)$ in (\ref{qkN}).
It is because of this exact interpretation of the restricted Newton step that it also drives the equation
(\ref{station}) to be eventually satisfied. In this way, we establish the global convergence to
the $\eta$-stationarity.
However, there are still a number of technical hurdles to overcome.
We will tackle those difficulties in the next section.

\section{Newton Hard-Thresholding Pursuit and Its Convergence} \label{Section-NHTP}

In this main section, we present our Newton Hard-Thresholding Pursuit (\NHTP) algorithm,
which largely follows the general framework (\ref{General-Framework}), but with distinctive features.
We already discussed the choice of $T_k$ (Step 1 in (\ref{General-Framework})) and the quadratic
approximation function $q_k^N$ in (\ref{qkN}) (Step 2 in (\ref{General-Framework})).
Since $|T_k | = s$ and $\widetilde{\bx}^{k+1}$ obtained is restricted to the subspace $\bx|_{T_k^c} = 0$,
hence $\support( \widetilde{\bx}^{k+1} ) \subseteq T_k$ and the pruning step is not necessary.
Instead, we replace it with the globalization step:
\begin{equation} \label{Globalization-Step}
\mbox{Step 3' (globalization)} \quad
\left\{
 \begin{array}{l}
   \bx^{k+1} = \G( \widetilde{\bx}^{k+1} ) \ \mbox{such that} \\ [0.6ex]
   \support(  \bx^{k+1} ) \subseteq T_k \ \ \mbox{and} \ \ f(\bx^{k+1}) \le f( \bx^k),
 \end{array}
\right .
\end{equation}
where $\G$ symbolically represents a globalization process to generate $\bx^{k+1}$.
{{
We emphasize that globalization here refers to a process that will generate a sequence of iterates
from any initial point and the sequence converges to an $\eta$-stationary point.
The descent condition in (\ref{Globalization-Step}) will be realized by
the Armijo line search strategy (see \cite{NW1999}).
We also emphasize, however, that there are other strategies that may work for globalization.
}}

The rest of the section is to consolidate those three steps.
We first examine how good is the restricted Newton direction (\ref{Newton-Direction}) as well as
the restricted gradient direction. We note that both directions were proposed in \cite{bahmani2013greedy}.
But as far as we know, they are not theoretically studied. We then describe our \NHTP\ algorithm and present its global and quadratic convergence under the restricted strong convexity and smoothness.

\subsection{Descent properties of the restricted Newton and gradient directions}

Our first task is to answer whether the restricted Newton direction $\bd_N^k$ from (\ref{Newton-Direction}) provides
a ``good'' descent direction for $f(\bx)$ on the restricted subspace $\bx|_{T_k^c} =0$.
We have the following result.

\begin{lemma} \label{Lemma-Newton-Direction}
	(Descent inequality of the Newton direction)
	Suppose $f(\bx)$ is $m_{2s}$-restricted strongly convex and $M_{2s}$-restricted strongly smooth.
	Given a constant $\gamma \le m_{2s}$ and the step-size $\eta \le 1/(4M_{2s})$, we then have
	\begin{equation} \label{Descent-Inequality-Newton}
	 \left\langle \nabla_{T_k} f(\bx^k), \ ( \bd^k_N)_{T_k} \right\rangle
	 \le - \gamma \| \bd^k_N \|^2 + \frac{1}{4\eta} \| \bx^k_{T_k^c} \|^2.
	\end{equation}
\end{lemma}

We note that $T_k$ will eventually identify the true support
and $\bx^k_{T_k^c}$ should be close to zero when
this happens. Hence, the positive term $\| \bx^k_{T_k^c} \|^2/(4\eta)$ is eventually negligible and
the restricted Newton direction is able to provide a reasonably good descent
direction on the subspace $\bx_{T_k^c}  =0$.
But in general (e.g., $f(\bx)$ is not restricted strongly convex), the inequality (\ref{Descent-Inequality-Newton}) may not hold and hence $\bd^k_N$ may not provide a good descent direction at all.
In this case, we opt for the restricted gradient direction (denoted by $\bd^k_g$ to distinguish it from
$\bd^k_N$):
\begin{equation} \label{gradient-Direction}
  \bd_g^k := \left[
  \begin{array}{c}
   (\bd^k_g)_{T_k} \\ [1ex]
   (\bd^k_g)_{T_k^c}
  \end{array}
  \right]  = \left[
  \begin{array}{c}
  - \nabla_{T_k} f(\bx^k) \\ [0.6ex]
  - \bx_{T_k^c}
  \end{array}
  \right].
\end{equation}
This strategy of switching to the gradient direction whenever the Newton direction is not good enough
(by certain measure) 
appears very popular and practical in optimization,
see, e.g., \cite{NW1999,SWQ2002, QQS2003, QS2006, ZST2010}.
Therefore, our search direction $\bd^k$ for {{ the globalization step (Step 3' (\ref{Globalization-Step}))}} is defined as follows:
\begin{equation} \label{Search-Direction}
 \bd^k := \left\{
 \begin{array}{ll}
  \bd^k_N , \quad &  \mbox{if the condition  (\ref{Descent-Inequality-Newton}) is satisfied} \\ [0.6ex]
  \bd_g^k , \quad & \mbox{otherwise}.
 \end{array}
 \right .
\end{equation}
It is important to note that the choice of $\gamma$ and $\eta$ in Lemma~\ref{Lemma-Newton-Direction} {{is 
sufficient but not necessary}} for the Newton direction to be used. The inequality (\ref{Descent-Inequality-Newton}) may also hold if $\gamma$ and $\eta$ violate the required bounds. This has been experienced in our numerical
experiments.

Our next result further shows that the search direction $\bd^k$ is actually a descent direction for $f(\bx)$ at $\bx^k$ with respect to the full space $\mathbb{R}^n$ provided that $\eta$ is properly chosen.
Suppose we have three constants $\gamma$, $\sigma$ and $\beta$ such that
\begin{equation} \label{Choice-of-Constants}
 0 < \gamma \le \min\{1,\; 2M_{2s} \}, \ \ 0 < \sigma < 1/2, \ \ \mbox{and} \ \  0< \beta < 1.
\end{equation}
They will be used in our \NHTP\ algorithm. We note that this choice implies
$M_{2s}/\gamma > \sigma$.
Define two more constants based on them:
\begin{equation} \label{alpha-eta}
 \overline{\alpha} := \min \left\{ \frac{1-2\sigma}{ M_{2s}/\gamma - \sigma}, \ 1   \right\} \ \ \
 \mbox{and}\ \ \
 \overline{\eta} := \min\left\{ \frac{\gamma (\overline{\alpha } \beta ) }{ M_{2s}^2 },
 \  \overline{\alpha } \beta, \ \frac{1}{4M_{2s}}\right\}.
\end{equation}

\begin{lemma} \label{Lemma-dk}
	(Descent property of $\bd^k$)
Suppose $f(\bx)$ is $M_{2s}$-restricted strongly smooth.
Let $\gamma, \sigma$ and $\beta$ be chosen as in (\ref{Choice-of-Constants}).
Suppose $\eta < \overline{\eta}$
 and
$\support(\bx^{k}) \subseteq T_{k-1}$ (this will be automatically ensured by our algorithm).
 We then have
\begin{equation} \label{Descent-dk}
 \langle \nabla f(\bx^k), \; \bd^k \rangle
 \le -\rho \| \bd^k \|^2 - \frac{\eta}{2} \| \nabla_{T_{k-1}} f(\bx^k) \|^2,
\end{equation}
where $\rho >0$ is given by
\[
  \rho := \min\left\{
   \frac{2\gamma - \eta M^2_{2s} }{ 2 }, \ \frac{2 - \eta}{2}
  \right\} .
\]
\end{lemma}
Lemma~\ref{Lemma-dk} will ensure that our algorithm \NHTP\ is well defined.

\subsection{\NHTP\ and its convergence}

Having settled that $\bd^k$ is a descent direction of $f(\bx)$ at $\bx^k$, we compute the next
iterate along the direction $\bd^k$ but restricted to the subspace $\bx|_{T_k} = 0$:
$\bx^{k+1} = \bx^k(\alpha_k)$ with $\alpha_k$ being calculated through the Armijo line search  and
\begin{equation} \label{xk-alpha}
   \bx^k(\alpha) := \left[
   \begin{array}{c}
    \bx^k_{T_k} + \alpha \bd^k_{T_k} \\ [0.6ex]
    \bx^k_{T_k^c} + \bd^k_{T_k^c}
   \end{array}
   \right]
   = \left[
   \begin{array}{c}
   \bx^k_{T_k} + \alpha \bd^k_{T_k} \\ [0.6ex]
   0
   \end{array}
   \right]  , \qquad \alpha > 0.
\end{equation}
Our algorithm is described in Table~\ref{Alg-NHTP}.

\begin{table}[H]
	\caption{Framework of \NHTP \label{Alg-NHTP}}\vspace{-6 mm}
	{\renewcommand\baselinestretch{1.25}\selectfont
		{\centering\begin{tabular}{ ll }\\ \hline
				\multicolumn{2}{l}{\NHTP: Newton Hard-Thresholding Pursuit}\\\hline
				\textbf{Step 0}&Initialize $\bx^0$. Choose $\eta, \gamma>0, \sigma\in(0,1/2), \beta\in(0,1)$. Set $k\Leftarrow0$.   \\
				\textbf{Step 1}& Choose $T_{k} \in \T(\bx^k; \eta)$.\\
				\textbf{Step 2}& If ${\tt Tol}_{\eta}(\bx^k; T_k) =0$, then stop.
				Otherwise,  go to \textbf{Step 3}.\\
				\textbf{Step 3}& Compute the search direction $\bd^k$ by (\ref{Search-Direction}). \\
				\textbf{Step 4}& Find the smallest integer $\ell = 0,1,\ldots$ such that\\
				& \parbox{0.75\textwidth}{
					\begin{equation}\label{Armijo}
					f(\bx^{k}(\beta^{\ell}))\leq f(\bx^k)+\sigma\beta^{\ell} \langle \nabla  f(\bx^{k}),  \bd^{k} \rangle.
					\end{equation}}\\
				& Set $\alpha_k=\beta^{\ell}$, $\bx^{k+1}=\bx^{k}( \alpha_k)$ and $k\Leftarrow k+1$, go to {\bf Step 1}.
				\\\hline
			\end{tabular}\par} }
\end{table}

\noindent
{\bf Remark 5.} 
We will see that 
\NHTP\ has a fast computational performance because of two factors. 
One is that it terminates in a low number of iterations due to the quadratic convergence (to be proved)
and this has been experienced in our numerical experiments. 
The other is the low computational complexity of each step. 
For example,  for both CS  and sparse logistic regression problems, the computational complexity of each step is $\mathcal{O}(s^3+ms^2+mn+ms\ell)$, where $\ell$ is the smallest integer satisfying (\ref{Armijo})
and it often assumes the value $1$.
{{The way  $\bx^k(\alpha)$ is defined guarantees that }}$\support(\bx^{k+1}) \subseteq T_k$ for
all $k=0,1, \ldots,$. 
If ${\tt Tol}_{\eta}(\bx^k; T_k) =0$, then $\bx^k$ is already an $\eta$-stationary point and we should 
terminate the algorithm. Without loss of any generality, we assume that \NHTP\ generates an infinite sequence $\{ \bx^k\}$ and we will analyse its convergence properties.
The line search condition (\ref{Armijo}) is known as the Armijo line search and
ensures a sufficient decrease from $f(\bx^k)$ to $f(\bx^{k+1})$.
Therefore, the two properties in the globalization step (\ref{Globalization-Step})
is guaranteed, provided that the line
search in (\ref{Armijo}) is successful. This is the main claim of the following result.

\begin{lemma} \label{Lemma-alphak}
	(Existence and boundedness of $\alpha_k$)
Suppose$f(\bx)$ is $M_{2s}$-restricted strongly smooth.
Let the parameters $\gamma$, $\sigma$ and $\beta$ satisfy the conditions in (\ref{Choice-of-Constants})
and $\overline{\alpha}$ and $\overline{\eta}$ be defined in (\ref{alpha-eta}).
Suppose ${\tt Tol}_{\eta}(\bx^k; T_k) \not=0$.
For any $\alpha$ and $\eta$ satisfying
\[
  0 < \alpha \le \overline{\alpha}
  \qquad \mbox{and} \qquad
  0 < \eta < \min \left\{
   \frac{\alpha \gamma}{M^2_{2s}}, \ \alpha, \ \frac{1}{4M_{2s}}
  \right\},
\]
it holds
 \begin{eqnarray}\label{alpha-decreasing-property}
  f(\bx^k(\alpha)) \le f(\bx^k) + \sigma \alpha \langle \nabla f(\bx^k), \; \bd^k \rangle.
 \end{eqnarray}
Consequently, if we further assume that $\eta \le \overline{\eta}$, we have
\[
  \alpha_k \ge \beta \overline{\alpha} \qquad \forall \ k=0,1, \ldots, .
\]
\end{lemma}

It is worth noting that the objective function is only assumed to be restricted strongly smooth (not necessarily to be restricted strongly convex).
{{ Lemma~\ref{Lemma-alphak}  not only ensures}} the existence of $\alpha_k$ that satisfies the
line search condition (\ref{Armijo}), but also guarantees that $\alpha_k$ is always bounded away from zero
by a positive margin $\beta \overline{\alpha}$. This boundedness property will in turn ensure that \NHTP\
will converge. 
Our first result on convergence is about a few quantities approaching zero.

\begin{lemma} \label{Lemma-Converging-Quantities}
(Converging quantities)
{{ Suppose $f(\bx)$ is $M_{2s}$-restricted}} strongly smooth.
Let the parameters $\gamma$, $\sigma$ and $\beta$ satisfy the conditions in (\ref{Choice-of-Constants})
and $\overline{\eta}$ be defined in (\ref{alpha-eta}).
We further assume that $\eta \le \overline{\eta}$. Then the following hold.

\begin{itemize}
	
	\item[(i)] $\{ f(\bx^k)\}$ is a nonincreasing sequence and if $\bx^{k+1} \not = \bx^k$, then
	 $f(\bx^{k+1}) < f(\bx^k)$.
	
	\item[(ii)] $\| \bx^{k+1} - \bx^k \| \rightarrow 0$;
	
	\item[(iii)] $ \| F_\eta(\bx^k;\; T_k) \| \rightarrow 0 $;
	
	\item[(iv)] $ \| \nabla_{T_k} f(\bx^k) \| \rightarrow 0 $ and
	             $ \| \nabla_{T_{k-1} } f(\bx^k) \| \rightarrow 0 $.
	
	
\end{itemize}

\end{lemma}

Those converging quantities are the basis for our main results below. They also justify the halting conditions that we will use in our numerical experiments.

\begin{theorem} \label{Thm-Global-Convergence}
(Global convergence)	
Suppose$f(\bx)$ is $M_{2s}$-restricted strongly smooth.
Let the parameters $\gamma$, $\sigma$ and $\beta$ satisfy the conditions in (\ref{Choice-of-Constants})
and $\overline{\eta}$ be defined in (\ref{alpha-eta}).
We further assume that $\eta \le \overline{\eta}$. Then the following hold.

\begin{itemize}
	\item[(i)] Any accumulation point, say $\bx^*$, of the sequence $\{ \bx^k \}$ is an $\eta$-stationary point of (\ref{SCO}).
	If $f$ is a convex function, then for any given reference point $\bx$ we have
	\begin{equation} \label{epsilon-Optimality}
	  f(\bx^*) \le f(\bx) + \frac{ x^*_{(s)} }{ \eta } \| \bx_{\Gamma_*^c} \|_1,
	\end{equation}
	where $\Gamma_* := \support(\bx^*)$.
	
	\item[(ii)] If $\bx^*$ is isolated, then the whole sequence converges to $\bx^*$. Moreover, we have
	the following characterization on the support of $\bx^*$.
	
	 \begin{itemize}
	 	\item[(a)] If $\| \bx^* \|_0 = s$, then
	 	\[
	 	   \support(\bx^*) = \support(\bx^k) = T_k \quad \mbox{for all sufficiently large } k.
	 	\]
	 	
	 	\item[(b)] If $\| \bx^* \|_0 < s$, then
	 	\[
	 	\support(\bx^*) \subseteq \support(\bx^k) \cap T_k \quad \mbox{for all sufficiently large } k.
	 	\]
	 	
	 \end{itemize}

\end{itemize}

\end{theorem}

\noindent
{\bf Remark 6.}
Under the assumption of $f$ being restricted strongly smooth, \NHTP\ shares the most desirable
convergence property (i.e., to $\eta$-stationary point) of the iterative hard-thresholding algorithm
of \cite{Beck13}. If $f$ is assumed to be convex, then (\ref{epsilon-Optimality}) implies that
for any given $\epsilon>0$, there exists 
neighborhood $\N(\bx^*)$ of $\bx^*$ such that $f(\bx^*) \le f(\bx) + \epsilon$ 
for any $\bx \in \N(\bx^*)$. In particular, if $\| \bx^*\|_0 = s$, then $\bx^*$ is a local minimum of
(\ref{SCO}). If $\| \bx^*\|_0 < s$ (so that $x^*_{(s)} = 0$), then $\bx^*$ is a global optimum of
(\ref{SCO}).

It achieves more. If the generated sequence converges to $\bx^*$, the support
of $\bx^*$ is eventually identified as $T_k$ provided that the sparse level of $\bx^*$ is $s$.
If $\| \bx^* \|_0 <s$, its support would be eventually included in $T_k$.
When specialized to the CS problem (\ref{CS}) with $s$-regularity, the whole sequence $\{\bx^k\}$ will
convergence to one point $\bx^*$. This is because that any $\eta$-stationary point of the CS problem under
the $s$-regularity is isolated, see \cite[Lemma~2.1 and Corollary 2.1]{Beck13}.
Our next result implies that under the $2s$-regularity, the whole sequence $\{\bx^k\}$ converges to $\bx^*$ at a quadratic rate.

\begin{theorem} \label{Thm-Qudratic-Convergence}
(Quadratic convergence)
Suppose all conditions as in Thm.~\ref{Thm-Global-Convergence} hold. Let $\bx^*$ be one of the
accumulation points of $\{\bx^k\}$.
We further assume  $f(\bx)$ is  $m_{2s}$-restricted strongly convex in a neighborhood of $\bx^*$.
If $\gamma \le \min\{ 1, \ m_{2s} \}$ and $\eta \le \overline{ \eta}$,
 then the following hold.

\begin{itemize}
	\item[(i)] The whole sequence $\{\bx^k\}$ converges to $\bx^*$, which is necessarily an $\eta$-stationary point.
	
	\item[(ii)] The Newton direction is   accepted for sufficiently large $k$. 

	\item[(iii)] If we further assume that
	$f$ is locally restricted Hessian Lipschitz continuous at $\bx^*$ with the Lipschitz constant $L_f$. The line search steplength eventually becomes unity  and the convergence rate of $\{\bx^k \}$ to $\bx^*$ is quadratic.
	That is, there exists an iteration index $k_0$ such that
	\begin{equation} \label{Quadratic-Convergence-in-x}
	  \alpha_k\equiv1,\ \ \ \ \| \bx^{k+1} - \bx^* \| \le \frac{  L_f }{ 2m_{2s} } \| \bx^k - \bx^*\|^2,
	  \qquad \forall \ k \ge k_0.
	\end{equation}
	Moreover, for sufficiently large $k$, we have
	\[
	  \| F_\eta ( \bx^{k+1}; \; T_{k+1}) \| 
	  \le \frac{L_f \sqrt{M_{2s}^2+1}}{\min\{m^3_{2s},m_{2s}\}}  \| F_\eta ( \bx^{k}; \; T_{k}) \|^2.
	\]
	
\end{itemize}
		
\end{theorem}

\noindent{\bf Remark 7.} Taking into account of Lemma~\ref{Lemma-Converging-Quantities}(iii) that
$\| F_\eta(\bx^k; \; T_k) \|$ converges to $0$, Thm.~\ref{Thm-Qudratic-Convergence}(iii) asserts that
it converges at a quadratic rate. 
Compared with the quadratic convergence (\ref{Quadratic-Convergence-in-x}), the quadratic convergence
in $\| F_\eta(\bx^k; \; T_k) \|$ has the advantage that it is computationally verifiable.
The quantity is also a major part of our stopping criterion in monitoring $ \mbox{Tol}_{\eta}(\bx^k;\; T_k)$
of (\ref{Tolerance-Function}), see Sect.~\ref{Section-Numerical}.
{{ In addition, we proved the existence of $k_0$. The proof of   Thm.~\ref{Thm-Qudratic-Convergence}(iii) suggests that $k_0$ should be near to the iteration when the support sets of the sequence start  to be identified to be
		the correct support set at its limit. 
However,  deriving an explicit form of $k_0$ is somehow difficult and would require extra conditions.
  }}

\subsection{The case of CS} \label{Subsect-CS}

We use this part to demonstrate the application and implication of our main convergence results to
the CS problem (\ref{CS}). We will also discuss the similarities to and differences from the
existing algorithms, in particular the \HTP\ family of \cite{foucart2011hard}.
The purpose is to show that there is a wide range of choices for the parameters that will lead to
quadratic convergence.
This is best done in terms of the restricted isometry constant (RSC) of the sensing matrix $A$.
We recall from \cite{candes2005decoding} that RSC $\delta_s$ is the smallest $\delta\ge 0$ such that
\[
  (1- \delta) \| \bx\|^2 \le \| A \bx \|^2 \le (1+\delta) \| \bx\|^2 \qquad \forall \ \| \bx\|_0 \le s.
\]
We will use $\delta_{2s}$, which is assumed to be positive throughout. For this setting, we have
\[
   m_{2s} = 1 - \delta_{2s}, \quad M_{2s} = 1+ \delta_{2s}, \quad \mbox{and} \quad
   \mu_{2s} := \frac{M_{2s}}{m_{2s}}>1,
\]
where, $\mu_{2s}$ is known as the $2s$-restricted stable Hessian coefficient of $A$ in \cite{bahmani2013greedy}. For simplicity, we choose a particular set of parameters used in our \NHTP\ to illustrate our results
(many other choices are also possible).{{ Let
\[  
 \beta = \frac 14, \quad \gamma = m_{2s}, \quad \sigma =  \frac{1-w}{2-w/\mu_{2s}} \ \ \mbox{with} \ \ 0 < w< 1.
\] 
It is easy to see that $\sigma \in (0, 1/2)$ for any choice $w$ between
$0$ and $1$. This set of parameter choices certainly satisfies the condition (\ref{Choice-of-Constants}).
We now calculate $\overline{\alpha}$ and $\overline{\eta}$ defined in (\ref{alpha-eta}).
The definition of $\overline{\alpha}$ chooses
\[
  \overline{ \alpha } \overset{(\ref{alpha-eta})}{=} \frac{1-2\sigma}{M_{2s} / \gamma - \sigma} =  \frac{1-2\sigma}{\mu_{2s}- \sigma} = \frac{w}{\mu_{2s} }\in (0,1).
\]
Since $\gamma/M_{2s}^2 = 1/(\mu_{2s} M_{2s}) < 1$ and
$\beta = 1/4$, we have
\begin{eqnarray*}
  \overline{ \eta}
  &\overset{(\ref{alpha-eta})}{=}&  \frac{\gamma}{M_{2s}^2} \overline{\alpha} \beta
  = \frac 14 \times \frac{1}{\mu_{2s}} \times \frac{1}{M_{2s}}  \times \overline{\alpha} \\
  &\ge& \frac 14 \times \frac{1}{\mu_{2s}} \times \frac{1}{2} \times \frac{w}{ \mu_{2s} }
         \qquad (\mbox{because} \ M_{2s} \le 2) \\
  &=& \frac{w}{ 8 \mu^2_{2s}  }.
\end{eqnarray*}}}
Direct application of Thm.~\ref{Thm-Qudratic-Convergence} yields the following corollary.

\begin{corollary} \label{Cor-CS}
Suppose the RIC $\delta_{2s} > 0$ and the parameters of \NHTP\ are chosen as follows:
\begin{equation} \label{CS-Parameters}
 \beta = \frac 14, \ \ {{ \sigma =\frac{1-w}{2-w/\mu_{2s}}}}, \ \ \gamma = m_{2s},
 \ \ \eta \le \frac{w}{ 8 \mu^2_{2s} } \ \ \mbox{with} \ \ 0< w < 1.
\end{equation}
Then \NHTP\ is well-defined. In particular, the Newton direction $d^k_N$ is always accepted as the search direction in (\ref{Search-Direction}) at
each iteration. Moreover, \NHTP\ enjoys all the three convergence results in Thm.~\ref{Thm-Qudratic-Convergence}.
\end{corollary}

\noindent
{\bf Remark 8.} (On RIC conditions)
In the literature of CS, a benchmark condition (for theoretical investigation) often takes the form
$\delta_{t} \le \delta_*$ with $t$ being an integer.
Suppose $\delta_{2s} \le \delta_*$. It is easy to define and derive the following.
\begin{eqnarray*}
 m_{2s}^* &:=& 1 - \delta_* \le 1- \delta_{2s} = m_{2s} \\
 \mu_{2s}^* &:=& \frac{1+ \delta_*}{1-\delta_*} \ge \frac{ 1+ \delta_{2s} }{ 1 - \delta_{2s} } = \mu_{2s} \\
 \overline{ \eta }^* &:=& \frac{w}{ 8 (\mu^*_{2s})^2 } \le \frac{w}{ 8 \mu^2_{2s} }.
\end{eqnarray*}
Therefore, in the selection of the parameters in (\ref{CS-Parameters}), $\mu_{2s}$ and $m_{2s}$ can be
respectively replaced by $\mu^*_{2s}$ and $m^*_{2s}$, and $\eta$ can be chosen to satisfy
$\eta \le \overline{ \eta }^*$. In the scenario of \cite{garg2009gradient} where $\delta_* = 1/3$, with $w = 0.5$ we could choose the parameters as $\beta = 1/4$, $\gamma = 2/3$, $\sigma = 2/7$ and $\eta = 1/64$. This set of
choices would ensure \NHTP\ converges quadratically under the RIP condition $\delta_{2s} \le \delta_*= 1/3$.\\

\noindent{\bf Remark 9.} (On Newton's direction)
That the Newton direction is always accepted at each iteration is because the inequality (\ref{Descent-Inequality-Newton}) is always satisfied with the parameter selection in (\ref{CS-Parameters}) (its proof can be patterned after that for Thm.~\ref{Thm-Qudratic-Convergence}(iii)).
Therefore, the Newton direction $\bd^k_N$ at each iteration takes the form:
\begin{eqnarray*}
 \Big(\bd^k_N \Big)_{T_k} &=& \Big( A^\top_{T_k} A_{T_k}  \Big)^{-1} \Big(
  A_{T_k}^\top A_{T_k^c}  \bx^k_{T_k^c} - A_{T_k}^\top ( A \bx^k - \bb )
 \Big) \\
 &=& \Big( A^\top_{T_k} A_{T_k}  \Big)^{-1} \Big(
 A_{T_k}^\top A_{T_k^c}  \bx^k_{T_k^c} - A_{T_k}^\top ( A_{T_k} \bx^k_{T_k} + A_{T_k^c} \bx^k_{T_k^c}  - \bb )
 \Big) \\
 &=& - \bx^k_{T_k} + \Big( A^\top_{T_k} A_{T_k}  \Big)^{-1}  A_{T_k}^\top \bb.
\end{eqnarray*}
Since the unit line search steplength $\alpha_k = 1$ is always accepted for all $k$ sufficiently large
(say, $k \ge k_0$), we have
\[
  \bx^{k+1} = \bx^k(\alpha_k) = \bx^k(1) = \left[
   \begin{array}{c}
     \bx^k_{T_k} + (\bd^k_N)_{T_k} \\
     0
   \end{array}
  \right]
  = \left[
  \begin{array}{c}
  \Big( A^\top_{T_k} A_{T_k}  \Big)^{-1}  A_{T_k}^\top \bb \\
  0
  \end{array}
  \right]
\]
Equivalently,
\[
  \bx^{k+1} = \arg\min \left\{  \| \bb - A \bz \|: \ \ \support(\bz) \subseteq T_k  \right\}.
\]
Consequently, \NHTP\ eventually (when $k \ge k_0$) becomes $\HTP^\eta$ of \cite{foucart2011hard}:
\[
  \HTP^\eta: \quad \left\{
   \begin{array}{ll}
    T_k & = \left\{\mbox{the best $s$ support of} \ (\bx^k - \eta \nabla f(\bx^k) )  \right\} , \quad
        (\mbox{i.e.,} \ T_k \in \T(\bx^k; \eta) )  \\ [1ex]
    \bx^{k+1} & = \arg\min \left\{  \| \bb - A \bz \|: \ \ \support(\bz) \subseteq T_k  \right\}.
   \end{array}
  \right .
\]
\cite[Prop.~3.2]{foucart2011hard} states that $\HTP^\eta$ will converge provided that $\eta \| A \|_2^2 < 1$, which is ensured when $\eta < 1/M_{2s}$. Our choice $\eta \le w/(8\mu_{2s}^2)$ apparently
satisfies this condition. Hence, \NHTP\ eventually enjoys all the good properties stated for $\HTP^\eta$
under the same conditions assumed in \cite{foucart2011hard} as long as the $\eta$ (note: $\mu$ is used in
\cite{foucart2011hard} instead of $\eta$) used there does not clash with our choice.

Since the Newton direction is always accepted as the search direction every iteration, one may wonder why
we did not just use the unit steplength $\alpha_k = 1$.
We note that \NHTP\ does not just seek for the next iterate satisfying $f(\bx^{k+1}) \le f(\bx^k)$,
it also  requires it to deduce a sufficient decrease by the quantity $\alpha_k \sigma \langle \nabla f(\bx^k), \ \bd^k \rangle$, which is proportional to the steplength $\alpha_k$.
Newton's direction $\bd^k_N$ with the unit steplength may not provide this proportional decrease and hence the unit steplength cannot be accepted in this case
(but the unit steplength will be eventually accepted).
In contrast, the \HTP\ family algorithms of \cite{foucart2011hard} only require a decrease
$f(\bx^{k+1}) \le f(\bx^k)$.
It is interesting to note that, in optimization,
one of the guidelines in designing a descent algorithm
is to ensure it deduces a sufficient decrease every iteration (see \cite{NW1999}) in order to achieve
desirable convergence properties. \\

\noindent {\bf Remark 10.} (On the gradient direction)
When the information on $\mu_{2s}$ and $m_{2s}$ is difficult to estimate, the choice of (\ref{CS-Parameters}) may not be possible.  On the one hand, those are the sufficient conditions for the
Newton direction to be accepted. Numerical experiments show that Newton's direction is often accepted with
a wide range of parameter choices. On the other hand, we have the restricted gradient direction to
rescue if the Newton direction is not deemed to be good enough in terms of the condition (\ref{Descent-Inequality-Newton}). The resulting algorithm still enjoys the global convergence in
Thm.~\ref{Thm-Global-Convergence} even if all search directions are of gradients.
It is interesting to note that a restricted gradient method was also proposed in \cite{foucart2011hard}
and is referred to as fast \HTP. We describe this algorithm (with just one gradient iteration each step)
in terms of our technical terminologies.
\[
 \FHTP^\eta: \quad \left\{
 \begin{array}{ll}
  \widetilde{\bx}^{k+1} &= \ \P_s( \bx^k - \eta \nabla f(\bx^k) ) \\
  T_{k+1}  & \in \ \T( \widetilde{\bx}^{k+1}; \eta ) \\
  \bx^{k+1}_{T_{k+1}} &= \
  \Big( \widetilde{\bx}^{k+1} - t_{k+1} \nabla f(  \widetilde{\bx}^{k+1} ) \Big)_{T_{k+1}} \ \
  \mbox{and} \ \ \bx^{k+1}_{T_{k+1}^c} = 0,
 \end{array}
 \right .
\]
where $t_{k+1}$ can be set to $1$ or chosen adaptively.
Despite it being also shown to enjoy similar convergence properties as $\HTP^\eta$ in \cite{foucart2011hard}, it does not fall within the framework (\ref{General-Framework})
and (\ref{Globalization-Step}).
A noticeable difference is that $\FHTP^\eta$ solves two optimization problems each step: one for
$\widetilde{\bx}^{k+1}$ and the other for $\bx^{k+1}_{T_{k+1}}$.
It would be interesting to see how the convergence analysis conducted in this paper can be extended to
$\FHTP^\eta$.

\section{Numerical Experiments} \label{Section-Numerical}

In this part, we show experimental results of \NHTP\ in CS (Sect.~\ref{Subsection-CS}) and
sparse logistic regression (Sect.~\ref{Subsection-SLR}) on both synthetic and real data.
A general conclusion is that \NHTP\ is capable of producing solutions of high quality and
is very fast when benchmarked against  six leading solvers  from compressed sensing 
and three solvers from sparse logistic regression. All experiments were conducted by using MATLAB (R2018a)  on a desktop of 8GB memory and Inter(R) Core(TM) i5-4570 3.2Ghz CPU.

We first describe how \NHTP\ was set up.
We  initialize  \NHTP\ with $\bx^0=0$ if $\nabla f(0)\neq0$ and $\bx^0=\textbf{1}$ if $\nabla f(0)=0$. Parameters are set as $\sigma=10^{-4}/2, \beta=0.5$.  For  $\gamma$, theoretically any positive $\gamma\leq m_{2s}$ is fine, but in practice to guarantee more steps using Newton directions, it is supposed to be relatively small \citep{de1996semismooth,facchinei1997nonsmooth}. Thus we choose $\gamma=\gamma_k$ with updating  
 $$ \gamma_k=\left\{\begin{array}{ccc}
                      10^{-10}, & {\rm if}&\bx^k_{  T^c_{k}}=0,  \\
                       10^{-4}, &   {\rm if}& \bx^k_{  T^c_{k}}\neq0.
                     \end{array}
 \right.$$
 For parameter $\eta $, in spite of that Theorem \ref{Thm-Qudratic-Convergence} has suggested to set $0<\eta <\overline{\eta }$, it is still difficult to fix a proper one since $M_{2s}$ is not easy to compute in general. Overall, we choose to update $\eta $ adaptively. Typically, we use the following rule: starting $\eta $ with a fixed scalar associated with the dimensions of a problem and then  update it as,
\begin{eqnarray*}
\eta _0&=&\frac{10(1+s/n)}{\min\{10,\ln(n)\}}>1,\\
\eta _{k+1}&=&\left\{
\begin{array}{ll}
\eta _k/1.05,&\text{if}~ {\rm mod} (k,10)=0~\text{and}\ \|F_{\eta _k}(\bx^{k};T_{k})\|>  k^{-2}, \\
1.05\eta _k,&\text{if}~ {\rm mod} (k,10)=0~\text{and}\ \|F_{\eta _k}(\bx^{k};T_{k})\|\leq  k^{-2}, \\
\eta _k,&\text{otherwise}.
\end{array}
\right.
\end{eqnarray*}
where {\rm mod} $(k,10)=0$ means $k$ is a multiple of $10$.
We terminate our method if at $k$th step it meets one of the following conditions:
 \begin{itemize}
 \item $\mbox{Tol}_{\eta_k}(\bx^k;\; T_k) \leq 10^{-6}$, where $\mbox{Tol}_{\eta}(\bx;\; T) $ is defined as (\ref{Tolerance-Function});
 \item $|f(\bx^{k+1})-f(\bx^{k})|<10^{-6}(1+|f(\bx^{k})|)$.
 \item $k$ reaches the maximum number (e.g., 2000) of iterations.
 \end{itemize}

\subsection{Compressed Sensing} \label{Subsection-CS}

Compressed sensing (CS) has seen revolutionary advances both in theory and
algorithms over the past decade. Ground-breaking papers that pioneered the advances are \citep{donoho2006compressed, candes2006robust,candes2005decoding}. The model   is described as in (\ref{CS})

{\bf a) Testing examples.}  We will focus on the exact recovery $\bb=A \bx$ by utilizing the sensing matrix $A$  chosen as in \citep{yin2015minimization,zhou2016null}.
\begin{example}[Gaussian matrix]\label{ex-gau}
Let $A\in\mathbb{R}^{m\times n}$  be a random Gaussian matrix  with each column $A_j, j\in N_n$ being identically and independently generated from the standard normal distribution. We then normalize each column such that $\|A_j\|=1$. Finally, the `ground  truth' signal $\bx^*$ and the measurement  $\bb$ are produced by the following pseudo Matlab codes:
\begin{eqnarray}\label{gen-x*}
\bx^{*} = \verb"zeros"(n, 1),~~\Gamma = \verb"randperm"(n),~~
\bx^{*}(\Gamma(1 : s)) = \verb"randn"(s, 1),~~\bb = A \bx^*.
\end{eqnarray}
\end{example}
\begin{example}[Partial DCT  matrix]\label{ex-dct}
Let $A\in\mathbb{R}^{m\times n}$  be a random partial discrete cosine transform (DCT) matrix generated by
\begin{eqnarray*}
A_{ij}= \cos(2\pi (j-1) \psi_i),~~i=1,\ldots, m,~~j=1,\ldots, n
\end{eqnarray*}
where $\psi_i,i=1,\ldots, m$ is uniformly and independently sampled from $[0,1]$. We then normalize each column such that $\|A_j\|=1$ with $\bx^*$ and  $\bb$ being generated the same way as in
 Example \ref{ex-gau}.
\end{example}

{\bf b) Benchmark methods.} There exists a large number of numerical methods for
 the CS problem (\ref{CS}). It is beyond the scope of this paper to compare them all. 
We selected six state-of-the-art methods.  
They are {\tt HTP}\; \citep{foucart2011hard}\setcounter{footnote}{0}\footnote{
{\tt HTP} is available at:\emph{  https://github.com/foucart/HTP.}}, {\tt NIHT} \citep{blumensath2010normalized}\footnote{{\tt NIHT}, {\tt GP} and {\tt OMP} are available at  \emph{https://www.southampton.ac.uk/engineering/about/staff /tb1m08.page$\#$software}. We use the version {\tt sparsify\_0\_5} in which {\tt NIHT}, {\tt GP} and {\tt OMP} are called {\tt hard\_l0\_Mterm}, {\tt greed\_gp} and {\tt greed\_omp}.},  {\tt GP} \citep{blumensath2008gradient}\footref{2}, {\tt OMP} \citep{pati1993orthogonal,tropp2007signal}\footref{2}, {\tt CoSaMP} \citep{needell2009cosamp}\footref{3} and  {\tt SP} \citep{dai2009subspace}\footnote{
{\tt CoSaMP} and  {\tt SP} are available at:\emph{  http://media.aau.dk/null\_space\_pursuits/2011/07/a-few-corrections-to-cosamp-and-sp-matlab.html.}}.  For {\tt HTP}, set {\tt MaxNbIter}=1000 and {\tt mu}={\tt `NHTP'}. For {\tt NIHT}, the maximum iteration {\tt `maxIter'} is set as 1000 and ${\tt M}=s$.  For {\tt GP} and {\tt OMP}, the {\tt `stopTol'} is set as $1000$. For {\tt CoSaMP} and  {\tt SP}, set {\tt tol}$=10^{-6}$ and {\tt maxiteration}$=1000$. Notice that the first three methods prefer solving sensing matrix $A$ with unit columns, which is the reason for us to normalize each generated $A$ in  Example \ref{ex-gau} and Example \ref{ex-dct}. Let $\bx$ be the solution produced by a method. We say a recovery of this method is successful if $\|\bx-\bx^*\|<0.01\|\bx^*\|.$

\begin{figure}[H]
\centering
\begin{subfigure}{0.48\textwidth}
\caption{Gaussian Matrix}
  \includegraphics[width=.9\linewidth]{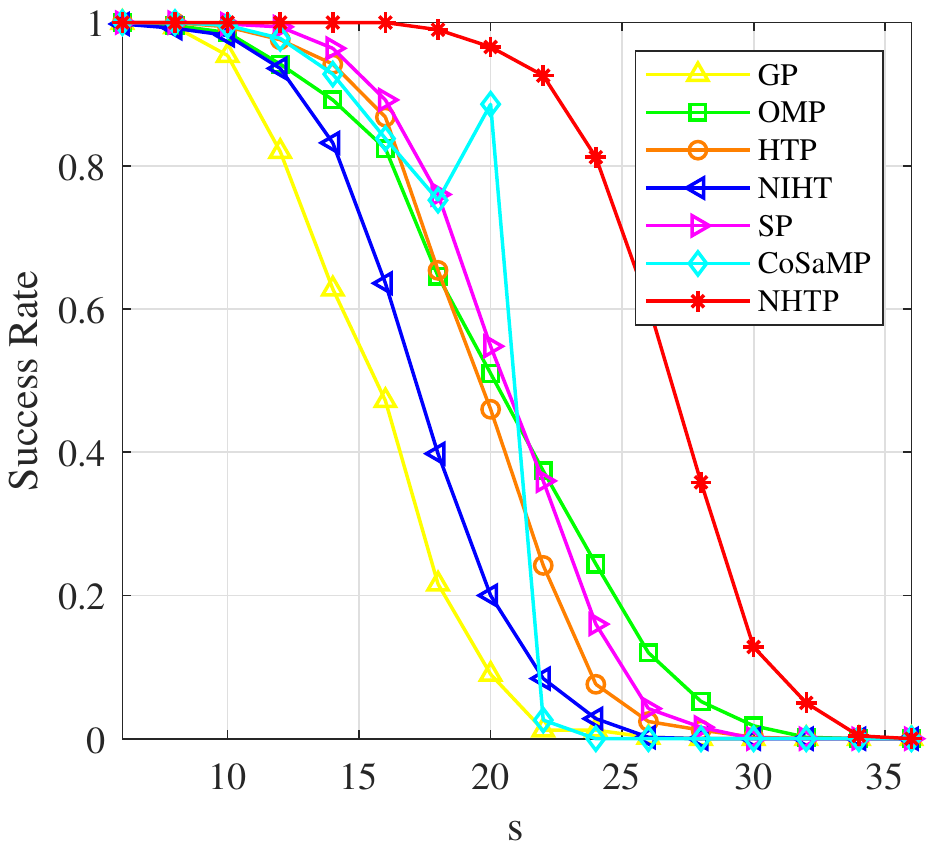}
\end{subfigure}
\begin{subfigure}{0.48\textwidth}
 \caption{Partial DCT Matrix}
  \includegraphics[width=.9\linewidth]{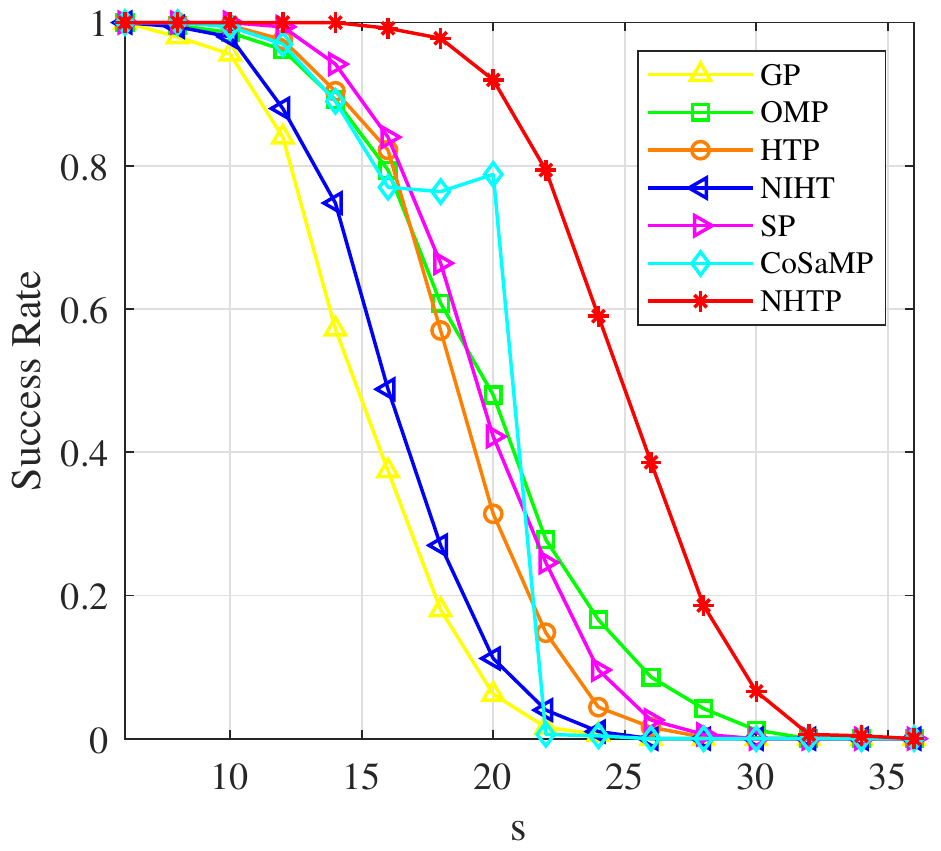}
\end{subfigure}
\caption{Success rates. $n = 256, m=\lceil n/4\rceil,  s\in\{6,8,\cdots, 36\}$.}
\label{fig:SuccRate}
\end{figure}
\begin{figure}[H]  
\centering
\begin{subfigure}{0.48\textwidth}
 \caption{Gaussian Matrix}
 \includegraphics[width=.9\linewidth]{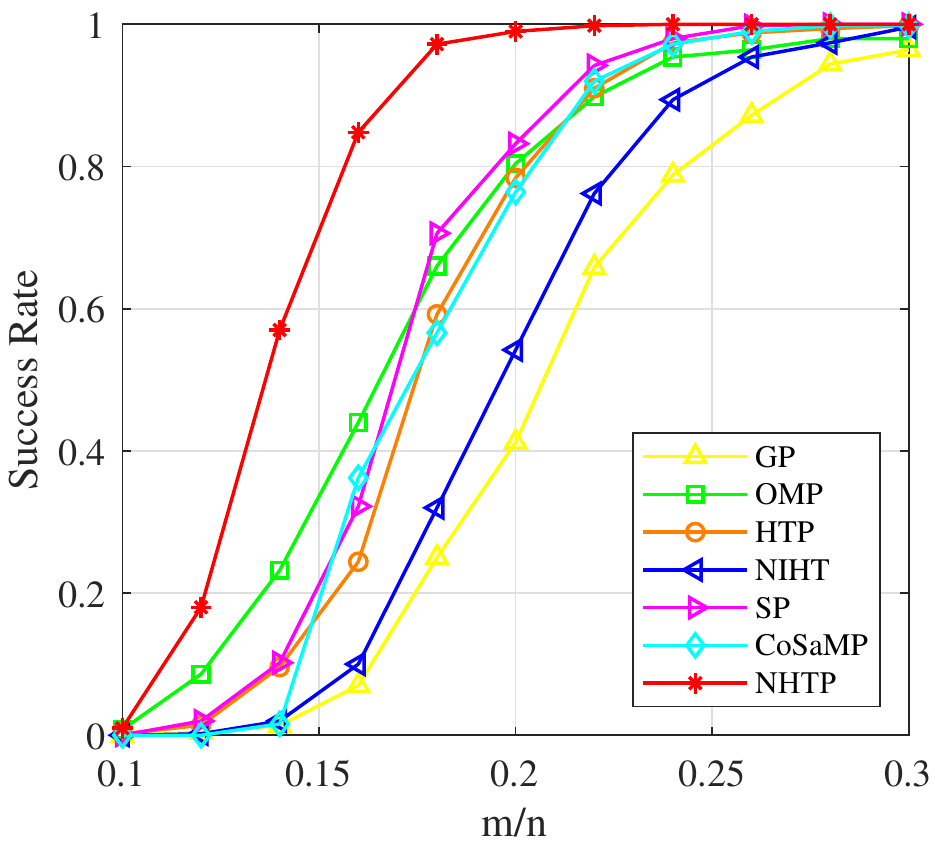}
\end{subfigure}
\begin{subfigure}{0.48\textwidth}
  \caption{Partial DCT Matrix}
 \includegraphics[width=.9\linewidth]{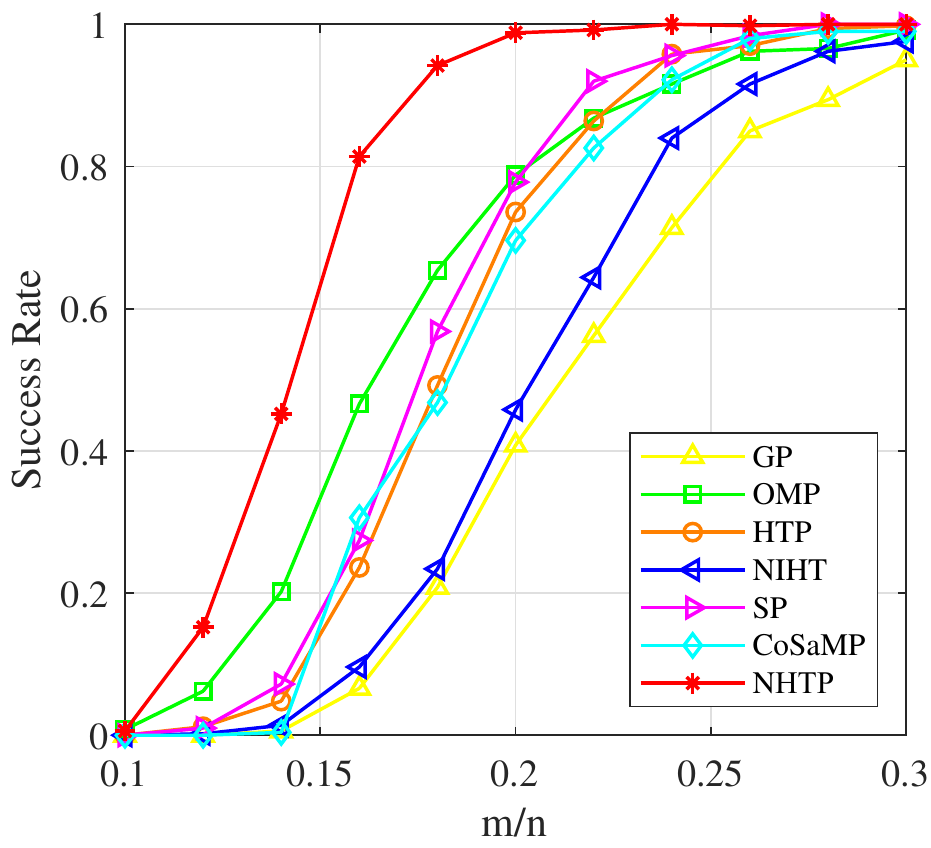}
\end{subfigure}
\caption{Success rates. $n = 256, s=\lceil 0.05 n\rceil, m=\lceil rn\rceil$ with $r\in\{  0.1, 0.12, \cdots, 0.3\}$.}
\label{fig:SuccRate-mn}
\end{figure}

{\bf c) Numerical comparisons.} We begin with running $500$ independent trials with fixed $n=256, m=\lceil n/4\rceil$ and recording the corresponding success rates (which is defined by the percentage of the number of successful recoveries over all trails) at sparsity levels $s$ from 6 to 36, where $\lceil a\rceil$ is the smallest integer that is no less than $a$. From Fig. \ref{fig:SuccRate}, one can observe that for both Example \ref{ex-gau} and Example \ref{ex-dct}, \NHTP\ yielded the highest success rate for each $s$. 
For example, when $s=22$ for Gaussian matrix, our method still obtained $90\%$ successful recoveries while the other methods only guaranteed less than $40\%$ successful ones. 
Moreover, {\tt OMP}, {\tt SP} and {\tt HTP} generated similar results, and {\tt GP} and {\tt NIHT} always came the last. Next we run $500$ independent trials with fixing  $n=256, s=\lceil 0.05 n\rceil$ but varying $m=\lceil r n\rceil$ where $r\in\{0.1,0.12,\cdots,0.3\}$. 
It is clearly to be seen that  the larger $m$ is, the easier the problem becomes to be solved.
This is illustrated by Fig. \ref{fig:SuccRate-mn}.  Again  \NHTP\ outperformed the others due to highest success rate for each $s$, and {\tt GP} and {\tt NIHT} still came the last. 

To see the accuracy of the solutions and the speed of these seven methods, we now run 50 trials for each kind of matrices with higher dimensions $n$ increasing from 5000 to 25000 and keeping $m=\lceil n/4\rceil, s=\lceil 0.01n\rceil, \lceil0.05n\rceil$.   Specific results produced by these seven methods are recorded in Tables \ref{tab:error} and \ref{tab:time}. Our method  \NHTP\ always obtained the most accurate recovery, with accuracy order of $10^{-14}$ or higher, followed by {\tt HTP}. {\tt NIHT} was stable at achieving the solutions with  accuracy of order $10^{-7}$.
Moreover, {\tt GP} and {\tt OMP} rendered solutions as accurate as those by  \NHTP\ when  $s=\lceil0.01n\rceil$, but yielded inaccurate ones when $s=\lceil0.05n\rceil$, which means that these two methods worked well when the solution is very sparse. 
In contrast,  {\tt SP} and {\tt CoSaMP} always generated results with worst accuracy.  When it comes to the computational speed in Table  \ref{tab:time},   \NHTP\ is the fastest for most of the cases. 
The fast convergence of \NHTP\ becomes more superior in high dimensional data setting.
For example, when $n=25000$ and $s=\lceil0.05n\rceil$,  6.58 seconds  by \NHTP\ against 36.93 seconds by {\tt HTP}, which is the fastest method among the other five methods.   {\tt GP} and {\tt OMP} always ran the slowest. In addition, we also compared seven algorithms on Example  \ref{ex-gau}, but omitted all the related results since they were similar to those of Example   \ref{ex-dct}

 \begin{table}[H]
 \caption{Average absolute error $\|\bx-\bx^*\|$ for Example   \ref{ex-dct}. \label{tab:error}}\vspace{-3 mm}
{\renewcommand\baselinestretch{1.15}\selectfont
{\centering\begin{tabular}{ p{0.8cm}cccccccc }\\ \hline
$s$&$n$&{\tt GP} &  {\tt OMP}&{\tt HTP}&{\tt NIHT}&{\tt SP}&{\tt CoSaMP} &\NHTP \\\hline
\multirow{5}{*}{$\lceil0.01n\rceil$}  &	5000	&	2.78e-15	&	2.40e-15	&	2.97e-15	&	2.42e-7	&	1.12e-5	&	1.12e-5	&	4.59e-16	\\
	&	10000	&	5.21e-15	&	4.75e-15	&	5.70e-15	&	3.26e-7	&	3.59e-5	&	3.59e-5	&	1.10e-15	\\
	&	15000	&	7.05e-15	&	7.07e-15	&	7.36e-15	&	4.28e-7	&	4.25e-5	&	4.25e-5	&	1.39e-15	\\
	&	20000	&	9.49e-15	&	9.06e-15	&	9.47e-15	&	4.88e-7	&	6.56e-5	&	6.56e-5	&	1.88e-15	\\
	&	25000	&	1.15e-14	&	1.12e-14	&	1.11e-14	&	5.32e-7	&	1.78e-4	&	1.78e-4	&	2.47e-15	\\\hline
\multirow{5}{*}{$\lceil0.05n\rceil$}	&	5000	&	1.28e-03	&	1.40e-03	&	1.26e-14	&	4.80e-7	&	9.07e-5	&	9.07e-5	&	5.94e-15	\\
	&	10000	&	7.91e-04	&	3.56e-04	&	2.44e-14	&	6.86e-7	&	1.77e-4	&	1.77e-4	&	1.18e-14	\\
	&	15000	&	1.10e-03	&	6.20e-04	&	3.57e-14	&	8.54e-7	&	2.11e-4	&	2.11e-4	&	1.76e-14	\\
	&	20000	&	9.43e-04	&	3.33e-04	&	4.87e-14	&	9.80e-7	&	3.53e-4	&	3.53e-4	&	2.39e-14	\\ 
	&	25000	&	1.24e-03	&	5.57e-04	&	5.94e-14	&	1.01e-6	&	2.59e-4	&	2.59e-4	&	2.86e-14	\\\hline
    \end{tabular}\par} }
\end{table}

\begin{table}[H]
 \caption{Average CPU time (in seconds) for Example  \ref{ex-dct}. \label{tab:time}}\vspace{-3 mm}
{\renewcommand\baselinestretch{1.1}\selectfont
{\centering\begin{tabular}{ ccccccccc }\\ \hline
$s$&$n$&{\tt GP} &  {\tt OMP}&{\tt HTP}&{\tt NIHT}&{\tt SP}&{\tt CoSaMP} &\NHTP \\\hline
\multirow{5}{*}{$\lceil0.01n\rceil$}	&	5000	&	0.69	&	0.48	&	0.09	&	0.30	&	0.07	&	0.05	&	0.06	\\
	&	10000	&	4.47	&	3.70	&	0.33	&	1.21	&	0.31	&	0.25	&	0.16	\\
	&	15000	&	14.57	&	13.41	&	0.74	&	2.96	&	0.96	&	0.86	&	0.37	\\
	&	20000	&	32.70	&	30.46	&	1.34	&	5.53	&	2.30	&	2.00	&	0.65	\\
	&	25000	&	68.94	&	67.13	&	2.49	&	37.03	&	20.11	&	4.18	&	1.13	\\\hline
\multirow{5}{*}{$\lceil0.05n\rceil$}	&	5000	&	3.52	&	3.22	&	0.23	&	1.29	&	0.90	&	1.43	&	0.28	\\
	&	10000	&	19.84	&	23.55	&	1.52	&	4.63	&	6.02	&	15.56	&	0.79	\\
	&	15000	&	67.79	&	77.30	&	7.25	&	10.43	&	23.03	&	60.87	&	2.20	\\
	&	20000	&	151.28	&	177.00	&	18.02	&	18.70	&	58.20	&	148.83	&	3.49	\\ 
	&	25000	&	312.57	&	363.44	&	36.93	&	78.69	&	153.52	&	307.53	&	6.58	\\\hline

 \end{tabular}\par} }
\end{table}

\subsection{Sparse Logistic Regression} \label{Subsection-SLR}

Sparse logistic regression (SLR) has drawn extensive attention since it was first proposed by \citet*{tibshirani1996regression}. Same as \citep{bahmani2013greedy}, we will address the so-called  $\ell_2$ norm regularized sparsity constrained logistic regression (SCLR) model, namely,
\begin{equation}\label{SLR-L2}
\min_{\|\bx\|_0\leq s}~~\ell(\bx) +\mu\|\bx\|_2^2 ~~~~~~~~ {\rm with}~~~~~~ \ell(\bx):= \frac{1}{m}\sum_{i=1}^{m}\left\{\ln(1+ e^{\langle\ba_i,\bx\rangle})-b_i\langle\ba_i,\bx\rangle\right\},
\end{equation}
where $\ba_{i} \in\mathbb{R}^{n}, b_{i}\in\{0,1\}, i=1,\ldots, m$ are respectively given $m$ features and responses/labels, and $\mu>0$ (e.g. $\mu=10^{-6}/m$). 
The employment of a regularization was well justified because otherwise 
\textit{`one can achieve arbitrarily small loss values by tending the parameters to infinity along certain directions'} (see \citep{bahmani2013greedy}).
This is the reason why we will only focus on  (\ref{SLR-L2}).

{\bf d) Testing examples.} We will test three types of data sets. 
The first two are synthetic and the last one is from a real database.  
One synthetic data is adopted from \citep{lu2013sparse}, \citep{pan2017convergent}
with the features $[\ba_{1}~\cdots~\ba_{m}]$ being generated 
identically and independently.
The other is the same as  \citet{agarwal2010fast} or \citet{ bahmani2013greedy} who have considered independent features with each $\ba_{i}$ being generated by an autoregressive process \citep{hamilton1994time}.

\begin{example}[Independent Data \citep{lu2013sparse,pan2017convergent}] \label{log-EX1} To generate data labels $\bb\in\{0,1\}^m$, we first randomly separate $\{1,\ldots,m\}$ into two parts $I$ and $ I^c$ and  set $b_i=0$ for $i\in I$ and $b_i=1$ for $i\in I^c$. Then the feature data is produced by
$$\ba_i=y_iv_i\textbf{1}+\bw_i,\ \ \ \ i=1,\ldots, m$$
with $\mathbb{R}\ni~v_i\sim \mathcal{N} (0,1)$, $\mathbb{R}^n\ni\bw_i\sim \mathcal{N} (0,\mathcal{I}_n)$ and $\mathcal{N} (0,\mathcal{I}_n)$ is the normal distribution with zero mean and the identity
covariance. Since the  sparse parameter $\bx^*\in\mathbb{R}^n$ is unknown, 
different sparsity level$s$ will be tested. 
\end{example}

\begin{example}[Correlated Data \citep{agarwal2010fast, bahmani2013greedy}]\label{log-EX2}
The\\ sparse parameter $\bx^*\in\mathbb{R}^n$  has $s$ nonzero entries drawn independently from the standard Gaussian distribution. Each data sample $\ba_i=[a_{i1}~\cdots~a_{in}]^\top, i=1,\ldots, m$ is an independent instance of the random vector
generated by an autoregressive process  \citep[see][]{hamilton1994time} 
$$
   a_{i(j+1)}=\theta a_{ij}+\sqrt{1-\eta^2}v_{ij}, \quad j=1,\ldots, n-1,
$$
with $a_{i1}\sim \mathcal{N} (0,1)$, $v_{ij}\sim \mathcal{N} (0,1)$ and $ \theta\in[0,1]$ being the correlation parameter. The data labels $\by\in\{0,1\}^m$
are then drawn randomly according to the Bernoulli distribution with
$${\rm Pr}\{y_i = 0 | \ba_i\} = \left[1+e^{\langle \ba_i,\bx^*\rangle}\right]^{-1},\ \ \ \ i=1,\ldots, m.$$ 
\end{example}
\begin{example}[Real data]\label{log-EX3} This example comprises of seven real data sets for binary classification. They are
\setcounter{footnote}{0}
\texttt{colon-cancer\footnote{$\texttt{https://www.csie.ntu.edu.tw/}\sim\texttt{cjlin/libsvmtools/datasets/}$\label{1}}},
\texttt{arcene\footnote{$\texttt{http://archive.ics.uci.edu/ml/index.php}$\label{2}}}, 
\texttt{newsgroup\footnote{$\texttt{https://web.stanford.edu/}\sim\texttt{hastie/glmnet\_matlab/}$\label{3}}}, 
\texttt{news20.binary\footref{1}}, 
\texttt{duke breast-} \texttt{cancer\footref{1}}, 
\texttt{leukemia\footref{1}}, 
\texttt{rcv1.binary\footref{1}}, which are summarized in the following table, where the last three data sets have testing data.
Moreover, as described in the website\footref{1},  for the four data with small sample sizes: \texttt{colon-cancer}, \texttt{arcene}, \texttt{duke breast-cancer} and \texttt{leukemia},  sample-wise normalization has been conducted so that each sample has  mean zero and variance one, and then feature-wise normalization has been conducted so that each feature has  mean zero and variance one. 
For the rest four  data with larger sample sizes, they are feature-wisely scaled to $[-1,1]$. 
All $-1$s in  classes $\bb$ are replaced by 0. 
\begin{table}[H]
\begin{center}
\begin{tabular}{lcccc}\\\hline
Data name&$m$ samples&$n$ features&training size $m_1$ &testing size $m_2$\\ \hline
\texttt{colon-cancer}&62 &2000&62&0\\
\texttt{arcene}&100 &10000&100&0\\
\texttt{newsgroup}&11314 &777811 &11314&0\\
\texttt{news20.binary}&19996&1355191&19996&0\\
\texttt{duke breast-cancer}&42	&7129&38&4\\
\texttt{leukemia}&72		&7129&38&34\\
\texttt{rcv1.binary}&	40242&47236&20242&20000\\\hline
\end{tabular}
\end{center}
\end{table}
\end{example}
{\bf e) Benchmark methods.} Since there are numerous leading solvers that have been proposed to solve SLR problems, we again only focus on those dealing with the $\ell_2$ norm regularized SCLR. 
We select three solvers: 
{\tt GraSP}\;\citep{bahmani2013greedy}\footnote{$\texttt{http://sbahmani.ece.gatech.edu/GraSP.html}$\label{4}}, {\tt NTGP}\;\citep{yuan2014newton} 
and {\tt IIHT}\;\citep{pan2017convergent}. 
Notice that all those methods are used to solve $\ell_2$ norm regularized SCLR model (\ref{SLR-L2}) with  $\mu=10^{-6}/m$. 
Except for  {\tt IIHT}, which only used the first order information such as objective values or gradients,  the other three methods exploit second order information of the objective function. 
{\tt NTGP} integrates  Newton directions into some steps, 
and {\tt GraSP} takes advantage of the Matlab built-in function: {\tt minFunc} which calls a Quasi-Newton strategy. 
For {\tt GraSP}, if we use its defaults parameters,
it would be less likely to meet its stopping criteria before the number of iteration reaching the maximal one. 
Compared with other three methods, which all generate a sequence with decreasing  objective function values,  the objective function value at each iteration by  {\tt GraSP} fluctuated greatly. 
Therefore, we set an extra stopping criterion for {\tt GraSP}: $f(\bx^{k})-f(\bx^{k+1})<10^{-6}$. 
And if $f(\bx^{k})<f(\bx^{k+1})$, then terminate it and output $\bx^{k}$. For {\tt NTGP}, to facilitate its computational speed, we set {\tt maxIter}=20 for outer loops, and {\tt maxIter$\_$sub}=50 and {\tt optTol$\_$sub} $= 10^{-3}$ for inner loops.  For   {\tt IIHT}, we keep its default parameters.

For both Example \ref{log-EX1} and Example \ref{log-EX2}, we run $500$ independent trials if $n<10^3$ and $50$ independent trials otherwise, and report the average  logistic loss $\ell(\bx)$ and CPU time to demonstrate the performance of each method.  

\begin{figure}[H]
\centering
\begin{subfigure}{0.48\textwidth}
 \caption{~}\vspace{-3mm}
  \includegraphics[width=.9\linewidth]{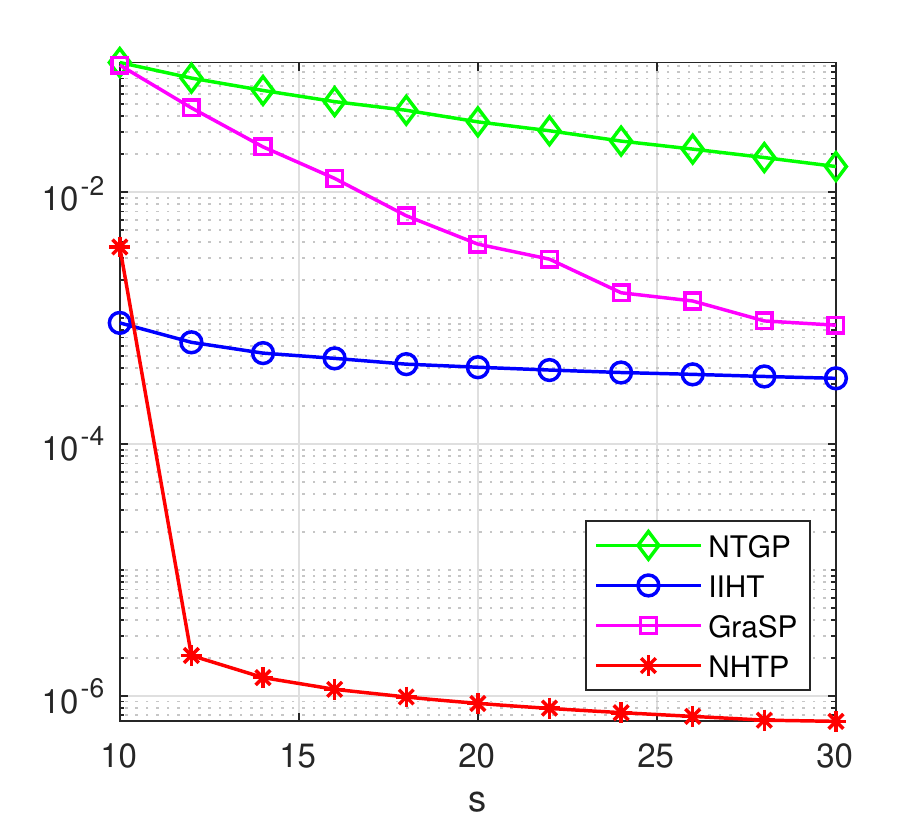}
\end{subfigure}
\begin{subfigure}{0.48\textwidth}
 \caption{~}\vspace{-3mm}
  \includegraphics[width=.9\linewidth]{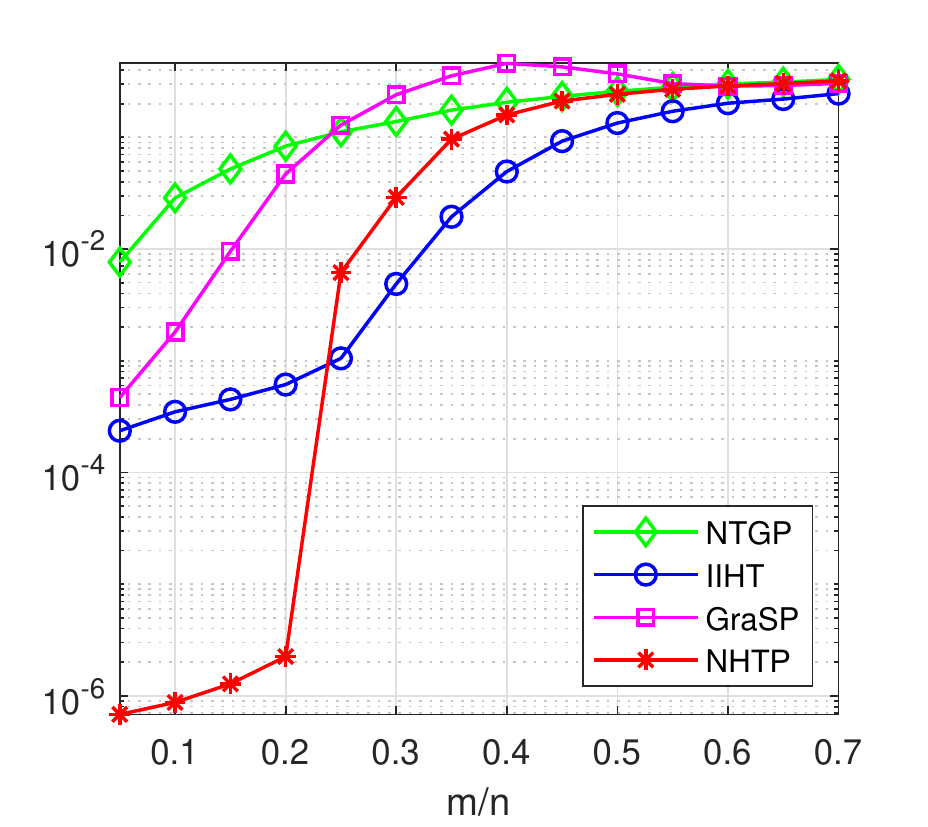}
\end{subfigure}
\caption{Average logistic loss $\ell(\bx)$ of four methods for Example \ref{log-EX1}.}
\label{fig:log-ex1-s-mn}
\end{figure}

{\bf f) Numerical comparisons.} For Example \ref{log-EX1}, we begin with testing each method 
for the case
$n=256$ and $m=\lceil n/5\rceil$ with  varying sparsity levels $s$ from 10 to 30.  
From Fig.~\ref{fig:log-ex1-s-mn}(a), one can observe that  {\tt IIHT} rendered the best $\ell(\bx)$  when $s=10$ and  \NHTP\ performed the best $\ell(\bx)$ when $s> 10$.  
And importantly, the value $\ell(\bx)$ produced by \NHTP\ for each instance is far smaller than others, with order about $10^{-6}$.    
We then test the case  $n=256, s=\lceil 0.05n\rceil$ and $m=\lceil rn\rceil$ with varying $r\in\{0.05,0.1,\cdots,0.7\}$. 
From Fig.~\ref{fig:log-ex1-s-mn}(b), $\ell(\bx)$ generated by \NHTP\  is the lowest  when the sample size was relatively small, and it gradually approached to the values similar to those obtained by the others. 
{\tt IIHT} performed the best in terms of $\ell(\bx)$ when $m/n>0.2$ 
and {\tt GraSP} always rendered the highest loss.  

 \begin{table}[H]
 \caption{Average logistic loss $\ell(\bx)$ and CPU time (in seconds) for Example   \ref{log-EX1}. \label{tab:lx-time-1}}\vspace{-3 mm}
{\renewcommand\baselinestretch{1.1 }\selectfont
{\centering\begin{tabular}{p{1cm} p{1cm}|p{1.25cm}p{1.25cm}p{1.25cm}p{1.25cm}|p{.8cm}p{.8cm}p{.8cm}p{.8cm}   }\\ \hline
\multirow{2}{*}{ $s$}&\multirow{2}{*}{ $n$}  &\multicolumn{4}{c|}{$\ell(\bx)$} &\multicolumn{4}{c}{CPU Time}\\\cline{3-10}
&&{\tt NTGP} &  {\tt IIHT}&{\tt GraSP} &\NHTP
&{\tt NTGP} &  {\tt IIHT}&{\tt GraSP} &\NHTP \\\hline
\multirow{7}{*}{$\lceil0.01n\rceil$}	&	10000	&	2.39e-1	&	1.43e-1	&	2.44e-1	&	2.26e-1	&	8.403	&	1.723	&	0.488	&	0.313	\\
	&	15000	&	2.48e-1	&	1.37e-1	&	2.39e-1	&	2.28e-1	&	17.81	&	3.307	&	0.974	&	0.457	\\
	&	20000	&	2.35e-1	&	1.36e-1	&	2.36e-1	&	2.20e-1	&	32.61	&	6.245	&	1.862	&	0.842	\\
	&	25000	&	2.25e-1	&	1.29e-1	&	2.30e-1	&	2.11e-1	&	52.99	&	8.913	&	3.006	&	1.372	\\
	&	30000	&	2.24e-1	&	1.24e-1	&	2.30e-1	&	2.07e-1	&	76.31	&	14.15	&	4.309	&	2.140	\\
	&	35000	&	2.21e-1	&	1.23e-1	&	2.29e-1	&	2.08e-1	&	149.7	&	21.84	&	16.08	&	2.875	\\
	&	40000	&	2.18e-1	&	1.21e-1	&	2.32e-1	&	2.05e-1	&	466.1	&	29.12	&	804.2	&	3.923	\\\hline
\multirow{7}{*}{$\lceil0.05n\rceil$}	&	10000	&	4.58e-2	&	4.76e-4	&	4.97e-3	&	6.50e-7	&	9.931	&	3.094	&	1.795	&	0.987	\\
	&	15000	&	4.05e-2	&	4.69e-4	&	7.77e-3	&	3.32e-7	&	26.34	&	6.218	&	4.069	&	2.442	\\
	&	20000	&	4.10e-2	&	4.80e-4	&	8.24e-3	&	6.32e-7	&	51.29	&	10.69	&	5.695	&	4.315	\\
	&	25000	&	4.56e-2	&	4.90e-4	&	6.06e-3	&	4.77e-7	&	54.96	&	15.93	&	8.964	&	7.004	\\
	&	30000	&	4.17e-2	&	4.92e-4	&	6.49e-3	&	6.89e-7	&	85.22	&	23.54	&	11.79	&	11.15	\\
	&	35000	&	3.95e-2	&	4.89e-4	&	6.46e-3	&	6.65e-7	&	182.1	&	35.97	&	24.34	&	17.25	\\
	&	40000	&	3.84e-2	&	4.92e-4	&	7.54e-3	&	5.81e-7	&	551.1	&	55.00	&	619.5	&	25.41	\\\hline
    \end{tabular}\par} }
\end{table}
 When the size of example is becoming relatively large, the picture is significant different. Hence we now run 50  independent trials with higher dimensions $n$ increasing from 10000 to 40000 and keeping $m=\lceil n/5\rceil, s= \lceil0.01n\rceil, \lceil0.05n\rceil$. As presented in Table \ref{tab:lx-time-1}, when $s= \lceil0.01n\rceil$, {\tt IIHT} produced the lowest $\ell(\bx)$, followed by \NHTP\ which was the fastest.  
But when $s= \lceil0.05n\rceil$, \NHTP\ outperformed others in terms of $\ell(\bx)$ with order of $10^{-7}$ which was much better than others.
The time used by \NHTP\ is also significantly less than the others,
for example, 25.41s by \NHTP\ vs. 619.5s by {\tt GraSP} when $n=40000$.

\begin{figure}[H]
	\centering
	\begin{subfigure}{0.48\textwidth}
		\caption{~}\vspace{-3mm}
			\includegraphics[width=.95\linewidth]{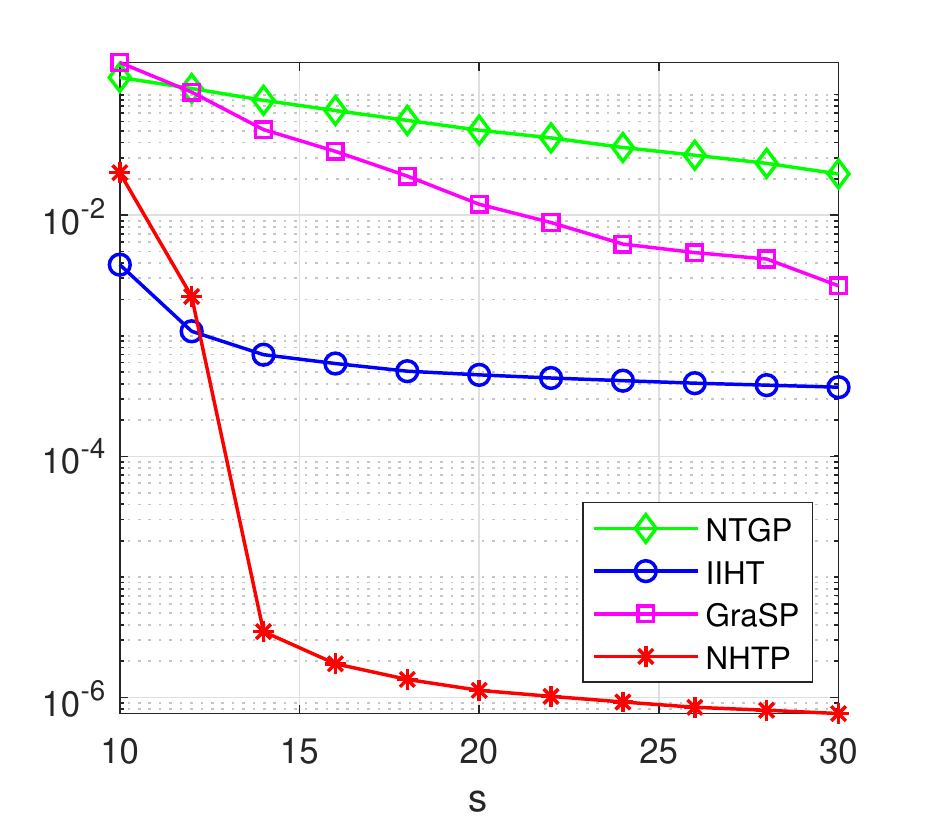}
	\end{subfigure}
	\begin{subfigure}{0.48\textwidth}
		\caption{~}\vspace{-3mm}
		\includegraphics[width=.95\linewidth]{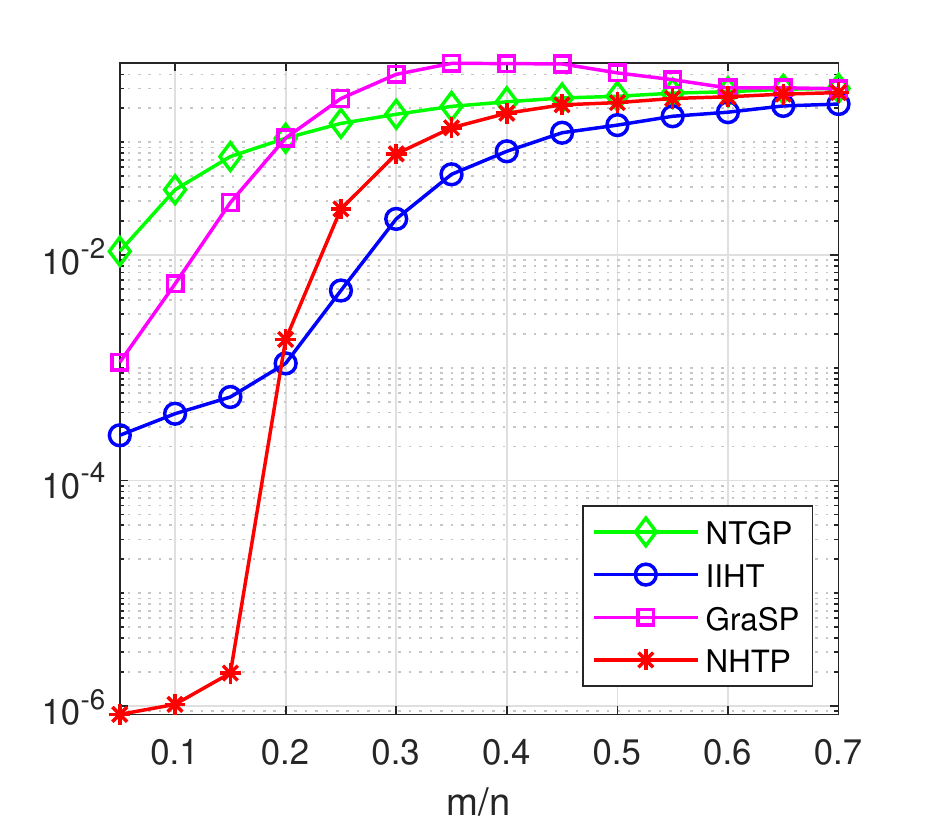}
	\end{subfigure}
	\caption{Average logistic loss $\ell(\bx)$ of four methods for Example \ref{log-EX2}.}
	\label{fig:log-ex2-s-mn}
\end{figure}

 \begin{table}[H]
 \caption{Average logistic loss $\ell(\bx)$ and CPU time (in seconds) for Example   \ref{log-EX2}. \label{tab:lx-time-2}}\vspace{-3 mm}
{\renewcommand\baselinestretch{1.1 }\selectfont
{\centering\begin{tabular}{p{1cm} p{1cm}|p{1.25cm}p{1.25cm}p{1.25cm}p{1.25cm}|p{.8cm}p{.8cm}p{.8cm}p{.8cm}   }\\ \hline
\multirow{2}{*}{ $s$}&\multirow{2}{*}{ $n$}  &\multicolumn{4}{c|}{$\ell(\bx)$} &\multicolumn{4}{c}{CPU Time}\\\cline{3-10}
&&{\tt NTGP} &  {\tt IIHT}&{\tt GraSP} &\NHTP
&{\tt NTGP} &  {\tt IIHT}&{\tt GraSP} &\NHTP \\\hline
\multirow{7}{*}{$\lceil0.01n\rceil$}	&	10000	&	1.87e-1	&	5.68e-2	&	1.93e-1	&	1.51e-1	&	8.338	&	4.394	&	0.471	&	0.245	\\
	&	15000	&	1.81e-1	&	4.07e-2	&	1.73e-1	&	1.25e-1	&	19.72	&	7.156	&	1.403	&	0.702	\\
	&	20000	&	1.61e-1	&	3.39e-2	&	1.64e-1	&	9.94e-2	&	36.68	&	10.74	&	2.370	&	1.194	\\
	&	25000	&	1.62e-1	&	2.62e-2	&	1.61e-1	&	9.84e-2	&	54.37	&	16.51	&	3.800	&	1.922	\\
	&	30000	&	1.63e-1	&	2.75e-2	&	1.63e-1	&	9.59e-2	&	124.4	&	40.11	&	18.83	&	6.067	\\
	&	35000	&	1.58e-1	&	2.09e-2	&	1.52e-1	&	8.73e-2	&	179.1	&	44.47	&	199.2	&	8.257	\\
	&	40000	&	1.59e-1	&	2.14e-2	&	1.57e-1	&	8.87e-2	&	423.4	&	46.13	&	639.4	&	19.47	\\\hline																			
\multirow{7}{*}{$\lceil0.05n\rceil$}	&	10000	&	7.59e-2	&	6.02e-4	&	2.18e-2	&	1.54e-6	&	9.101	&	3.426	&	1.875	&	0.880	\\
	&	15000	&	7.95e-2	&	6.15e-4	&	2.02e-2	&	1.67e-6	&	20.40	&	7.426	&	4.316	&	2.140	\\
	&	20000	&	7.84e-2	&	5.93e-4	&	2.34e-2	&	1.55e-6	&	34.91	&	12.51	&	6.394	&	4.015	\\
	&	25000	&	7.96e-2	&	5.97e-4	&	2.44e-2	&	1.65e-6	&	54.41	&	19.03	&	8.921	&	6.590	\\
	&	30000	&	7.76e-2	&	6.00e-4	&	2.04e-2	&	1.58e-6	&	107.2	&	29.95	&	16.57	&	10.09	\\
	&	35000	&	7.74e-2	&	6.01e-4	&	2.18e-2	&	1.61e-6	&	137.3	&	45.71	&	26.05	&	16.10	\\
	&	40000	&	7.89e-2	&	5.90e-4	&	2.41e-2	&	1.58e-6	&	305.8	&	70.83	&	721.0	&	22.46	\\\hline
    \end{tabular}\par} }
\end{table}

 For Example \ref{log-EX2}, it is related to the parameter $\theta$. 
We only report the results for $\theta=1/2$ since the comparisons of all methods are similar for each fixed $\theta\in(0,1)$. 
Again we first fix $n=256, m=\lceil n/5\rceil$ and  vary sparsity levels $s$ from 10 to 30. As shown in Fig. \ref{fig:log-ex2-s-mn} (a),   \NHTP\ yielded the smallest logistic loss  when $s>12$, followed by   {\tt IIHT}.
We then  fix  $n=256, s=\lceil 0.05n\rceil$ and change the sample size $m=\lceil rn\rceil$, where $r\in\{0.05,0.1,0.15,\cdots,0.7\}$. From Fig.~\ref{fig:log-ex2-s-mn}(b), \NHTP\ outperformed others  when the sample size was relatively small such as $m/n<0.2$, while {\tt IIHT} performed  best
in terms of $\ell(\bx)$ when $m/n\geq0.2$.  

When the size of the example is becoming relatively large, the picture again is significant different.
We run 50  independent trials with higher dimensions $n$ increasing from 10000 to 40000 and keeping $m=\lceil n/5\rceil, s= \lceil0.01n\rceil, \lceil0.05n\rceil$. 
As presented in Table \ref{tab:lx-time-2},  when $s= \lceil0.01n\rceil$ {\tt IIHT} indeed provided the best logistic loss and comparable to ours.
However, \NHTP\ was significantly faster than {\tt IIHT}.
Clearly, under the case of $s= \lceil0.05n\rceil$, \NHTP\ offered the far lowest $\ell(\bx)$ with order of $10^{-6}$ and CPU time with 22.46 seconds against 721 seconds from {\tt GraSP} when $n=40000$.

\begin{figure} 
\centering
\begin{subfigure}{0.43\textwidth}
\caption{$\ell(\bx)$} \vspace{-2mm}
  \includegraphics[width=.8\linewidth]{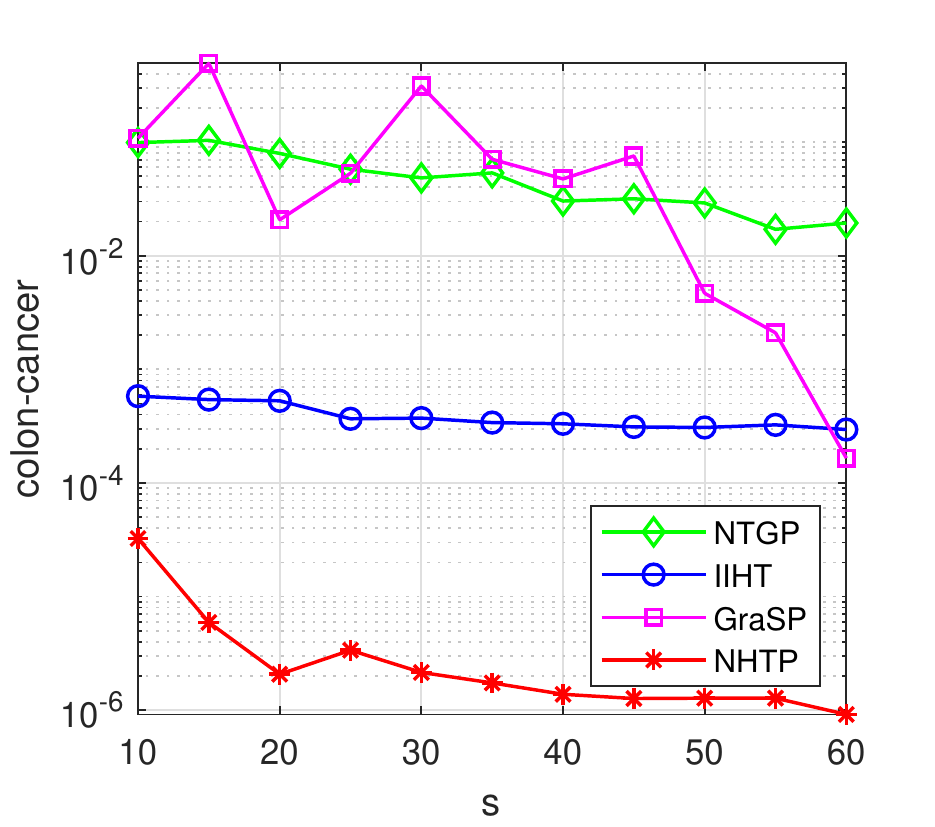}
\end{subfigure}
\begin{subfigure}{0.43\textwidth}
\caption{CPU time}\vspace{-2mm}
  \includegraphics[width=.8\linewidth]{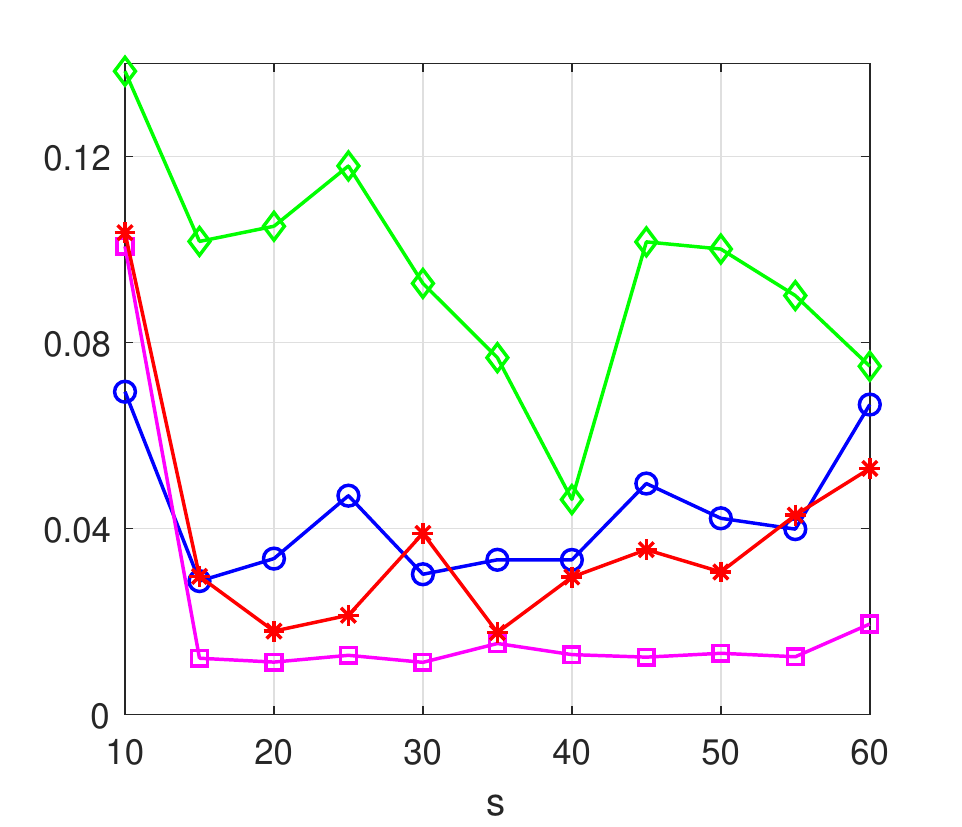}
\end{subfigure}\\\vspace{3mm}
\begin{subfigure}{0.43\textwidth}
\caption{$\ell(\bx)$}\vspace{-2mm}
  \includegraphics[width=.8\linewidth]{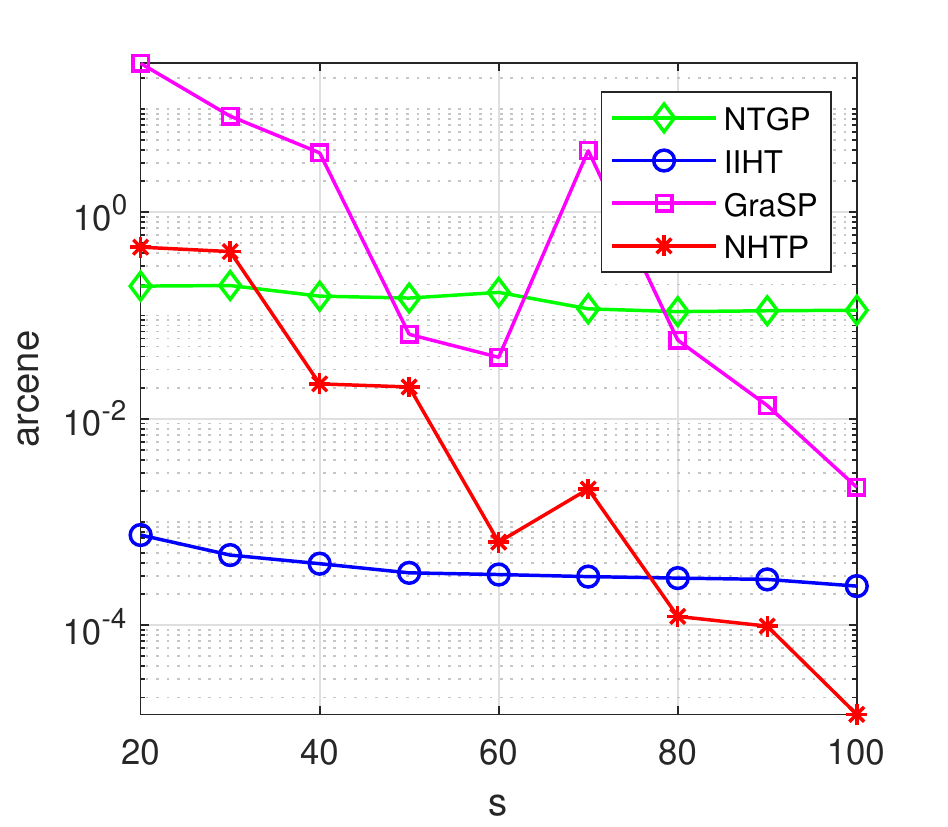}
\end{subfigure}
\begin{subfigure}{0.43\textwidth}
 \caption{CPU time}\vspace{-2mm}
 \includegraphics[width=.8\linewidth]{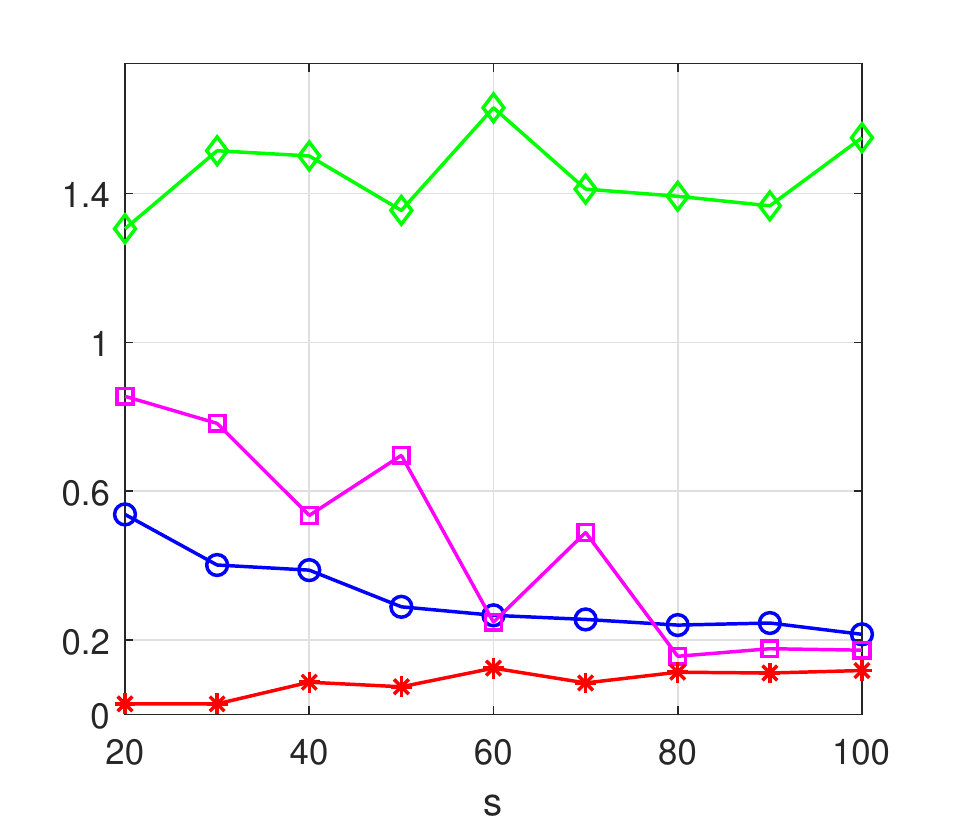}
\end{subfigure}\\\vspace{3mm}
\begin{subfigure}{0.43\textwidth}
 \caption{$\ell(\bx)$}\vspace{-2mm}
 \includegraphics[width=.8\linewidth]{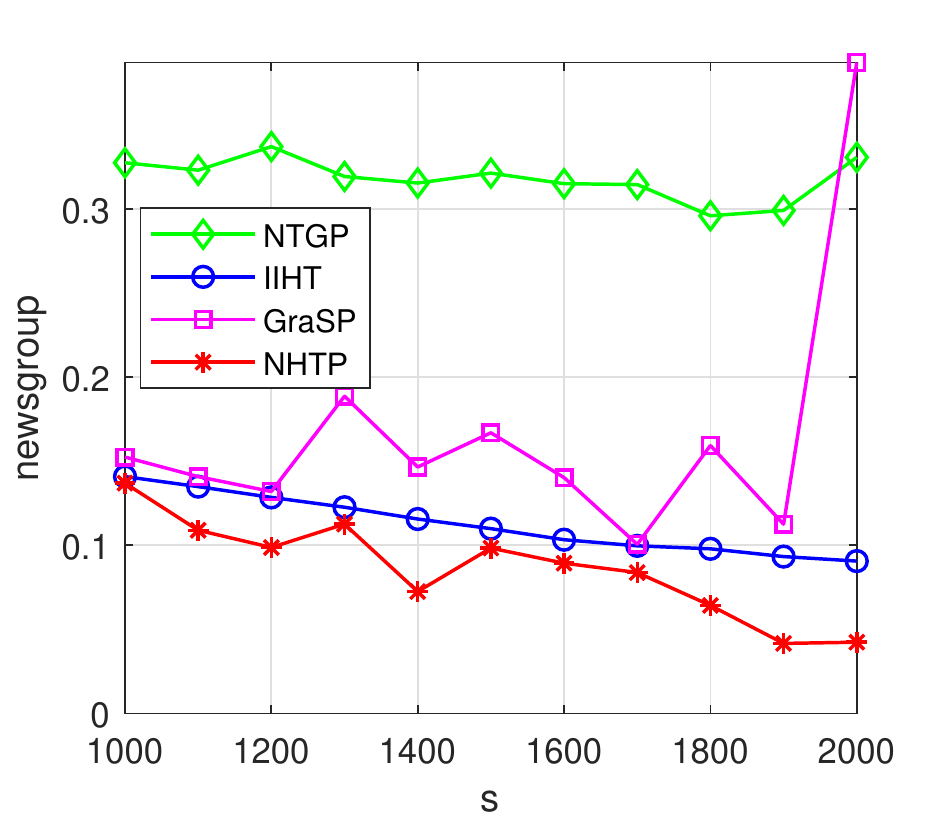}
\end{subfigure}
\begin{subfigure}{0.43\textwidth}
\caption{CPU time}\vspace{-2mm}
  \includegraphics[width=.8\linewidth]{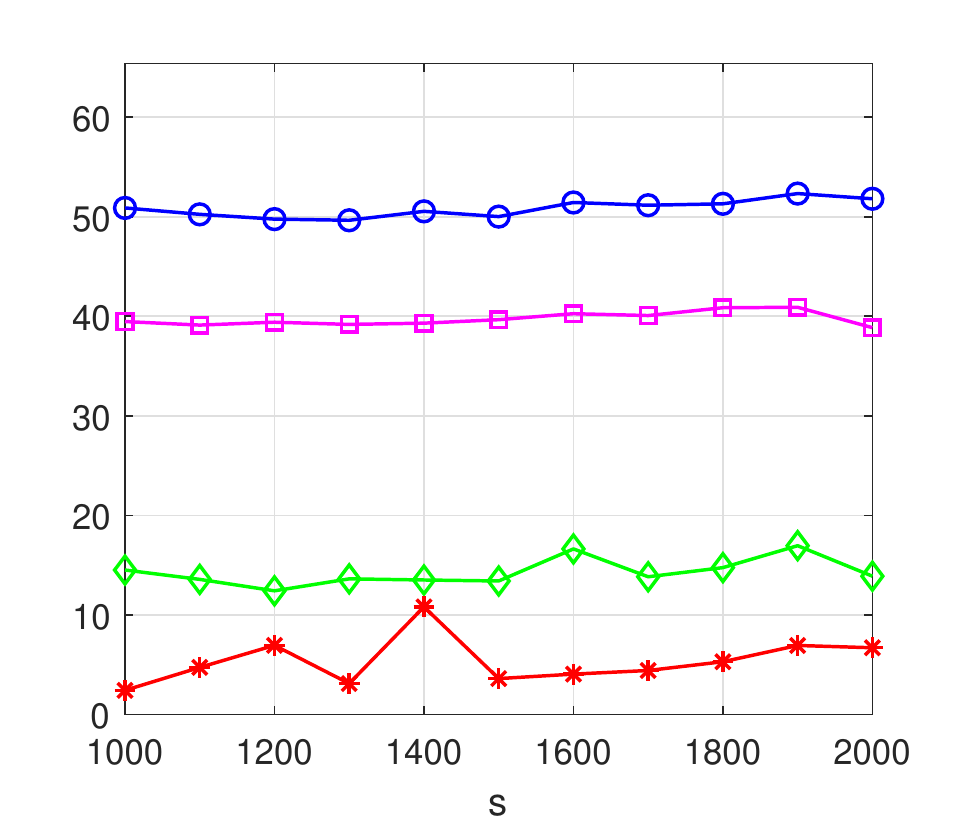}
\end{subfigure}\\\vspace{3mm}
\begin{subfigure}{0.43\textwidth}
\caption{$\ell(\bx)$}\vspace{-2mm}
  \includegraphics[width=.8\linewidth]{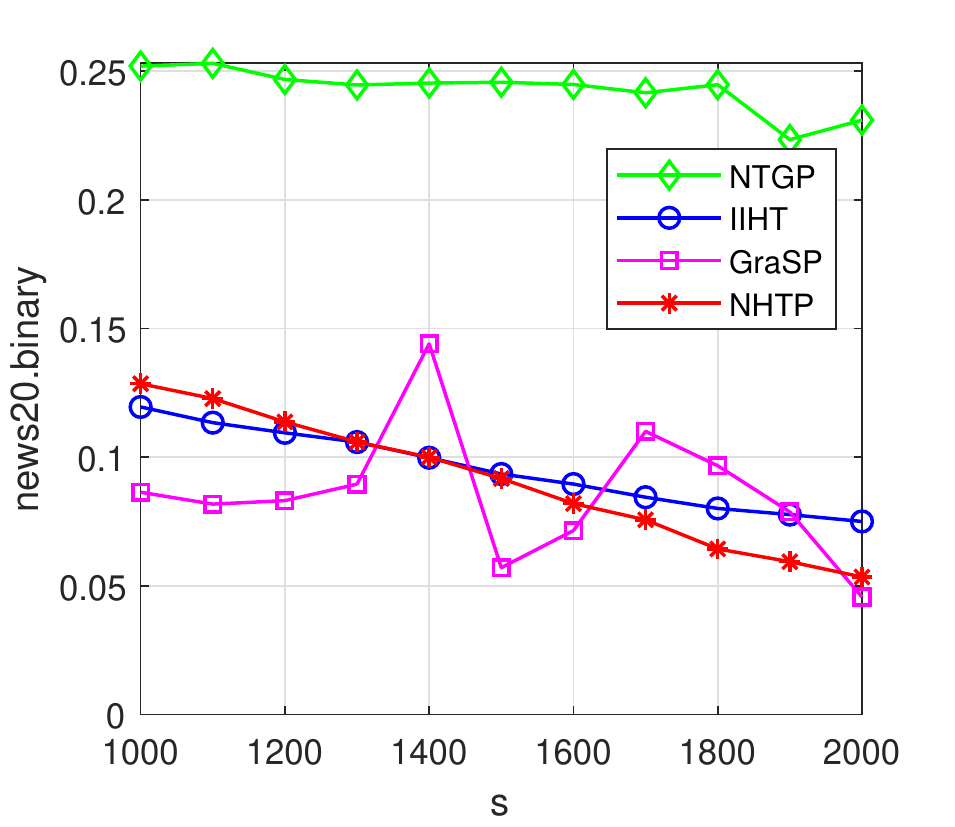}
\end{subfigure}
\begin{subfigure}{0.43\textwidth}
\caption{CPU time}\vspace{-2mm}
  \includegraphics[width=.8\linewidth]{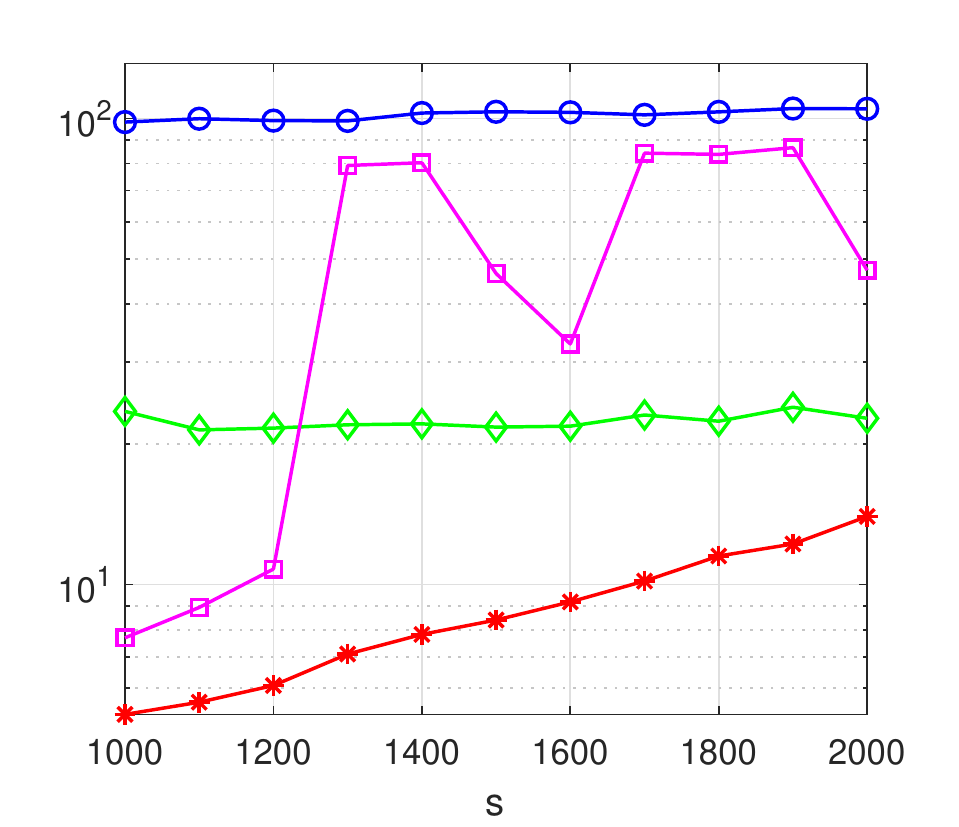}
\end{subfigure}
\caption{Logistic loss $\ell(\bx)$ and CPU time of four methods for Example \ref{log-EX3}.}
\label{fig:log-ex3-s}
\end{figure}

Now we compare these four methods on solving real data in Example \ref{log-EX3}. 
For each method, we demonstrate its performance on instances with varying $s$.
We first illustrate the performance of each method on solving those data without testing data sets.  As presented in Fig. \ref{fig:log-ex3-s}, we have the following observations:

 \begin{itemize}
 	
\item For \texttt{colon-cancer},  \NHTP\ obtained the smallest $\ell(\bx)$ followed by {\tt IIHT}. While {\tt GraSP} ran the fastest and {\tt NTGP} performed the slowest.

\item For \texttt{arcene},  {\tt IIHT} and \NHTP\ generated best $\ell(\bx)$ when $s<80$ and $s\geq 80$ respectively. And the latter consumed the smallest CPU time.

\item For \texttt{newsgroup},  \NHTP\ outperformed others in terms of the smallest $\ell(\bx)$ and CPU time. {\tt NTGP} rendered the worst logistic loss and {\tt IIHT} ran the slowest.

\item For \texttt{news20.binary},  {\tt GraSP} performed unstably, yet achieving best $\ell(\bx)$ 
for some cases such as $s\leq1300$.  
{\tt NTGP} still produced the highest logistic loss. As for computational speed, \NHTP\ was the fastest and {\tt IIHT} was the slowest.

 \end{itemize}

\begin{figure}[H]
\centering
\begin{subfigure}{0.32\textwidth}
\caption{Training $\ell(\bx)$}\vspace{-2mm}
  \includegraphics[width=1\linewidth]{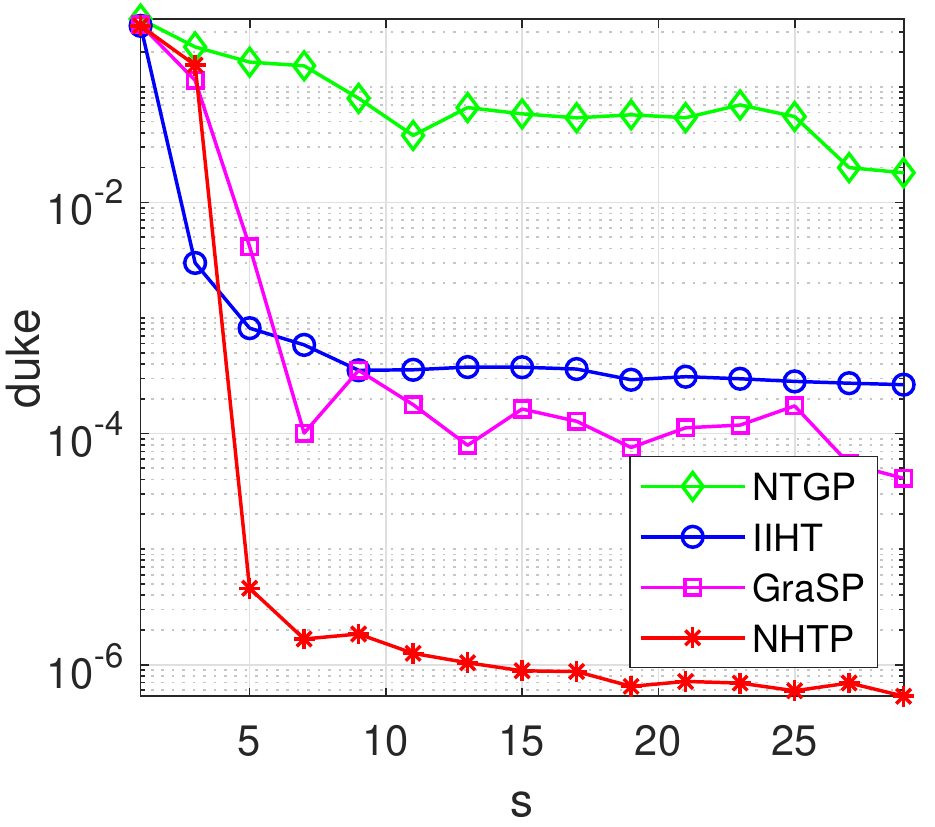}
\end{subfigure}
\begin{subfigure}{0.3\textwidth}
 \caption{Testing $\ell(\bx)$}\vspace{-2mm}
 \includegraphics[width=1\linewidth]{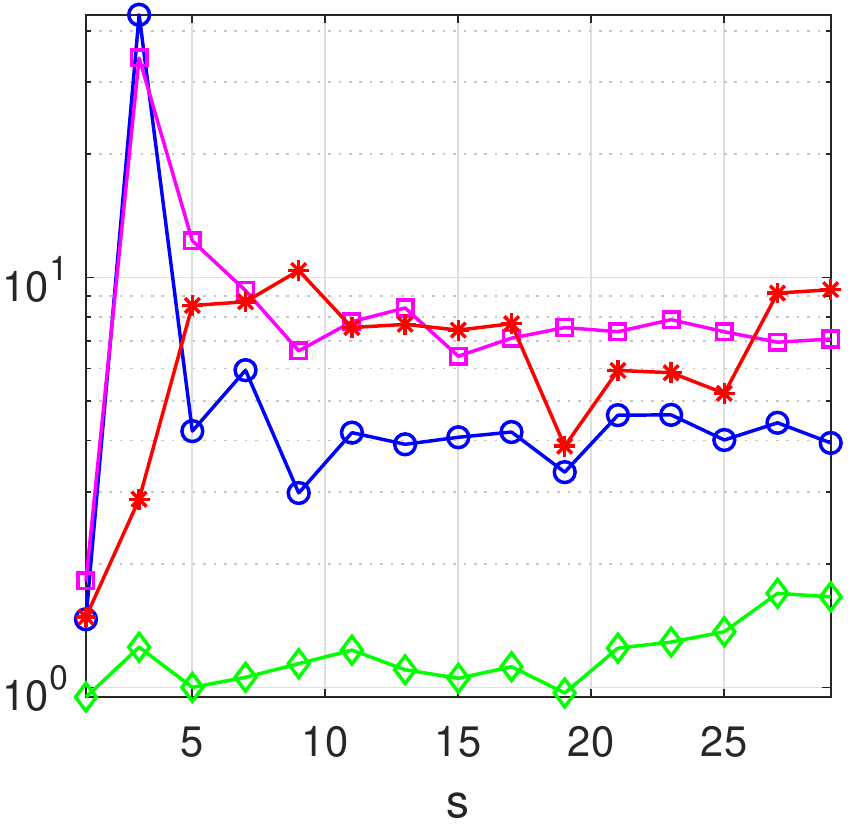}
\end{subfigure}
\begin{subfigure}{0.3\textwidth}
\caption{CPU time}\vspace{-2mm}
  \includegraphics[width=1\linewidth]{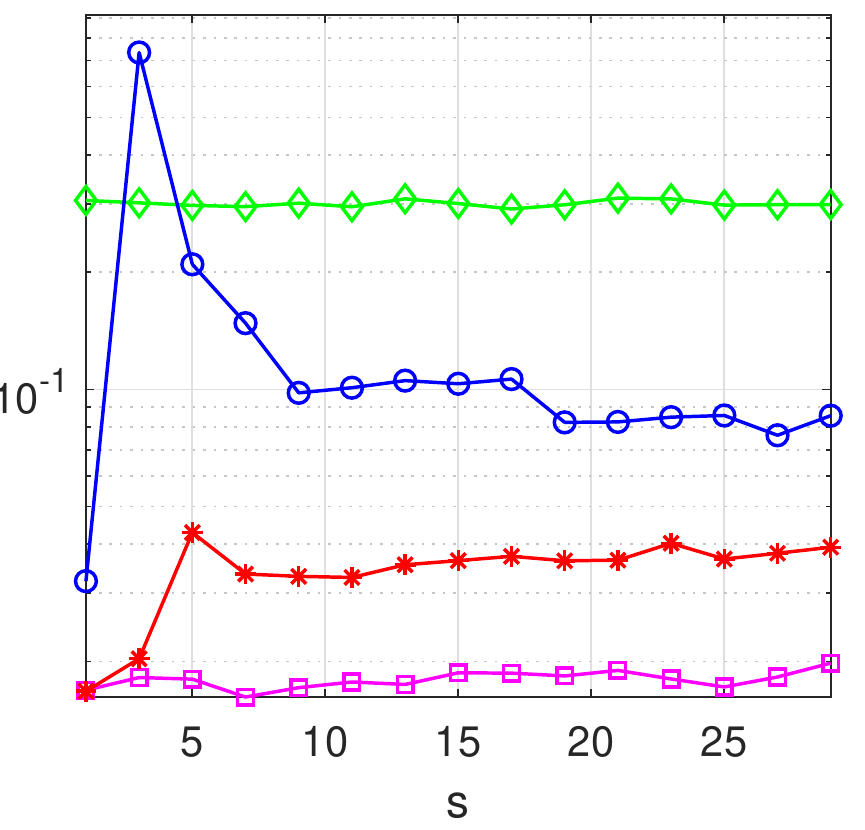}
\end{subfigure}\\
\begin{subfigure}{0.32\textwidth}
\caption{Training $\ell(\bx)$}\vspace{-2mm}
  \includegraphics[width=1\linewidth]{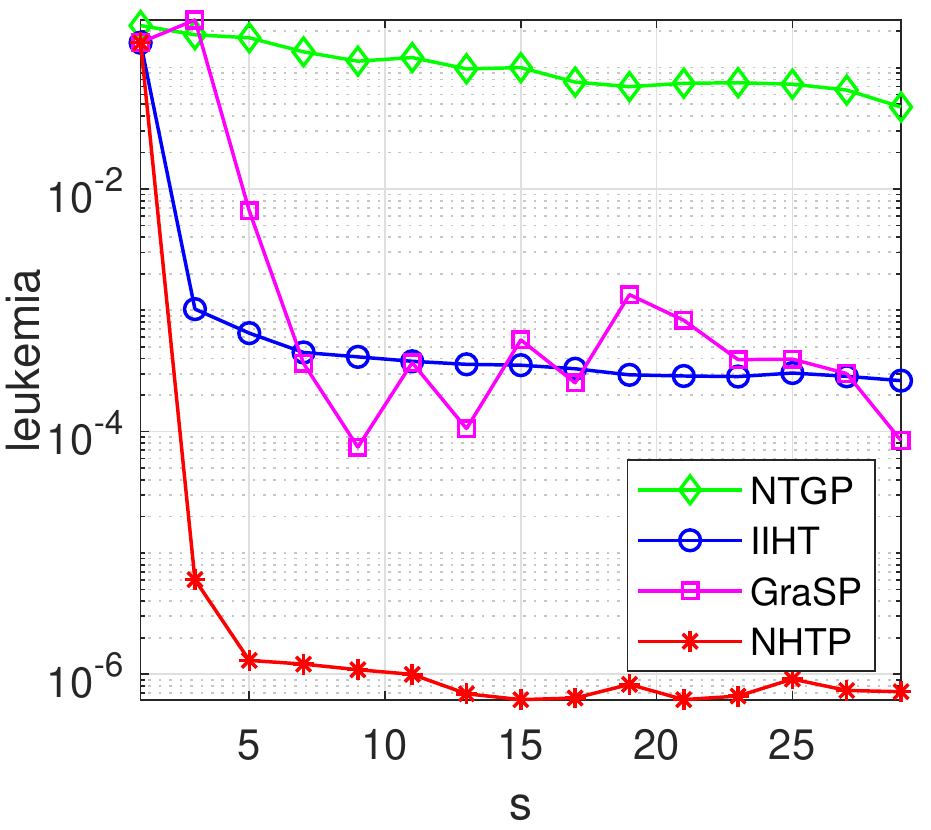}
\end{subfigure}
\begin{subfigure}{0.3\textwidth}
\caption{Testing $\ell(\bx)$}\vspace{-2mm}
  \includegraphics[width=1\linewidth]{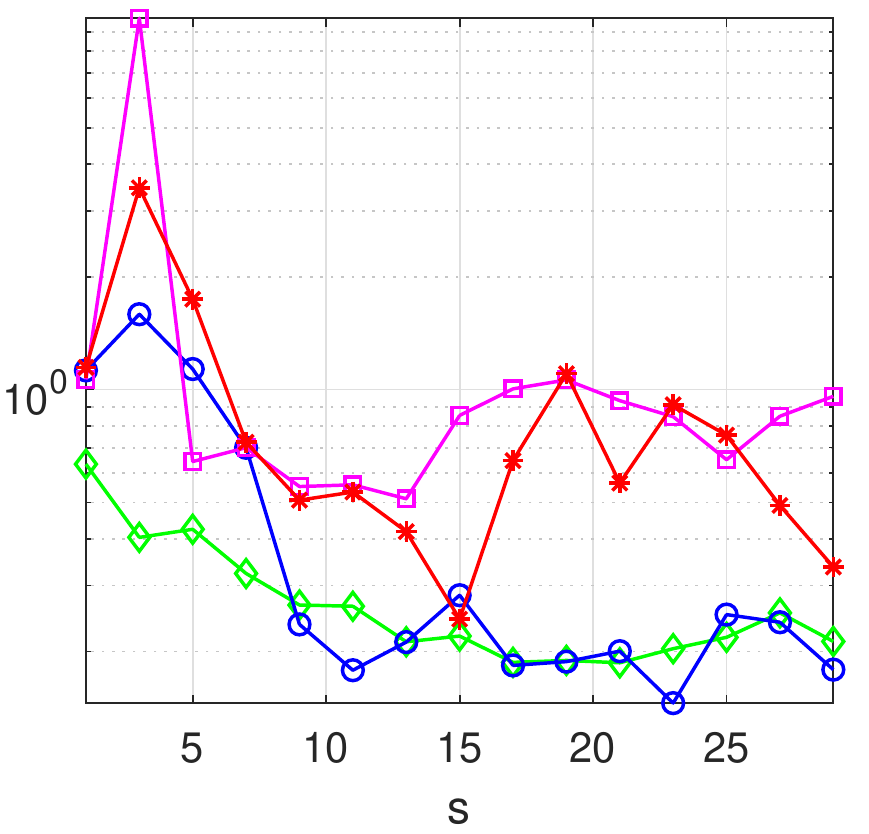}
\end{subfigure}
\begin{subfigure}{0.3\textwidth}
\caption{CPU time}\vspace{-2mm}
  \includegraphics[width=1\linewidth]{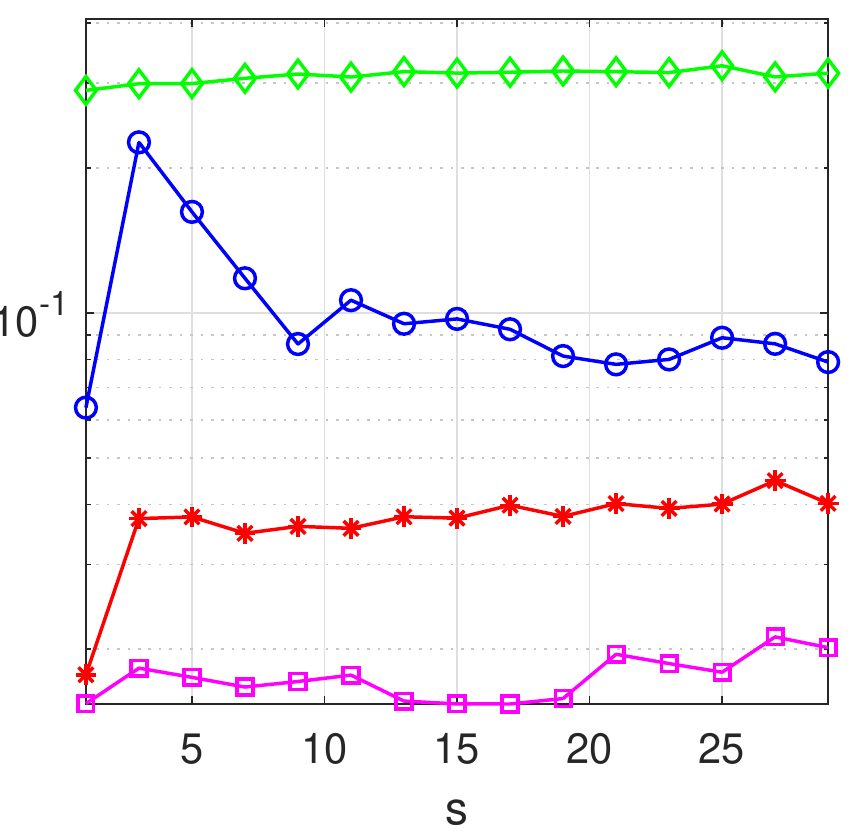}
\end{subfigure}\\
\begin{subfigure}{0.32\textwidth}
\caption{Training $\ell(\bx)$}\vspace{-2mm}
  \includegraphics[width=1\linewidth]{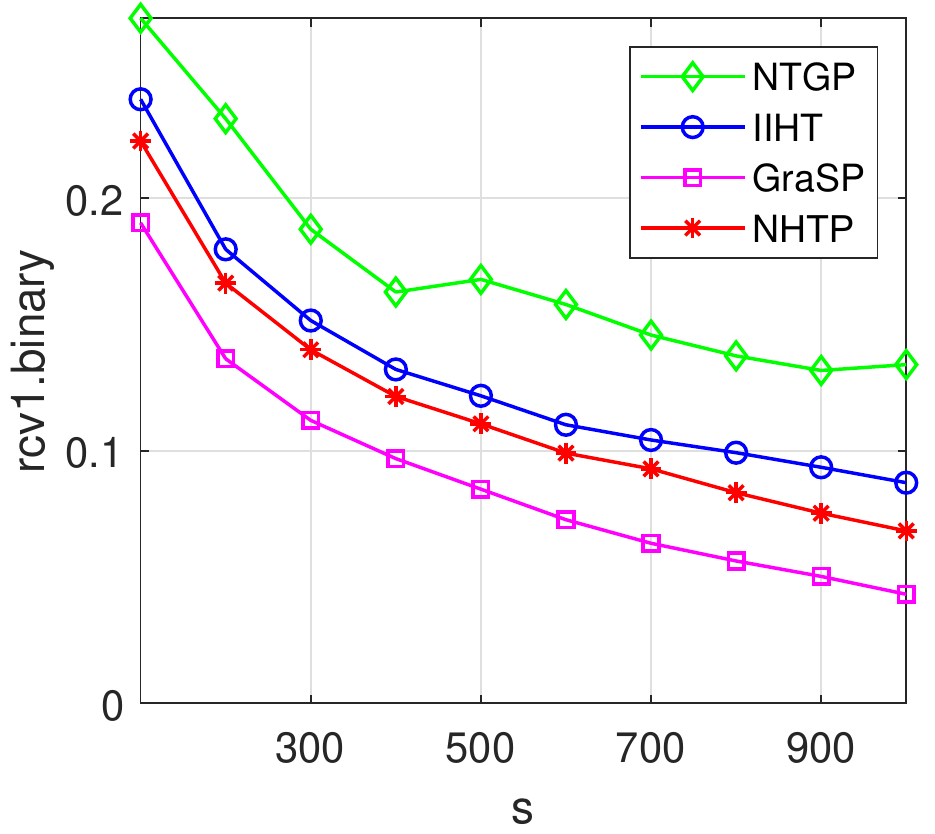}
\end{subfigure}
\begin{subfigure}{0.3\textwidth}
\caption{Testing $\ell(\bx)$}\vspace{-2mm}
  \includegraphics[width=1\linewidth]{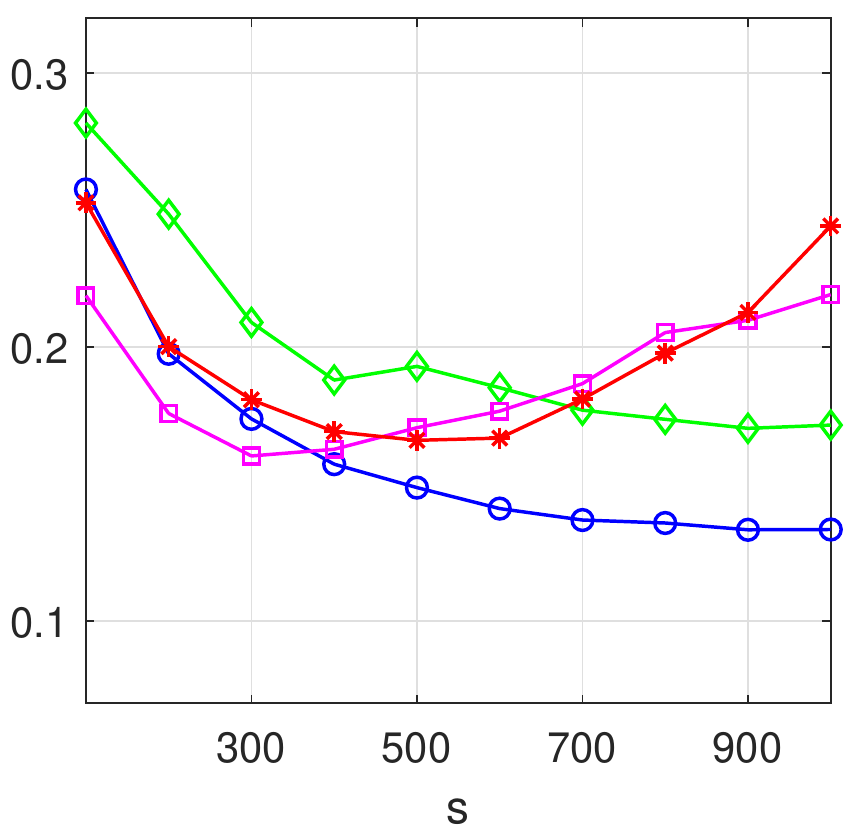}
\end{subfigure}
\begin{subfigure}{0.3\textwidth}
\caption{CPU time}\vspace{-2mm}
  \includegraphics[width=1\linewidth]{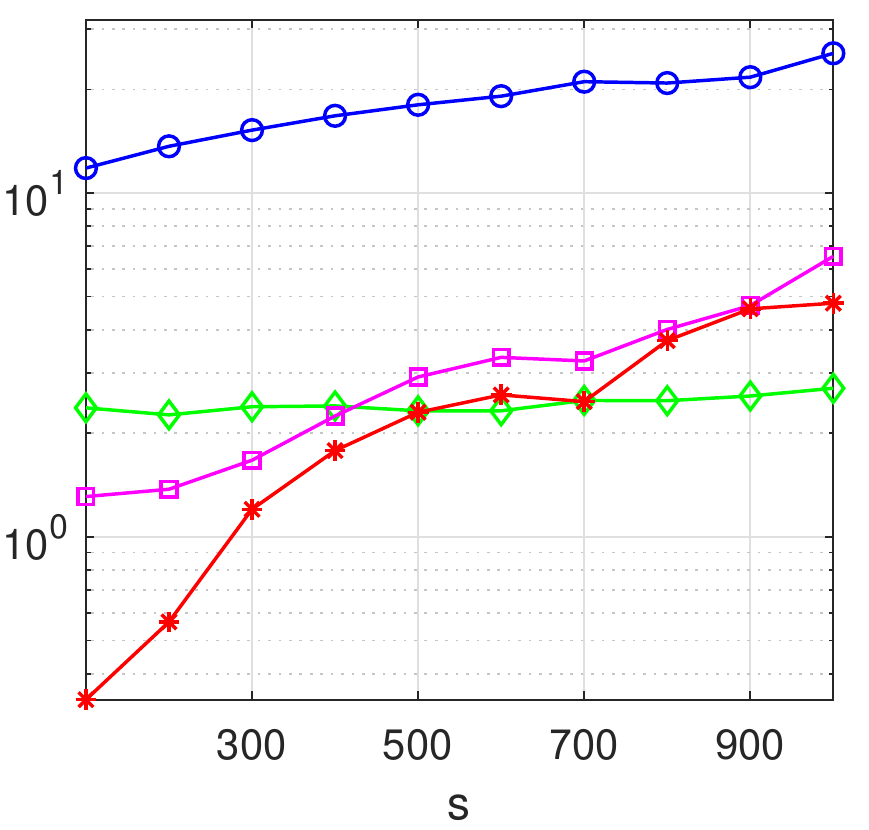}
\end{subfigure}
\caption{Logistic loss $\ell(\bx)$ and CPU time of four methods for Example \ref{log-EX3}.}
\label{fig:log-ex3-s-test}
\end{figure}

Next we illustrate the performance of each method on solving those data with testing data sets. As shown in Fig. \ref{fig:log-ex3-s-test}, some comments are able to be made as follows:

\begin{itemize}
	
	\item For \texttt{duke breast-cancer}, along with increasing $s$, $\ell(\bx)$ on training data obtained
	by \NHTP\ dropped significantly, with order $10^{-6}$. By contrast, {\tt NTGP} stabilized at above  $10^{-2}$. When it comes to the testing data, apparently {\tt NTGP} yielded the best $\ell(\bx)$, followed by {\tt IIHT}. It seems that the higher $\ell(\bx)$ on training data  was solved by a method, the lower $\ell(\bx)$ on testing data would be provided. For CPU time, {\tt GraSP} behaved the fastest, followed by  \NHTP, {\tt IIHT} and {\tt NTGP}.
	
	\item For  \texttt{leukemia}, the performance of each method was similar to that on  \texttt{duke breast-} \texttt{cancer} data. A slightly difference  was that {\tt NTGP}  no more offered the best $\ell(\bx)$ on testing data as {\tt IIHT} generated the best ones for some $s$.
	
	\item For  \texttt{rcv1.binary}, {\tt GraSP} performed the best  $\ell(\bx)$ on training data, followed by our method. Again  {\tt NTGP} came the last. It is obvious that {\tt IIHT} got the smallest $\ell(\bx)$ on testing data when $s\geq400$, while {\tt GraSP} produced the best ones otherwise. For CPU time,  \NHTP\ and {\tt NTGP} was the most efficient when $s\geq600$ and $s>600$ respectively.
	
\end{itemize} 

\section{Conclusion and Future Research} \label{Section-Conclusion}

There exists numerous papers that use a restricted Newton step to accelerate methods belonging
to hard-thresholding pursuits.
This results in the method of Newton hard-thresholding pursuit.
On the one hand, existing empirical experience shows significance acceleration when Newton's step is
employed .
On the other hand, existing theory for such methods does not offer any better statistical guarantee than
the simple hard thresholding counterparts.
The discrepancy between the superior empirical performance and the no-better theoretical guarantee 
has been well documented in the case of CS problem (\ref{CS}) and it invites further theory for 
justification.

In this paper, we develop a new \NHTP, which makes use of the strategy 
``{\em approximation and restriction}'' to obtain the truncated approximation within a subspace.
This is in contrast to the popular strategy ``{\em restriction and approximation}''.
We note that both strategies lead to the same Newton step in the case of CS.
We further cast the resulting Newton step as a Newton iteration for a nonlinear equation.
This new interpretation of the Newton step provides a new route for establishing its quadratic
convergence. Finally, we used the Armijo line search to globalize the method.
Extensive numerical experiments confirm the efficiency of the proposed method.
The global and quadratic convergence theory for \NHTP\ offers a theoretical justification 
why such methods are more efficient than their simple hard thresholding counterparts.
There are a few of topics that are worth exploring further.

\begin{itemize}

{{
\item[(i)] We expect that our algorithmic framework will make it possible to study quadratic convergence
of existing \NHTP\ based on the strategy of ``{\em restriction and approximation}''.
A plausible approach would be to regard such method as an inexact version of our \NHTP.
Technically, it would involve quantifying/controlling the inexactness so as to ensure the 
quadratic convergence to hold. 

\item[(ii)] As rightly pointed out by one referee, ``{\em the proof technique revolves 
	around providing sufficient conditions for descent (and reverting to standard gradient descent when descent does not hold). However, there are many methods that are non-descent and convergence still does hold. Blumensath's method
	in \cite{Blu12} (or acceleration in general) is one such method wherein it has been observed in practice in subsequent works that the existence of a `ripple' effect, akin to other accelerated methods wherein descent is not required for overall convergence.}'' 
In our numerical experiments, we also observed that the provided sufficient conditions for descent are not necessary. 
But it will be curious to see if  the descent itself is necessary, while enjoying the stated convergence. 

\item[(iii)] We proved in Thm.~\ref{Thm-Qudratic-Convergence} that quadratic convergence takes place after
certain $k_0$ iterations. It would be nice to estimate and quantify how big this $k_0$ would be.
Such research belongs to computational complexity in optimization.
We plan to investigate all of those in future.
}}
\end{itemize} 

\noindent
{\bf Acknowledgement.}
We sincerely thank the associate editor and the  two referees for their detailed comments that have helped us to 
improve the paper. We particularly thank the referee who went through our technical proofs and offered us valuable 
suggestions on the condition (\ref{Extra-Condition}). 
We also thank Prof Ziyan Luo of Beijing Jiaotong University, who helped us to improve the proof of Thm.~\ref{Thm-Qudratic-Convergence}.

\appendix
\section{Identities and Inequalities for Proofs}

Due to the restricted fashion of \NHTP, we need to keep tracking the indices belonging to the subspace
$\bx|_{T_k} = 0 $ and also those fall out of this subspace. To simplify our proofs, we will use a few more
abbreviations and derive some identities and inequalities associated with the Newton direction
$\bd^k_N$. The sequence $\{ \bx^k\}$ used is generated by \NHTP.
\\

\noindent{\bf (a) Simplification of Newton's equation (\ref{Newton-Direction}).}  We first define
 \begin{eqnarray}\label{HTT}
J_{k} := T_{k-1}\setminus T_{k}, \qquad
H_k:=\nabla_{{T_{k}}}^{2} f (\bx^k), \qquad
G_k :=\nabla_{{T_{k}},{J_k}}^{2} f (\bx^k).
\end{eqnarray}
We also have the following easy observation:
\begin{equation} \label{Index-T}
\support(\bx^k) \subseteq T_{k-1}, \quad
| T_{k} | = | T_{k-1} | = s, \quad \mbox{and} \quad
| T_k \setminus T_{k-1} | = | T_{k-1} \setminus T_k | = | J_k|.
\end{equation}
It is important to note that $\bx^k_{T_k^c}$ is also $s$-sparse.
This is because for any $i \not\in T_{k-1}$,
$x^k_i = 0$ (because $\support(\bx^k) \subseteq T_{k-1}$),
\begin{equation} \label{xktc}
\bx^k_{T_k^c} = \left[ \begin{array}{c}
 \bx^k_{ T_k^c \cap T_{k-1}} \\ [0.6ex]
 0
\end{array} \right]
= \left[ \begin{array}{c}
\bx^k_{ T_{k-1} \setminus T_k} \\ [0.6ex]
0
\end{array} \right] = \left[ \begin{array}{c}
\bx^k_{J_k} \\ [0.6ex]
0
\end{array} \right],
\end{equation}
and $| J_k | \le | T_{k-1} \setminus T_k | \le | T_{k-1} | = s$.
We emphasize that $J_k$ captures all nonzero elements in $\bx^{k}_{T_k^c}$.
Therefore, we will see more $J_k$ instead of $T_k^c$ being used in our derivation below.
This observation leads to the simplified Newton  equation of (\ref{Newton-Direction}):
\begin{equation}\label{Newton-Direction-2}
\left\{
\begin{array}{rcl}
H_k (\bd_N^k)_{T_{k}} &=& G_k \bx^{k}_{J_k} - \nabla_{ T_{k}}  f (\bx^k) \\ [1.5ex]
(\bd_N^{k})_{ T^c_{k}}&=&-\bx^{k}_{{ T }_{k}^c} = - \left[ \begin{array}{c}
\bx^k_{J_k} \\ [0.6ex]
0
\end{array} \right].
\end{array}
\right.
\end{equation}
An important feature to note is that the vectors $(\bd_N^k)_{T_{k}}$, $(\bd_N^{k})_{ T^c_{k}}$,
$\bx^{k}_{J_k}$ are all $s$-sparse. Putting together, at each iteration, we only involve vectors
that do not exceed $2s$-sparsity. This is the reason why our assumptions are always on
$2s$-restricted properties of $f$.\\

\noindent {\bf (b) An identity on the Newton direction.} This involves a string of equalities as
follows. We write $\bd^k$ for $\bd^k_N$ because there is no danger to cause any confusion.

\begin{eqnarray*}
&& \langle \bd^{k}_{T_{k}\cup J_{k} },  \nabla^2_{T_{k}\cup J_{k} }
    f(\bx^k)\bd^{k}_{T_{k}\cup J_{k} } \rangle
    \qquad \quad (\mbox{note} \ {T_{k}\cap J_{k} } = \emptyset)
    \\ [0.6ex]
&=& \left[
    \begin{array}{cc}
       \bd^{k}_{T_{k}  }  \\ [0.6ex]
       \bd^{k}_{  J_{k} }
    \end{array}
\right] ^\top \left[
    \begin{array}{cc}
       H_k, \ & \ G_k  \\ [0.6ex]
       G_k ^\top,\ & \ \nabla^2_{ J_{k} }f(\bx^k)
    \end{array}
\right]  \left[
    \begin{array}{cc}
       \bd^{k}_{T_{k}  }  \\ [0.6ex]
       \bd^{k}_{  J_{k} }
    \end{array}
\right] \\ [0.6ex]
&=&  \left[
    \begin{array}{cc}
       \bd^{k}_{T_{k}  }  \\ [0.6ex]
       \bd^{k}_{  J_{k} }
    \end{array}
\right] ^\top \left[
    \begin{array}{c}
      H_k\bd^{k}_{T_{k} }+ G_k \bd^{k}_{  J_{k} } \\ [0.6ex]
      G_k ^\top  \bd^{k}_{T_{k}  }+   \nabla^2_{ J_{k} }f(\bx^k)  \bd^{k}_{  J_{k} }
    \end{array} \right]
                        \\ [0.6ex]
&\overset{(\ref{Newton-Direction-2})}{=}& \left[
    \begin{array}{cc}
       \bd^{k}_{T_{k}  }  \\ [0.6ex]
       \bd^{k}_{  J_{k} }
    \end{array}
\right]^\top \left[
\begin{array}{c}
-\nabla_{T_{k} } f(\bx^{k}) \\ [0.6ex]
 G_k ^\top  \bd^{k}_{T_{k}  }+   \nabla^2_{ J_{k} }f(\bx^k)  \bd^{k}_{  J_{k} }
\end{array}
\right]  \nonumber\\ [0.6ex]
&=&-\langle \nabla_{T_{k} } f(\bx^{k}),  \bd^{k}_{T_{k} }\rangle+
\langle G_k \bd^{k}_{  J_{k} }, \bd^{k}_{  T_{k}}\rangle+
\langle \bd^{k}_{ J_{k} },   \nabla^2_{ J_{k} }f(\bx^k )\bd^{k}_{ J_{k} }\rangle\nonumber\\ [0.6ex]
&\overset{(\ref{Newton-Direction-2})}{=}&-\langle \nabla_{T_{k} } f(\bx^{k}),  \bd^{k}_{T_{k} }\rangle-
\langle H_k \bd^{k}_{ T_{k}}+\nabla_{ T_{k}}  f (\bx^k), \bd^{k}_{  T_{k}}\rangle+
\langle \bd^{k}_{ J_{k} },   \nabla^2_{ J_{k} }f(\bx^k )\bd^{k}_{ J_{k} }\rangle\nonumber\\ [0.6ex]
&=&-2\langle \nabla_{T_{k} } f(\bx^{k}),  \bd^{k}_{T_{k} }\rangle-
\langle H_k \bd^{k}_{ T_{k}}, \bd^{k}_{  T_{k}}\rangle+
\langle \bd^{k}_{ J_{k} },   \nabla^2_{ J_{k} }f(\bx^k )\bd^{k}_{ J_{k} }\rangle .\nonumber
\end{eqnarray*}
This leads to our identity:
\begin{eqnarray}
\label{bounded-H}
2\langle \nabla_{T_{k} } f(\bx^{k}),\  \bd^{k}_{T_{k} }\rangle   &=& -\langle \bd^{k}_{T_{k}\cup J_{k} }, \nabla^2_{T_{k}\cup J_{k} } f(\bx^k)\bd^{k}_{T_{k}\cup J_{k} }\rangle\nonumber  \\ [0.6ex]
&&-\langle H_k \bd^{k}_{ T_{k}},\ \bd^{k}_{  T_{k}}\rangle+
\langle \bd^{k}_{ J_{k} }, \  \nabla^2_{ J_{k}  }f(\bx^k )\bd^{k}_{ J_{k} }\rangle.
\end{eqnarray}

\noindent
{\bf (c) An inequality on the gradient sequence.}
The role of $T_k$ is like a working active set that is designed to identify the true support of
an optimal solution.
Its complementary set $T_k^c$ is handled in such a way to make sure the next iterate $\bx^{k+1}$ has
zeros on $T_k^c$. To achieve this, in both the Newton direction $\bd^k_N$ and the gradient direction
$\bd^k_g$ we set
\[
  \Big( \bd^k_N \Big)_{T_k^c} = \Big( \bd^k_g \Big)_{T_k^c} = - \bx^k_{T_k^c}.
\]
Let $\bd^k$ be either $\bd_N^k$ or $\bd^k_g$. It follows from (\ref{xktc}) that
\begin{equation}\label{xk-Tc-T1}
\begin{array}{l}
\|\bx^{k}_{T^c_{k}}\| = \|\bx^{k}_{ J_{k} }\|=\|\bd^{k}_{ J_k}\| = \|\bd^{k}_{T^c_{k}}\|,
\quad
\| \bd^k \|=\| \bd^{k}_{T_{k}\cup J_{k} } \|\\ [2ex]
\langle \nabla_{ T^c_{k} } f(\bx^{k}), \ \bx^{k}_{T^c_{k} }\rangle
= \langle \nabla_{J_{k} } f(\bx^{k}),  \ \bx^{k}_{ J_{k} }\rangle
\end{array}
\end{equation}
By the definition of $T_k$ and the fact, $x^k_i = 0$ for $i \in T_k \setminus T_{k-1}$, we have
\[
|\eta \nabla_i f(\bx^{k})|^2=|x_i^k-\eta \nabla_i f(\bx^{k})|^2
\geq |x_j^k-\eta \nabla_j f(\bx^{k})|^2, \ \ \forall \ i \in T_k \setminus T_{k-1}, \ \ j\in J_{k}.
\]
The above inequality and the fact $|T_k \setminus T_{k-1}| = | J_k|$ in (\ref{Index-T}) imply
\begin{eqnarray*}
	&& \eta^2 \| \nabla_{T_{k}\setminus T_{k-1}} f(\bx^{k})\|^2
	   =  \sum_{i \in T_k \setminus T_{k-1}} | \eta \nabla_i f(\bx^k) |^2
	  \ge \sum_{j \in J_k} | x^k_j - \eta \nabla_i f(\bx^k) |^2 \\ [0.6ex]
	&\geq& \|\bx^k_{J_{k}}-\eta \nabla_{J_{k}} f(\bx^{k})\|^2
	  =  \|\bx^k_{J_{k}}\|^2-2\eta  \langle\bx^k_{J_{k}},\ \nabla_{J_{k}} f(\bx^{k})\rangle+ \eta^2\| \nabla_{J_{k}} f(\bx^{k})\|^2\\
	&\overset{(\ref{xk-Tc-T1})}{=}&  \|\bx^k_{T^c_{k}}\|^2-2\eta \langle\bx^k_{J_{k}},\ \nabla_{J_{k}} f(\bx^{k})\rangle+ \eta^2\| \nabla_{J_{k}} f(\bx^{k})\|^2,
\end{eqnarray*}
which together with  \begin{eqnarray*}
	\| \nabla_{T_{k}} f(\bx^{k})\|^2&=&\| \nabla_{T_{k} \cap T_{k-1}} f(\bx^{k})\|^2+\| \nabla_{T_{k} \setminus T_{k-1}} f(\bx^{k})\|^2\\
	\| \nabla_{T_{k-1}} f(\bx^{k})\|^2&=&\| \nabla_{T_{k} \cap T_{k-1}} f(\bx^{k})\|^2+\| \nabla_{J_{k}} f(\bx^{k})\|^2\end{eqnarray*}
results in the following inequality on the gradient $\nabla f(\bx^k)$
\begin{eqnarray}\label{existence-alpha-facts-3}
\eta  \| \nabla_{T_{k}} f(\bx^{k})\|^2 - \eta  \| \nabla_{T_{k-1}} f(\bx^{k})\|^2-
\|\bx^k_{ T_k^c}\|^2/\eta\geq - 2\langle \bx^k_{J_{k}},\ \nabla_{J_{k}} f(\bx^{k})\rangle.
\end{eqnarray}

\section{Proofs of all results}

\subsection{Proof of Lemma~\ref{Lemma-Feta}}

\begin{proof} {{The first claim is obvious and we onely prove the second one.}}
The proof for the ``only if'' part is straightforward. Suppose $\bx$ satisfies (\ref{Fix-Point-Eq}).
We have $\bx = \P_s(\bx - \eta \nabla f(\bx))$. 
By the definition of $\P_s(\cdot)$ and $T \in \T(\bx, \eta)$, we have $\bx_{T^c}=0$ and
\[
  \bx_T = \Big(  \P_s(\bx - \eta \nabla f(\bx)) \Big)_T = (\bx - \eta \nabla f(\bx))_T 
  = \bx_T - \eta \nabla_Tf(\bx), 
\]
which implies $\nabla_T f(\bx) =0$.

We now prove the ``if'' part. 
Suppose we have $F_{\eta }(\bx;T)=0~~ \text{for~all}~~T\in \T (\bx;\eta)$, namely,
\begin{equation}\label{zdt}
  \nabla_{ T}  f (\bx)=0, ~~  \bx_{  T^c }=0.
\end{equation}
We consider two cases. 
Case I: $  \T (\bx;\eta )$ is a singleton.
By letting $T$ {{ be the only element of}} $  \T (\bx;\eta )$, then
$$
\bx-\P_{s}(\bx-\eta   \nabla  f (\bx))=\left[\begin{array}{c}
\bx_{ T}\\
\bx_{   T^c }
\end{array}\right]
-\left[\begin{array}{c}
\bx_{ T}-\eta   \nabla_{ T}  f (\bx) \\
0
\end{array}\right]
\overset{(\ref{zdt})}{=}
\left[\begin{array}{c}
\bx_{ T}-\bx_{ T}\\
0-0
\end{array}\right]=0,$$
which means $\bx$ satisfies the fixed point equation (\ref{Fix-Point-Eq}). 

Case II: $  \T (\bx;\eta )$ has multiple elements.
Then by the definition (\ref{Tu}) of $\T(\bx;\eta )$ we have two claims:
$$(x-\eta   \nabla  f (\bx))_{(s)}=(x-\eta   \nabla  f (\bx))_{(s+1)}>0~~~~{\rm or}~~~~(x-\eta   \nabla  f (\bx))_{(s)}=0.$$
 Now we exclude the first claim. Without loss of any generality, we assume
$$|x_1-\eta  \nabla_1  f (\bx)|\geq\cdots\geq|x_s-\eta  \nabla_s  f (\bx)|=|x_{s+1}-\eta  \nabla_{s+1}  f (\bx)|=(x-\eta   \nabla  f (\bx))_{(s)}.$$
 Let $T_1=\{1,2,\cdots,s\}$ and $T_2=\{1,2,\cdots,s-1,s+1\}$. Then $F_{\eta }(\bx;T_1)=F_{\eta }(\bx;T_2)=0$ imply that $ \nabla _{T_1} f (\bx)= \nabla_{T_2}  f (\bx)=0$ and $\bx_{  T^c _1}=\bx_{  T^c _2}=0$, which lead to
\begin{eqnarray*}|x_1|&=&|x_1-\eta  \nabla_1  f (\bx)|\geq\cdots\geq|x_s|=|x_s-\eta  \nabla_s  f (\bx)|=\\
|x_{s+1}|&=&|x_{s+1}-\eta  \nabla_{s+1}  f (\bx)|=(x-\eta   \nabla  f (\bx))_{(s)}>0.\end{eqnarray*}
This is contradicted with $\bx_{  T^c _1}=0$ because of $(s+1)\in  T^c _1$. 
Therefore, we have $(x-\eta   \nabla  f (\bx))_{(s)}=0$.
 This together with the definition (\ref{Tu}) of $\T(\bx;\eta )$ yields $0=(x-\eta   \nabla  f (\bx))_{(s)}\geq |x_i-\eta  \nabla_i  f (\bx)|=|\eta   \nabla_i  f (\bx)|$ for any $i\in  T^c$, which combining $ \nabla_{T} f (\bx)=0$ renders $ \nabla  f (\bx)=0$. Hence $x_{(s)}=(x-\eta   \nabla  f (\bx))_{(s)}=0$, yielding $\|\bx\|_0<s$. 
Consequently, $\bx = \bx - \eta \nabla f(\bx)$ (because $\nabla f(\bx) =0$ and
$\bx = \P_s(\bx) = \P_s( \bx - \eta \nabla f(\bx)  )$ (because $\| \bx\|_0 < s$).
That is $\bx$ also satisfies the fixed point equation
(\ref{Fix-Point-Eq}).
\end{proof}

\subsection{Proof of Lemma~\ref{Lemma-Newton-Direction}}
\begin{proof} For  simplicity, we write $\bd^{k}:=\bd^{k}_N$. Since $f(\bx)$ is $m_{2s}$-restricted strongly convex and $M_{2s}$-restricted strongly smooth. 
For any $\|\bx\|_0\leq s$ , it follows from Definition \ref{def-RSCS} that
 \begin{equation}
m_{2s} I_{2s} \preceq \nabla^2_{{T}} f(\bx) \preceq M_{2s} I_{2s} \hspace{10mm} {\rm for~any~} |T|\leq 2s.\end{equation}
Clearly, $|T_{k}\cup J_{k}| \leq 2s$ due to $|T_{k}| \leq s$ and $|J_{k}| \leq s$.  This together with (\ref{bounded-H}) implies
      \begin{eqnarray}
2\left\langle \nabla_{T_{k} } f(\bx^{k}),  \bd^{k}_{T_{k} }\right\rangle
   &=& -\left\langle \bd^{k}_{T_{k}\cup J_{k} }, \nabla^2_{T_{k}\cup J_{k}  } f(\bx^k)\bd^{k}_{T_{k}\cup J_{k} }\right\rangle\nonumber  \\
 &&-\left\langle H_k \bd^{k}_{ T_{k}}, \bd^{k}_{  T_{k}}\right\rangle+
\left\langle \bd^{k}_{ J_{k} },   \nabla^2_{ J_{k} }f(\bx^k )\bd^{k}_{ J_{k} }\right\rangle \nonumber  \\
 & \leq &- m_{2s}\left[\|\bd^{k}_{T_{k}\cup J_{k} }\|^2+\|\bd^{k}_{T_{k} }\|^2\right] + M_{2s} \|\bx^{k}_{ T^c_{k} }\|^2 \nonumber\\
  &=&- m_{2s}\left[\|\bd^{k}_{T_{k}\cup J_{k} }\|^2+  \|\bd^{k}_{T_{k} }\|^2+  \|\bd^{k}_{ J_{k} }\|^2-  \|\bd^{k}_{ J_{k} }\|^2\right] + M_{2s} \|\bx^{k}_{ T^c_{k} }\|^2 \nonumber\\
   & \overset{(\ref{xk-Tc-T1})}{=}&- 2m_{2s}\|\bd^{k} \|^2 + m_{2s}\|\bx^{k}_{ T^c_{k} }\|^2+ M_{2s} \|\bx^{k}_{ T^c_{k} }\|^2 \nonumber\\
   & \leq& - 2m_{2s}\|\bd^{k} \|^2 + 2M_{2s} \|\bx^{k}_{ T^c_{k} }\|^2 \nonumber\\
\label{varify-cond}&\leq& -2\gamma \|\bd^k \|^2+  \|\bx^{k}_{ T^c_{k}}\|^2/(2\eta),
\end{eqnarray}
where   the last inequality  is owing to that $\gamma\leq m_{2s}$ and  $\eta \leq 1/(4M_{2s})$.\end{proof}

\subsection{Proof of Lemma~\ref{Lemma-dk}}

\begin{proof} 
It follows from the fact $\eta < \overline{ \eta }$ that	
$$
\eta<\overline{\eta }\leq \min \left\{ \frac{\gamma (\overline{\alpha } \beta ) }{ M_{2s}^2 },
 \  \overline{\alpha } \beta \right\} <  \min \left\{ \frac{\gamma   }{ M_{2s}^2 },
 \   1 \right\},
$$
where the last strict inequality used $\overline{ \alpha } \le 1$ and $\beta < 1$.
Therefore, $\rho$ is well defined and $\rho > 0$.
Since $\support(\bx^k) \subseteq T_{k-1}$, the relationships in (\ref{Index-T})-(\ref{existence-alpha-facts-3}) all hold.   
We now prove the claim by two cases.\\

\noindent
{\bf Case 1:} If $\bd^k=\bd^k_N$, then it follows  from (\ref{Search-Direction}) that
\begin{eqnarray}\label{fd-TT} 2\langle \nabla_{T_{k} } f(\bx^{k}),  \bd^{k}_{T_{k} }\rangle&\leq&-2\gamma \|\bd^{k} \|^2+ \|\bx^{k}_{ T^c_{k}  }\|^2/(2\eta ).
\end{eqnarray}
In addition,
\begin{eqnarray}\label{nabla-T-f-0}
\|\nabla_{ T_{k}}  f (\bx^k)\|^2 
&\overset{(\ref{Newton-Direction-2})}{=}& 
\|H_k  \bd^k_{T_{k}}-G_k  \bx^{k}_{ J_{k} } \|^2
\overset{(\ref{Newton-Direction-2})}{=} 
\|[H_k,\; G_k ] \bd^{k}_{T_{k}\cup J_{k} } \|^2 \\
\label{nabla-T-f}&\leq&M_{2s}^2 \| \bd^{k}_{T_{k}\cup J_{k} } \|^2\overset{(\ref{xk-Tc-T1})}{=}M_{2s}^2 \| \bd^k \|^2,
 \end{eqnarray}
where the inequality holds because $\|[H_k, \; G_k ]\|_2\leq \|\nabla^2_{T_{k}\cup J_{k} } f(\bx^k)\|_2$ due to $f$ being  $M_{2s}$-restricted strongly smooth and $|T_{k}\cup J_{k}|\leq 2s$. This together with (\ref{existence-alpha-facts-3}) derives
 \begin{eqnarray}\label{xk-T1-T1}
 - 2\langle \bx^k_{J_{k}},\nabla_{J_{k}} f(\bx^{k})\rangle \leq  \eta  M_{2s}^2 \| \bd^k\|^2-  \|\bx^k_{ T^c_{k} }\|^2/\eta -\eta   \| \nabla_{T_{k-1}} f(\bx^{k})\|^2.
 \end{eqnarray}
Direct calculation yields the following chain of inequalities,
\begin{eqnarray*}
   2\langle \nabla f(\bx^{k}), \bd^{k} \rangle
&=&  2\langle \nabla_{T_{k} } f(\bx^{k}),  \bd^{k}_{T_{k} }\rangle- 2\langle \nabla_{ T^c_{k}  } f(\bx^{k}),  \bx^{k}_{ T^c_{k}  }\rangle\\
 &\overset{(\ref{xk-Tc-T1})}{=}&   2 \langle \nabla_{T_{k} } f(\bx^{k}),  \bd^{k}_{T_{k} }\rangle- 2\langle \nabla_{J_{k} } f(\bx^{k}),  \bx^{k}_{ J_{k} }\rangle\\
&\overset{(\ref{fd-TT}, \ref{xk-T1-T1})}{\le}&-\left[2\gamma - \eta  M_{2s}^2 \right] \| \bd^k \|^2-  \|\bx^k_{ T^c_{k}  }\|^2/(2\eta ) -\eta   \| \nabla_{T_{k-1}} f(\bx^{k})\|^2\\
&\leq& -2\rho \| \bd^k\|^2-\eta   \| \nabla_{T_{k-1}} f(\bx^{k})\|^2
 \end{eqnarray*}
 {\bf Case 2:} If $\bd^k=\bd^k_g$, then it follows from (\ref{gradient-Direction}) (namely, $\bd^{k}_{T_{k} }=-\nabla_{T_{k} } f(\bx^{k})$) that
 \begin{eqnarray*}
 2\langle \nabla f(\bx^{k}), \bd^{k} \rangle&=& 2\langle \nabla_{T_{k} } f(\bx^{k}),  \bd^{k}_{T_{k} }\rangle- 2\langle \nabla_{ T^c_{k}  } f(\bx^{k}),  \bx^{k}_{ T^c_{k}  }\rangle\\
&\overset{(\ref{xk-Tc-T1}, \ref{existence-alpha-facts-3})}{\leq}& -2 \| \bd^k_{T_{k}}\|^2  +\eta \|\nabla_{T_{k} } f(\bx^{k})\|^2-\|\bx^k_{ T^c_{k}  }\|^2/ \eta   -\eta   \| \nabla_{T_{k-1}} f(\bx^{k})\|^2\\
&=& -(2-\eta ) \| \bd^k_{T_{k}}\|^2- \|\bd^k_{ T^c_{k}  }\|^2/ \eta   -\eta   \| \nabla_{T_{k-1}} f(\bx^{k})\|^2\\
&\leq & -(2-\eta ) (\| \bd^k_{T_{k}}\|^2+ \|\bd^k_{ T^c_{k}  }\|^2)  -\eta   \| \nabla_{T_{k-1}} f(\bx^{k})\|^2\\
& \leq& -2\rho \| \bd^k\|^2-\eta   \| \nabla_{T_{k-1}} f(\bx^{k})\|^2,
 \end{eqnarray*}
where the second inequality used the fact $\eta(2-\eta) \le 1$. This finishes the proof.
\end{proof}

\subsection{Proof of Lemma~\ref{Lemma-alphak}}
\begin{proof} If $0<\alpha\leq\overline{\alpha}$ and $0<\gamma\leq\min \{1,2M_{2s}\}$, we have
$$
\alpha \leq \frac{ 1-2\sigma }{M_{2s}/ \gamma  -\sigma} \leq \frac{ 1-2\sigma }{M_{2s}  -\sigma}.$$
Since $f$ is  $M_{2s}$-restricted strongly smooth, we have
 \begin{eqnarray}\label{2f-2f}
 2f(\bx^{k}(\alpha))- 2f(\bx^{k})
 &\overset{(\ref{Eq-Lipschitz})}{\leq}&
 2\langle \nabla f(\bx^{k}), \bx^{k}(\alpha)-  \bx^{k}\rangle+  M_{2s}  \|\bx^{k}(\alpha)-  \bx^{k}\|^2\nonumber\\
  &=&2\langle \nabla f(\bx^{k}), \bx^{k}(\alpha)-  \bx^{k}\rangle+  M_{2s}  \|\bx^{k}(\alpha)-  \bx^{k}\|^2\nonumber\\
  & & -\; 2\alpha\sigma\langle \nabla f(\bx^{k}), \bd^{k}\rangle+2\alpha\sigma\langle \nabla f(\bx^{k}), \bd^{k}\rangle\nonumber\\
 &\overset{(\ref{xk-alpha})}{=}&  \alpha(1-\sigma)2\langle \nabla_{T_{k} } f(\bx^{k}),  \bd^{k}_{T_{k} }\rangle-(1-\alpha\sigma)2\langle \nabla_{ T^c_{k} } f(\bx^{k}),  \bx^{k}_{ T^c_{k} }\rangle\nonumber\\
 & & +\;  M_{2s} \left[\alpha^2\|\bd^{k}_{T_{k} }\|^2+ \|\bx^{k}_{ T^c_{k} }\|^2\right]+2\alpha\sigma\langle \nabla f(\bx^{k}), \bd^{k}\rangle\nonumber\\
  &\overset{(\ref{xk-Tc-T1})}{=} &  \Delta + 2\alpha\sigma\langle \nabla f(\bx^{k}), \bd^{k}\rangle , \nonumber
%
 \end{eqnarray}
where
\[
   \Delta := \alpha(1-\sigma)2\langle \nabla_{T_{k} } f(\bx^{k}),  \bd^{k}_{T_{k} }\rangle-(1-\alpha\sigma)2\langle \nabla_{ J_{k} } f(\bx^{k}),  \bx^{k}_{J_{k}}\rangle
   +\;  M_{2s} \left[\alpha^2\|\bd^{k}_{T_{k} }\|^2+ \|\bx^{k}_{ T^c_{k} }\|^2\right]
\]
To conclude the conclusion,  we only need to show $ \Delta \leq0$. We  prove it by two cases.\\

 \noindent {\bf Case 1:} If $\bd^k:=\bd^k_N$, then combining (\ref{fd-TT}) and (\ref{xk-T1-T1}) yields that
\begin{eqnarray*}
 \Delta &\leq&  \alpha(1-\sigma)2\langle \nabla_{T_{k} } f(\bx^{k}),  \bd^{k}_{T_{k} }\rangle-(1-\alpha\sigma)2\langle \nabla_{ J_{k} } f(\bx^{k}),  \bx^{k}_{J_{k}}\rangle +  M_{2s} \left[\alpha^2\|\bd^{k}\|^2+ \|\bx^{k}_{ T^c_{k} }\|^2\right]\\
&\leq& c_1\| \bd^k \|^2  +  c_2  \|\bx^k_{ T^c_{k} }\|^2 
 -(1-\alpha \sigma)\eta  \| \nabla_{T_{k-1}} f(\bx^k) \|^2,
 \end{eqnarray*}
 where
\begin{eqnarray*}
c_1&:=& -\alpha(1-\sigma)2\gamma +(1-\alpha\sigma)\eta  M_{2s}^2+ M_{2s}\alpha^2, \\
&\leq& -\alpha(1-\sigma)2\gamma +(1-\alpha\sigma)\gamma \alpha+ M_{2s}\alpha^2 \hspace{0.6cm} {\rm because~of}~\alpha \leq 1 ,\sigma \leq \frac{ 1}{2},\eta \leq \frac{ \alpha \gamma  }{ M_{2s}^2}  \\
&=& \alpha \left[( M_{2s}-\sigma\delta )\alpha-(1-2\sigma)\gamma \right]\leq0,\hspace{1.05cm} {\rm because~of}~\sigma\gamma\leq M_{2s}, \alpha\leq \frac{ 1-2\sigma }{ M_{2s}/\gamma -\sigma }\\
c_2&:=& \alpha(1-\sigma)/(2\eta )-(1-\alpha\sigma)/\eta  + M_{2s}\\
&\leq&(1-\alpha\sigma)/(2\eta )-(1-\alpha\sigma)/\eta  + M_{2s}\hspace{0.8cm}{\rm because~of}~\alpha \leq 1\\
&\leq&-(1-\alpha\sigma)/(2\eta ) + M_{2s}\leq 0,  \hspace{2cm} 
      {\rm because~of}~\alpha \leq 1 ,\sigma \leq \frac{ 1}{2}, \eta \leq\frac{1}{4M_{2s}}.
\end{eqnarray*}

\noindent
 {\bf Case 2:} If $\bd^k:=\bd^k_g$, then combining (\ref{gradient-Direction})  that $\bd^{k}_{T_{k} }=-\nabla_{T_{k} } f(\bx^{k})$ and  (\ref{existence-alpha-facts-3}) suffices to
 \begin{eqnarray*}
  \Delta &\leq&  c_3\| \bd^k_{T_{k}}\|^2  + c_4\|\bx^k_{ T^c_{k} }\|^2
                -(1 - \alpha \sigma ) \eta   \| \nabla_{T_{k-1}} f(\bx^k) \|^2,
 \end{eqnarray*}
  where
\begin{eqnarray*}c_3&:=& -2\alpha(1-\sigma)+(1-\alpha\sigma)\eta  + M_{2s}\alpha^2  \\
&\leq&\alpha\left[( M_{2s}-\sigma )\alpha-(1-2\sigma) \right] \hspace{1.6cm} {\rm because~of}~\alpha \leq 1 ,\sigma \leq \frac{ 1}{2},\eta \leq \alpha\\
&\leq&0. \hspace{5.7cm} {\rm because~of}~ \alpha\leq \frac{ 1-2\sigma }{ M_{2s}  -\sigma }\\
c_4&:=& -(1-\alpha\sigma)/\eta + M_{2s}\\
&\leq&-1/(2\eta )+M_{2s}\leq0,\hspace{2.7cm} {\rm because~of}~\alpha \leq 1 ,\sigma \leq \frac{ 1}{2},\eta \leq\frac{1}{4M_{2s}} \end{eqnarray*}
 which finishes proving the first claim.  If $\eta \in(0, \overline \eta )$ where $\overline \eta $ is defined as (\ref{alpha-eta}), then for any $\beta \overline \alpha \leq \alpha\leq \overline \alpha$, we have
$$0<\eta <\min\left\{ \frac{ \overline{\alpha} \gamma \beta}{M_{2s}^2}~~\overline{\alpha}\beta,~\frac{1}{4M_{2s}}  \right\}\leq\min\left\{\frac{ \alpha \gamma  }{M_{2s}^2},~\alpha,~ \frac{1}{4M_{2s}}\right\}.$$
 This together with  (\ref{alpha-decreasing-property}), namely,
$f(\bx^{k}(\alpha))- f(\bx^{k})  \leq  \sigma \alpha \langle \nabla  f(\bx^{k}),  \bd^{k} \rangle,$
and the Armijo-type step size rule  means  that  $\{\alpha_k\}$ is bounded from below by a positive constant, that is,
\begin{equation} \label{Positive-Lower-Bound}
\inf_{k\geq 0}\{\alpha_k\} \geq \beta \overline \alpha > 0.
\end{equation}
which finishes the whole proof. \end{proof}
\subsection{Proof of Lemma~\ref{Lemma-Converging-Quantities}}
\begin{proof} Lemma \ref{Lemma-alphak} shows the existence of $\alpha_k$, then  (\ref{Armijo}) in {\tt NHTP} (namely, (\ref{alpha-decreasing-property})) provides
\begin{eqnarray}
f(\bx^{k+1})-f(\bx^k)\leq  \sigma \alpha_k 
   \langle \nabla  f(\bx^{k}),  \bd^{k} \rangle 
 &\overset{(\ref{Descent-dk})}{\leq}& - \sigma \alpha_k 
  \left[ \rho\| \bd^k\|^2+\frac{\eta}{2} \| \nabla_{T_{k-1}} f(\bx^{k})\|^2  \right] \nonumber\\
\label{ffdfk}&\overset{(\ref{Positive-Lower-Bound})}{\leq}& -\sigma \overline \alpha \beta
  \left[ \rho\| \bd^k\|^2+\frac{\eta}{2}  \| \nabla_{T_{k-1}} f(\bx^{k})\|^2 \right].
\end{eqnarray}
Thus $f(\bx^{k+1})<f(\bx^k)$ if $\bx^{k+1}\neq \bx^k$. Then it follows from above inequality that
\begin{eqnarray*}
{\sigma \overline \alpha \beta }
\left[ \rho\sum^{\infty}_{k=0}\| \bd^k\|^2+ \frac{\eta}{2}  \sum^{\infty}_{k=0}\| \nabla_{T_{k-1}} f(\bx^{k})\|^2\right]
&\leq&
 \sum^{\infty}_{k=0}\left[ f(\bx^k)-f(\bx^{k+1})\right]\\
 &<& \left[f(\bx^0)-\lim _{k\rightarrow +\infty}f(\bx^k)\right]<+\infty,
\end{eqnarray*}
where the last inequality is due to $f$ being bounded from below.
Hence $${\lim}_{k\rightarrow \infty}\| \bd^k\|={\lim}_{k\rightarrow \infty}\| \nabla_{T_{k-1}} f(\bx^{k})\|=0$$
which
suffices to $\lim_{k\rightarrow \infty}\|\bx^{k+1}-\bx^k\|=0$ because of
\begin{eqnarray}\label{xxdk} \|\bx^{k+1}-\bx^k\|^2\overset{(\ref{xk-alpha})}{=}  \alpha_k^2 \| \bd^k_{T_{k}}\|^2+\| \bx^k_{ T^c_{k}}\|^2\leq \| \bd^k_{T_{k}}\|^2+\| \bd^k_{ T^c_{k}}\|^2=\| \bd^k\|^2.\end{eqnarray}
If $\bd^k=\bd^k_N$, then it follows from (\ref{nabla-T-f})
that $
\|\nabla_{ T_{k}}  f (\bx^k)\|  \leq   M_{2s} \|\bd^k\|.$
If $\bd^k=\bd^k_g$, it follows from   (\ref{gradient-Direction})   that $\|\nabla_{ T_{k}}  f (\bx^k)\|=\| \bd^k_{T_{k}}\|$. Those suffice  to ${\lim}_{k\rightarrow \infty} \|\nabla_{ T_{k}}  f (\bx^k)\|= 0$. Finally,  (\ref{station}) allows us to derive
$$\|F_{\eta }(\bx^{k};T_{k})\|^2= \|\nabla_{ T_{k}}  f (\bx^k)\|^2 + \| \bx^k_{ T^c_{k}}\|^2\leq (M_{2s}^2+1)\|\bd^k\|^2,$$
which is also able to claim  ${\lim}_{k\rightarrow \infty}\|F_{\eta }(\bx^{k};T_{k})\|=0.$
\end{proof}
\subsection{Proof of Theorem~\ref{Thm-Global-Convergence}}
\begin{proof} (i) We prove in Lemma \ref{Lemma-Converging-Quantities} (iv) that
 \begin{eqnarray}\label{gTf00}{\lim}_{k  \rightarrow \infty} \nabla_{ T_{k}}  f (\bx^{k+1}) =0.\end{eqnarray}
Let  $\{\bx^{k_\ell}\}$ be the convergent subsequence of $\{\bx^{k}\}$  that converges to $\bx^*$. 
Since there are only finitely many choices for $T_k$, 
(re-subsequencing if necessary)  we may without loss of any generality
assume that the sequence of the index sets $\{ \{T_{k_\ell-1}\} \}$ shares a same index set, denoted
as $T_{\infty}$. That is 
\begin{eqnarray}\label{TTTT}
   T_{k_\ell-1}=T_{k_{\ell+1}-1}=\cdots =T_\infty.
\end{eqnarray}
Since $\bx^{k_\ell}\rightarrow \bx^*$, $\supp(\bx^{k_\ell})\subseteq T_{k_\ell-1}= T_\infty$, we must have
  \begin{eqnarray*}
 T_\infty\left\{\begin{array}{ll}
=\Gamma_*:=\supp(\bx^*),& {\rm if } ~\|\bx^*\|_0=s,\\
\supset\Gamma_* ,& {\rm if } ~\|\bx^*\|_0<s.
\end{array}\right.
 \end{eqnarray*}
 which implies
 \begin{eqnarray}\label{gTf0}
  \nabla_{T_\infty} f (\bx ^{*})=\underset{{k_\ell}\rightarrow \infty}{\lim} \nabla_{ T_{k_\ell-1}}  f (\bx^{k_\ell})\overset{(\ref{gTf00})}{= }0.\end{eqnarray}
In addition, the definition (\ref{Tu}) of $\T(\bx^{k},\eta)$ means
 \begin{eqnarray}\label{Tu-k}
 |\bx^k_{i}-\eta \nabla_{i}  f (\bx^k)|\geq |\bx^k_{j}-\eta \nabla_{j}  f (\bx^k)|, ~~~~\forall~i\in T_{k},\ \forall ~j\in  T^c_{k}
 \end{eqnarray}
Again by Lemma \ref{Lemma-Converging-Quantities} (ii) that $ \lim_{k\rightarrow \infty} \|\bx^{k+1}-\bx^k\| =0,$ we obtain $ \lim_{k_\ell\rightarrow \infty} \bx^{k_\ell-1}=\bx^*$ due to $ \lim_{k_\ell\rightarrow \infty} \bx^{k_\ell}=\bx^*$. Now we have the following chain of inequalities
for any $ i\in T_{k_\ell-1}\overset{(\ref{TTTT})}{=}T_\infty,~j\in  T^c_{k_\ell-1}\overset{(\ref{TTTT})}{=}  T^c_\infty$
 \begin{eqnarray*}
 |\bx^*_{i}|&\overset{(\ref{gTf0})}{=}&|\bx^*_{i}-\eta \nabla_{i}  f (\bx^*)|~=~\lim_{{k_\ell}\rightarrow\infty} |\bx^{k_\ell-1}_{i}-\eta \nabla_{i}  f (\bx^{k_\ell-1})|\\
 &\overset{(\ref{Tu-k})}{\geq}&\lim_{{k_\ell}\rightarrow\infty} |\bx^{k_\ell-1}_{j}-\eta \nabla_{j}  f (\bx^{k_\ell-1})| ~=~ |\bx^*_{j}-\eta \nabla_{j}  f (\bx^*)|~=~\eta |\nabla_{j}  f (\bx^*)|,
 \end{eqnarray*}
 which leads to
 \begin{eqnarray}\label{xftf}
 \bx^*_{(s)} =\min_{i\in T_{\infty}}|\bx^*_{i}| \geq \eta |\nabla_{j}  f (\bx^*)|,~~ \forall ~~j\in   T^c_\infty, 
 \end{eqnarray}
If $\|\bx^*\|_0=s$, then $T_\infty=\Gamma_*$. Consequently, $x^*_{(s)} \geq \eta |\nabla_{j}  f (\bx^*)|, \forall~j\in \Gamma^c_*$.
If $\|\bx^*\|_0<s$, then $x^*_{(s)}=0$ and $\nabla  f (\bx^*)=0$ from (\ref{xftf}) and (\ref{gTf0}).
Those together with (\ref{agradient}) enable us to show  that  $\bx^*$ is an   $\eta$-stationary point.\\

If $f(\bx)$ is   convex,  letting $\Gamma_*:=\supp(\bx^*)$, then  
\begin{eqnarray*}
 f(\bx)&\geq &f(\bx^*)+\langle \nabla f(\bx^*), \bx-\bx^*\rangle \\
&>&f(\bx^*)+\sum_{i\in \Gamma_*} \nabla_i f(\bx^*)(x_i-x_i^*)+\sum_{i\notin \Gamma_*} \nabla_i f(\bx^*)(x_i-x_i^*)\\
&\overset{(\ref{agradient})}{=} &f(\bx^*)+ \sum_{i\notin\Gamma_*} \nabla_i f(\bx^*)x_i\geq f(\bx^*)- \sum_{i\notin \Gamma_*} |\nabla_i f(\bx^*)||x_i| \\ 
&\overset{(\ref{agradient})}{\geq} &f(\bx^*)- (\bx^*_{(s)}/\eta)\sum_{i\notin \Gamma_*}|x_i| \\
&=&f(\bx^*)- (\bx^*_{(s)}/\eta)\|x_{\Gamma_*^c}\|_1.
\end{eqnarray*}

(ii) The whole sequence converges because of  \citep[Lemma 4.10,][]{more1983computing} and  $\lim_{k\rightarrow \infty} \|\bx^{k+1}-\bx^k\| =0$ from Lemma \ref{Lemma-Converging-Quantities} (ii). If $\bx^*=0$, then the conclusion holds clearly due to $\support(\bx^*)=\emptyset$. We consider $\bx^*\neq0$. Since $\lim_{k\rightarrow\infty}\bx^k=\bx^*$, the for sufficiently large $k$ we must have $$\|\bx^k-\bx^*\|<\min_{i\in\support(\bx^*)} |x^*_i|=:t^*.$$
If $\supp(\bx^*)\nsubseteq\supp(\bx^k)$, then there is an $i_0\in\supp(\bx^*)\setminus\supp(\bx^k)$ such that
$$t^*>\|\bx^k-\bx^*\|\geq | x_{i_0}^k-x_{i_0}^* |=|  x_{i_0}^* |\geq t^*,$$
which is a contradiction. Therefore, $\supp(\bx^*)\subseteq\supp(\bx^k)$. By the updating rule (\ref{xk-alpha}), we have $\supp(\bx^k)\subseteq T_{k-1}$, where $|T_{k-1}|=s$ by (\ref{Tu}).  Therefore, if $\|\bx^{*}\|_0=s$ then $\supp(\bx^*)
 \equiv\supp(\bx^k)\equiv \supp(\bx^{k+1})\equiv T_{k}$. If $\|\bx^{*}\|_0<s$ then $\supp(\bx^*)
 \subseteq \supp(\bx^k),  \supp(\bx^*)
 \subseteq  \supp(\bx^{k+1})\subseteq T_{k}$.
 The whole proof is finished.\end{proof}

\subsection{Proof of Theorem~\ref{Thm-Qudratic-Convergence}}

\begin{proof} 
(i) We have proved in Theorem \ref{Thm-Global-Convergence} (i) that any limit $\bx^*$ of $ \{\bx^k\}$ is an $\eta$-stationary point. If $f(\bx)$ is  $m_{2s}$-restricted strongly convex in a neighborhood of $\bx^*$, then we can conclude that $\bx^*$ is   a strictly local minimizer of (\ref{SCO}). In fact,
\begin{eqnarray*}
 f(\bx)&\geq &f(\bx^*)+\langle \nabla f(\bx^*), \bx-\bx^*\rangle+(m_{2s}/2)\|\bx-\bx^*\|^2\\
&>&f(\bx^*)+\sum_{i\in \supp(\bx^*)} \nabla_i f(\bx^*)(x_i-x_i^*)+\sum_{i\notin \supp(\bx^*)} \nabla_i f(\bx^*)(x_i-x_i^*)\\
&\overset{(\ref{agradient})}{=} &f(\bx^*)+ \sum_{i\notin \supp(\bx^*)} \nabla_i f(\bx^*)x_i\\
&\overset{(\ref{agradient})}{=} &f(\bx^*)+ \left\{
\begin{array}{ll}
 \sum_{i\notin \supp(\bx^*)} \nabla_i f(\bx^*) \times 0,& {\rm}~\|\bx^*\|_0=s\\
 \sum_{i\notin \supp(\bx^*)} 0 \times x_i,& {\rm}~\|\bx^*\|_0<s.
\end{array}\right. \\
&=& f(\bx^*)
\end{eqnarray*}
for any $s$-sparse vector $\bx$, where the first inequality is from the $m_{2s}$-restricted strongly convexity.   This also shows $\bx^*$ is isolated and thus the whole sequence tends to $\bx^*$ by Theorem \ref{Thm-Global-Convergence} (ii).

(ii) The fact that $f(\bx)$ is  $m_{2s}$-restricted strongly convex in a neighborhood of $\bx^*$ 
and $\lim_{k\rightarrow\infty}\bx^k=\bx^*$ implies that 
 $f(\bx)$ is also  $m_{2s}$-restricted strongly convex in a neighborhood of $\bx^k$ for sufficiently large $k$. 
 By invoking Lemma~\ref{Lemma-Newton-Direction}, we see that the Newton direction $\bd^k_N$
 always satisfies the condition (\ref{Descent-Inequality-Newton}) and hence is accepted
 as the search direction when $k$ is sufficiently large.

(iii) By $\supp(\bx^*)\subseteq T_{k}$ for sufficiently large $k$ from \ref{Thm-Global-Convergence} (ii) and  $\bx^*$ is an $\eta$-stationary point, it follows from  Theorem (\ref{agradient}) that
\begin{eqnarray}\label{F0-s-1}
\bx^*_{ T^c_k }=0 \quad \mbox{and} \quad 
\left\{
\begin{array}{rcc}
 \nabla_{T_{k}}  f (\bx^*) = \nabla_{\supp(\bx^*)}  f (\bx^*) & =0 & {\rm if}~   \|\bx^*\|_0=s,\\
    \nabla f(\bx^*) & =0 & {\rm if}~   \|\bx^*\|_0<s.
\end{array}
\right. \end{eqnarray}
For any $0\leq t \leq1$, denote  $\bx(t):=\bx^* + t (\bx^k-\bx^*)$. Clearly, as $\bx^k$,  $\bx(t)$ is also in the neighbour of $\bx^*$ and  {{ $\supp(\bx^*)\subseteq  \supp(\bx(t)) \subseteq \supp(\bx^k)$ due to $\supp(\bx^*)\subseteq\supp(\bx^k)$.}}
So $f$ being locally restricted Hessian Lipschitz continuous at $\bx^*$ with the Lipschitz constant $L_f$ {{and $T_{k}\supseteq\supp(\bx^*)$ 
give  rise to
\begin{eqnarray}
\label{facts-0-4}  \|\nabla^{2}_{{T_{k}: }} f (\bx^k)-\nabla_{{T_{k} :} }^{2} f (\bx(t))\|
\leq  L_f\|\bx^k- \bx(t)\| = (1-t)L_f\|\bx^k- \bx^*\|.
\end{eqnarray}}}
 Moreover, by the Taylor expansion, we have
\begin{eqnarray}\label{facts-0-3}
\nabla  f (\bx ^k)-\nabla  f (\bx ^*)=\int_0^1\nabla^{2} f ( \bx(t))(\bx ^k-\bx ^*)dt.
\end{eqnarray}
We also have the following chain of inequalities
 \begin{eqnarray}
\|\bx^{k+1}-\bx^*\|&=&\left[\| \bx^{k+1}_{ T_{k} }-\bx^*_{ T_{k} }\|^2+\| \bx^{k+1}_{  T^c_{k} }-\bx^*_{  T^c_{k} }\|^2\right]^{1/2}\nonumber\\
&{=}&\| \bx^{k+1}_{ T_{k} }-\bx^*_{ T_{k} }\|\overset{(\ref{xk-alpha})}{=}\| \bx^{k}_{ T_{k} }-\bx^*_{ T_{k} }+\alpha_k \bd^{k}_{ T_{k} }\|\nonumber\\
\label{facts-0-5-1}&\leq&(1-\alpha_k)\| \bx^{k}_{ T_{k} }-\bx^*_{ T_{k} }\|+\alpha_k \|\bx^{k}_{ T_{k} }-\bx^*_{ T_{k} }+\bd^{k}_{ T_{k} }\| \\
\label{facts-0-5}&\overset{(\ref{Positive-Lower-Bound})}{\leq}&(1-\overline{\alpha}\beta )\| \bx^{k} -\bx^* \|+\overline{\alpha} \|\bx^{k}_{ T_{k} }-\bx^*_{ T_{k} }+\bd^{k}_{ T_{k} }\|,
\end{eqnarray}
where the second equality used the fact (\ref{xk-alpha}) and $\support(\bx^{k+1}) \subseteq T_k$.
Since $\bd^k=\bd^k_N$, we have
\allowdisplaybreaks  \begin{eqnarray}
&&\|\bx^{k}_{ T_{k} }-\bx^*_{ T_{k} }+\bd^{k}_{ T_{k} }\|\\
&\overset{(\ref{Newton-Direction})}{=}&\left\|H_k^{-1}\left(\nabla^2_{T_{k},  T^c_{k}  } 
  f (\bx ^k)\bx^k_{ T^c_{k} } -\nabla_{T_{k} } f (\bx ^k)\right)+\bx^{k}_{ T_{k} }-\bx^*_{ T_{k} }\right\|\nonumber\\
&=& \left\| 
    H_k^{-1} \left( \nabla^2_{T_{k},  T^c_{k}  } 
    f (\bx ^k)\bx^k_{ T^c_{k} } -\nabla_{T_{k} } f (\bx ^k)
    + \nabla^2_{T_k} f(\bx^k) \bx^k_{T_k} - \nabla^2_{T_k} f(\bx^k) \bx^*_{T_k} 
    \right) 
   \right\| \nonumber \\
&\leq&\frac{1}{m_{2s}}\left\|\nabla^2_{T_{k}  : } f (\bx ^k)\bx^k -\nabla_{T_{k} } f (\bx ^k)-\nabla_{T_{k}}^{2} f (\bx^k)\bx^*_{ T_{k} } \right\|\nonumber\\
&\overset{(\ref{F0-s-1})}{=}&\frac{1}{m_{2s}}\left\|\nabla^2_{T_{k}  : } f (\bx ^k)\bx^k -\nabla_{T_{k} } f (\bx ^k)-\nabla_{{T_{k}}  : }^{2} f (\bx^k)\bx^* +\nabla_{T_{k} } f (\bx ^*) \right\|\nonumber\\
&\overset{(\ref{facts-0-3})}{=}&\frac{1}{m_{2s}}\left\|\nabla_{{T_{k} }  : }^{2} f (\bx^k)(\bx^{k}-\bx^*)-\int_0^1\nabla_{{T_{k}}{ : }}^{2} f ( \bx(t))(\bx ^k-\bx^* )dt\right\|\nonumber\\
&=&\frac{1}{m_{2s}}\left\|\int_0^1\left[\nabla_{{T_{k} } : }^{2} f (\bx^k)-\nabla_{{T_{k} }{ : }}^{2} f ( \bx(t))\right](\bx ^k-\bx^*)dt\right\|\nonumber\\
&\leq&\frac{1}{m_{2s}}\int_0^1\left\|\nabla_{{T_{k}: }   }^{2} f (\bx^k)-\nabla_{{T_{k} }{ :}}^{2} f ( \bx(t))\right\| \|\bx ^k-\bx ^* \|dt\nonumber\\
&\overset{(\ref{facts-0-4})}{\leq}&\frac{ L_f}{m_{2s}}\|\bx ^k-\bx ^*\|^2 \int_0^1(1-t)dt  \nonumber\\
\label{facts-7}&=&L_f/(2m_{2s})\|\bx ^k-\bx ^*\|^2.
\end{eqnarray}
Now, we have obtained 
(fact 1) $\lim_{k\rightarrow\infty}\bx^k=\bx^*$,
(fact 2) $ \langle \nabla f(\bx^{k}), \bd^{k}\rangle  \leq -  \rho\| \bd^{k}\|^2 $ from Lemma \ref{Lemma-dk} 
and (fact 3)
  \begin{eqnarray*}\lim_{k\rightarrow\infty}\frac{\|\bx^k+\bd^k-\bx^*\|}{\|\bx^k-\bx^*\|}=\lim_{k\rightarrow\infty}\frac{\|\bx^k_{T_{k} }+\bd^k_{T_{k} }-\bx^*_{T_{k} }\|}{\|\bx^k-\bx^*\|}
\overset{( \ref{facts-7})}{\leq} \lim_{k\rightarrow\infty}\frac{L_f\|\bx ^k-\bx ^*\|^2}{2m_{2s}\|\bx^k-\bx^*\|}=0,
\end{eqnarray*}
where the first equality is because of $\bd^k_{ T^c_{k}}=-\bx^k_{ T^c_{k}}$ and (\ref{F0-s-1}). 
These three facts are exactly the same assumptions used in 
\citep[Theorem 3.3,][]{facchinei1995minimization}, which establishes that
eventually the step size $\alpha_k$ in the Armijo rule has to be $1$, namely $\alpha_k\equiv 1$.  
Therefore, for sufficiently large $k$, it follows from (\ref{facts-0-5-1}) that
  \begin{eqnarray}
\|\bx^{k+1}-\bx^*\|&\leq&\alpha_k \|\bx^{k}_{ T_{k} }-\bx^*_{ T_{k} }+\bd^{k}_{ T_{k} }\|+(1-\alpha_k)\| \bx^{k}_{ T_{k} }-\bx^*_{ T_{k} }\|\nonumber\\
&=& \|\bx^{k}_{ T_{k} }-\bx^*_{ T_{k} }+\bd^{k}_{ T_{k} }\| \nonumber\\
\label{FxT-xx} &\overset{( \ref{facts-7})}{\leq}& (L_f/ 2m_{2s} )\|\bx ^k-\bx ^*\|^2.
\end{eqnarray}
That is, we have proved that the sequence has a quadratic convergence rate. Finally, for sufficiently large $k$, it follows
\begin{eqnarray}\| F_{\eta}(\bx^{k+1};T_{k+1})\|^2
&\overset{(\ref{station})}{=}&
\|\nabla_{ T_{k+1}}f (\bx^{k+1})\|^2+
\|\bx^{k+1}_{ T^c_{k+1}}\|^2 \nonumber\\
 &\overset{(\ref{F0-s-1})}{=}& \|\nabla_{ T_{k+1}}f (\bx^{k+1})-\nabla_{ T_{k+1}}f (\bx^*)\|^2+
\|\bx^{k+1}_{ T^c_{k+1}} -\bx^{*}_{ T^c_{k+1}}
\|^2 \nonumber\\
&\overset{(\ref{Eq-Lipschitz})}{\leq}& (M_{2s}^2+1) \| \bx^{k+1}-\bx^* \|^2\nonumber\\
\label{Fk-1}&\overset{(\ref{FxT-xx})}{\leq}& (M_{2s}^2+1)(L_f/ 2m_{2s} )^2 \| \bx^{k}-\bx^* \|^4.
 \end{eqnarray}
 {{Since $f(\bx)$ is  $m_{2s}$-restricted strongly convex in a neighborhood of $\bx^*$, 
\begin{eqnarray*}
\nabla^2_{{T}} f(\bx^*) \succeq m_{2s} I_{2s}  \hspace{10mm} {\rm for~any~}~ T\supseteq \supp(\bx^*), |T| \leq 2s. 
\end{eqnarray*}
This together with $\supp(\bx^*)\subseteq T_k$ from ii) in Theorem \ref{Thm-Global-Convergence} indicates
\begin{eqnarray*}
\sigma_{\min}(F'_{\eta} (\bx^*; T_k)) = \sigma_{\min}\left(\left[
 \begin{array}{cc}
  \nabla^2_{T_k} f(\bx^*) \ & \ \nabla^2_{T_k, T_k^c} f(\bx^*)  \\ [0.6ex]
  0 & I_{n-s}
 \end{array}
\right]\right) \geq \min\{m_{2s},1\},
\end{eqnarray*}
where $\sigma_{\min}(A)$ denotes the smallest singular value of $A$. Then we have following Taylor expansion for a fixed $T_k$,
\begin{eqnarray*}
\|F_{\eta} (\bx^k; T_k)\|&\geq& \|F_{\eta} (\bx^*; T_k)+F'_{\eta} (\bx^*; T_k)(\bx^k-\bx^*)\| - o(\|\bx^k-\bx^*\|)\\
&= &\|F'_{\eta} (\bx^*; T_k)(\bx^k-\bx^*) \|- o(\|\bx^k-\bx^*\|),\\
&\geq &(1/\sqrt{2})\|F'_{\eta} (\bx^*; T_k)(\bx^k-\bx^*) \|,\\
&\geq &(\min\{m_{2s},1\}/\sqrt{2})\|\bx^k-\bx^*\|,
\end{eqnarray*}
where the first equation holds due to $F_{\eta} (\bx^*; T_k)=0$ by (\ref{F0-s-1}). Finally, we have
\begin{eqnarray*}
\|F_{\eta} (\bx^k; T_k)\|^2\geq 
(\min\{m^2_{2s},1\}/2)\|\bx^k-\bx^*\|^2
\overset{(\ref{Fk-1})}{\geq} \frac{\min\{m^3_{2s},m_{2s}\}}{L_f \sqrt{M_{2s}^2+1}} \| F_{\eta}(\bx^{k+1};T_{k+1})\|. 
\end{eqnarray*}
 This completes the whole proof.}} \end{proof}
\vskip 0.2in
\bibliographystyle{apalike}
\bibliography{refsco}

\end{document}